\documentclass{article}
\usepackage[utf8]{inputenc}
\usepackage{authblk}
\usepackage{setspace}
\usepackage{titlesec}
\usepackage[left=0.75in, right=0.75in, top=1in, bottom=1in]{geometry}
\usepackage{graphicx}
\graphicspath{ {./figures/} }
\usepackage{subcaption}
\usepackage{amsmath}
\usepackage{lineno}
\usepackage[raggedrightboxes]{ragged2e}
\usepackage{array}
\newcolumntype{L}[1]{>{\raggedright\let\newline\\\arraybackslash\hspace{0pt}}m{#1}}
\newcolumntype{C}[1]{>{\centering\let\newline\\\arraybackslash\hspace{0pt}}m{#1}}
\newcolumntype{R}[1]{>{\raggedleft\let\newline\\\arraybackslash\hspace{0pt}}m{#1}}
\usepackage{amsfonts}
\usepackage{multirow}
\usepackage{graphicx}
\usepackage{moresize}
\usepackage{rotating}
\usepackage{algorithm}
\usepackage{algpseudocode}
\usepackage{caption}
\usepackage{titlesec}
\titleformat{\paragraph}[block]
  {\normalfont\normalsize\bfseries}{\theparagraph}{1em}{}
\titlespacing*{\paragraph}{0pt}{1ex plus .2ex minus .2ex}{0.5ex}

\usepackage[colorlinks=true]{hyperref}
\usepackage[inkscapelatex=false]{svg}

\hypersetup{
  citecolor=blue,
  linkcolor=blue,
  urlcolor=blue}
\newtheorem{theorem}{Theorem}[section]
\usepackage[style=nejm, 
citestyle=numeric-comp,
sorting=none]{biblatex}
\usepackage{lipsum}
\usepackage{fancyhdr}  
\usepackage{titlesec}  
\usepackage{multicol}  

\usepackage{marvosym}
\title{\Large A Reinforcement-learning-based Column Generation Algorithm for Integrated Operating Room Planning and Scheduling}

  \author[a,b]{\normalsize Mahdi Dolatkhah}
  \author[a,b]{Hossein Hashemi Doulabi}
  \author[b,c]{Walter Rei}
  \author[b,d]{Michel Gendreau}

  \affil[a]{\tiny Department of Mechanical, Industrial and Aerospace Engineering, Concordia University, Montréal H3G 1M8, Québec, Canada}
  \affil[b]{\tiny CIRRELT, Université de Montréal, Montréal H3C 3J7, Québec, Canada}
  \affil[c]{\tiny Département d’Analytique, Opérations et Technologies de l’Information, École des Sciences de la Gestion, Université du Québec à Montréal H2X 3X2, Québec, Canada}
  \affil[d]{\tiny Département de mathématiques et de génie industriel, Polytechnique Montréal, Montréal H3C 3A7, Québec, Canada}
  \affil[$ $ $ $]{Emails: mahdi.dolatkhah@mail.concordia.ca, hossein.hashemi@concordia.ca, rei.walter@uqam.ca, michel.gendreau@polymtl.ca}

\date{}

\makeatletter
\let\oldaffillist\AB@affillist
\renewcommand{\AB@affillist}{\begin{flushleft}\begin{spacing}{1.5}\oldaffillist\end{spacing}\end{flushleft}}
\makeatother

\begin{document}

\maketitle
\begin{sloppypar}


\vspace{-2cm}

\vspace{10 pt}

\begin{abstract}
 Operating room planning and scheduling are vital components of hospital management, contributing to improved efficiency, patient satisfaction, staff well-being, and overall quality of care delivery. In this paper, we propose a novel mixed integer programming model to formulate integrated operating room planning and scheduling problems, where several mandatory and elective surgeries are to be assigned and scheduled in operating rooms on different days. Aside from the standard working hours in each operating room, we also take into account the potential for performing surgeries in overtime periods. In addition, our approach also takes into account the availability of surgeons by considering their allowed surgical time on each day. We propose a column generation (CG) algorithm to solve large-scale instances. In order to enhance the CG, we integrate the Reinforcement Learning Algorithm and the Genetic Algorithm and develop a hybrid algorithm to generate initial columns for the CG algorithm. For our analysis, we employed two sets of test instances: one consisting of synthetic data and the other based on real-world cases from a local hospital in Naples, Italy. Computational experiments demonstrate that our proposed model and methodology yields an average optimality gap of 1.23\% for synthetic instances and 1.49\% on real-world scenarios, significantly outperforming previous solution methodologies in the literature. Additionally, we demonstrate that the developed CG algorithm provides a high-quality solution for large-scale instances where other models and methods fail to obtain even a feasible solution. To further evaluate robustness under uncertainty, we examined scenarios with $\pm20\%$ variability in surgery durations. The results indicate that incorporating a 120-minute buffer time minimizes the overall cost. Moreover, we investigated the impact of emergency surgeries by either introducing additional cases or escalating surgical priorities. For synthetic instances, the inclusion of emergency surgeries increased the total rescheduling cost by 4.13\%, whereas in the real-world Naples cases, priority escalation led to only a 0.11\% increase, highlighting the resilience of our proposed model in practical hospital settings.
 
\vspace{0.1cm}
\end{abstract}



\normalsize
\section{Introduction}\label{Introduction}
The recent statistics provided by the Organization for Economic Co-operation and Development (OECD) show that the average health expenditure of all countries increased by 51.34\% over a five-year period between 2018 and 2022. The main factors driving recent increases are the growing population of senior individuals and the emergence of new diseases such as COVID-19, which significantly have challenged health care systems \cite{ghandehari2022mixed, khalilpourazari2021using}. In the realm of healthcare systems, three predominant models exist: public, private, and a hybrid combination of both \cite{niles2023basics}. Many nations, including Canada, Australia, Japan, and the United Kingdom, achieve comprehensive coverage for their entire populations through universal public health systems. Conversely, countries like the Netherlands and Switzerland predominantly rely on 100\% primary private health coverage, subsidized by the government. In the intermediate spectrum, 37\% of the population of the United States benefits from Medicare, a form of universal health coverage for individuals aged 65 and older, while 53\% are covered by private insurance \cite{niles2023basics}. In Canada, except for dental care, optometry, and certain specialized medical services, the majority of medical care is delivered through publicly funded hospitals and clinics \cite{fierlbeck2023boundaries, haghi2023integrated}. 

Operating rooms (ORs) and their related activities play a significant role in private medical clinics and hospitals, contributing to approximately 70\% of total revenue \cite{andam2021operating}. Despite this substantial contribution, they also account for about 33\% of the overall budget \cite{dodaro2020solving}. Proficient management of expenditures and cost reduction are crucial for both public and private healthcare systems, facilitating new investments such as adoption of novel technologies and recruitment of highly skilled staff \cite{li2021health}. Beyond the economic perspective, effective planning and scheduling of operating rooms (ORPS) is critical to ensure that surgeries are conducted at the appropriate times, potentially saving patients' lives. Consequently, ORPS not only plays a pivotal role in minimizing expenses and increasing hospitals' financial performance, but also contributes to saving more lives, all while maintaining the satisfaction of both staff and patients by mitigating the existing prolonged surgical waiting lists \cite{dodaro2020solving}.

In this work, we study an integrated operating room planning and scheduling problems (IORPS) problem with a finite planning horizon and the overtime option being available for each OR on each day. We formulate the IORPS problem as a new mixed integer programming (MIP) model which includes significantly fewer variables compared to the model proposed by Marques et al. \cite{marques2012integer} which is based on a four-index variable structure. While the four-indexed structure provides detailed modeling flexibility, it results in a large-scale problem that is computationally demanding, especially for real-world instances. Despite this limitation, the four-index formulation remains popular and continues to be used by recent studies including Boccia et al. \cite{boccia2024integrated}, who adopted it to develop an ILP-based IORPS model for the block scheduling strategy. In this study, all model parameters are assumed to be deterministic. This assumption is justified by the fact that incorporating stochastic elements significantly increases the complexity of the problem. Consequently, many previous studies have opted for deterministic formulations to ensure tractability and computational efficiency \cite{hashemi2014constraint, hashemi2016constraint, roshanaei2021solving, boccia2024integrated}. Despite the inherent uncertainty in real-world operating room environments, such as variability in surgery durations and emergency arrivals, deterministic models remain a practical and widely accepted approach for capturing essential system dynamics while maintaining solvability. In fact, we are focusing on obtaining more practical planning and more detailed scheduling which will increase the chance of implementation in real-world applications. Also, the IORPS is already a difficult problem to deal with and adding more complexity should be done step by step. For instance, Coban et al. \cite{coban2023effect} developed the model only for a single OR due to the complexity that arose from studying the stochastic problem. Therefore, we will study the non-deterministic version in future studies. Decision variables in our model include the start time of surgeries, the number of opened ORs on each day, and the list of surgeries postponed to the next planning horizon. We considered surgery postponement since it is a frequent occurrence in healthcare systems and can be attributed to different factors, such as the limited availability of resources. Based on the new mixed integer programming model, we develop a CG algorithm to tackle  large-scale instances of the problem. In this research, we make the following contributions:
\begin{itemize}
    \item We propose a novel alternative reformulation for the state-of-the-art mathematical programming model introduced by Marques et al. \cite{marques2012integer} for the IORPS problem. For a given amount of time, our MIP model results in smaller optimality gaps. We also provide a theorem to prove that our model finds an equivalent optimal solution of the existing IORPS model proposed by Marques et al. \cite{marques2012integer}
    \item We formulate a CG algorithm for the purpose of solving the new model. The subproblems of our developed CG algorithm mimic the structure of our proposed MIP model.
    \item To enhance the developed CG algorithm, we integrate a reinforcement learning algorithm (RLA) and a genetic algorithm (GA) and introduce a hybrid intelligent algorithm with a new local search enhancer mechanism.
    \item We provide computational results for three sets of synthetic instances with a range of 40 to 120, and 200 to 480 surgeries to demonstrate the efficiency of our model and CG algorithm compared to an existing model and branch-and-check algorithm in the literature.
    \item In addition to the synthetic instances, we provided extensive numerical results for the real-world instances from a local hospital in Naples, Italy.
\end{itemize}

The remainder of this paper is organized as follows. In the next section, a brief literature review of operating room planning and scheduling is presented. In Section \ref{Problem_Definition}, we present the problem definition, main assumptions made, some notation, and the proposed new formulation for the IORPS. In Section \ref{Existing_mathematical_model} we present the existing IORPS mathematical model and some of its extensions. In Section \ref{Proposed_Model} we propose the mathematical model developed to solve the IORPS problem. In Section \ref{CG_Algorithm} we explains the design of the CG algorithm as well as the enhancements. In section \ref{Computational_Experiments} we present extensive computational results. Finally, in Section \ref{Conclusion}, we provide concluding marks and directions for future research.

\section{Literature Review}\label{Literature_Review}
The management challenges of operating rooms are generally categorized into three levels: strategic (long-term), tactical (medium-term), and operational (short-term) \cite{hashemi2016constraint}. The available times of the ORs are divided among surgeons or different departments at the strategic level, which is referred to as case-mix planning in the literature \cite{hashemi2016constraint}. As previously stated, the planning horizon at this stage is long-term, ranging from several months to a year or longer. At the strategic stage, decision-making typically entails determining the appropriate number of resources, including ORs and surgical specialties \cite{zhu2019operating}. At the tactical stage, the master surgery scheduling problem (MSSP) needs to be solved. In this stage, the number of available ORs, hours of availability, and prioritizing the surgeons or surgical departments are specified \cite{hashemi2016constraint}. The tactical stage problem involves creating a cyclic schedule, which is typically optimized on a monthly or quarterly basis \cite{zhu2019operating}. The operational level is referred to as elective case scheduling in the literature. Elective case scheduling has been classified into planning and scheduling problems. In the context of OR planning, the goal is to allocate surgeries to available operating rooms or available days. Studies show that by optimizing the planning phase and performing just one additional surgery per day, could save an average hospital up to 7 million dollars per year \cite{andam2021operating}. The main goal in OR scheduling is to find the optimal sequence of the surgeries or the assignment of surgeries to the specific time slots. In the IORPS, planning and scheduling decisions are made simultaneously, necessitating sophisticated mathematical models and solution methodologies to find the optimal solution \cite{hashemi2016constraint}.

\subsection{Efforts to formulate the OR planning and scheduling problem}\label{Efforts_to_formulate}
In recent years, many researchers in the literature have worked on the ORPS problem (\cite{ghandehari2022mixed}, \cite{dodaro2020solving}, \cite{marques2012integer}, \cite{boccia2024integrated}, \cite{hashemi2016constraint}, \cite{roshanaei2021solving, marques2014scheduling, bargetto2023branch, aringhieri2022combining, mashkani2021operating, wang2022adaptive, eshghali2023machine, almoghrabi2025surgery, li2024scheduling, gur2024stochastic, zhu2024surgical,andam2026operatingroomplanningpooling}). Various mathematical formulations and optimization methods have been developed to address the complexity of OR scheduling, which involves multiple interacting constraints such as surgeon availability, operating room capacity, patient prioritization, and emergency case management.
Marques et al. \cite{marques2012integer} proposed a formulation for the IORPS that relies on the use of a four-digit index binary decision variable definition, originally proposed by Roland et al. \cite{roland2006operating}. Specifically, the binary variable $x_{crtd}$ takes the value 1 if surgery $c$ starts at the beginning of time period $t$ on day $d$ in room $r$. Their proposed integer linear programming (ILP) formulation can still be considered as a state-of-the-art model for the IORPS \cite{bargetto2023branch}. However, their formulation poses significant challenges when applied in real-world settings, as it involves $O(|C| \times |R| \times |T| \times |D|)$ variables and constraints. It is obvious that the complexity and dimension of the problem increase dramatically by increasing the number of surgeries, time intervals, planning horizon, or the number of available ORs.
Hashemi Doulabi et al. \cite{hashemi2016constraint} developed a pattern-based solution approach to address the IORPS. They reformulated the IORPS model initially proposed by Marques et al. \cite{marques2012integer} using constraint programming. Subsequently, they introduced a Constraint Programming-Based Branch-and-Price-and-Cut Algorithm, incorporating three distinct branching rules to achieve the optimal solution for the IORPS problem.
Roshanaei, and Naderi \cite{roshanaei2021solving} extended the model proposed by Doulabi et al. \cite{hashemi2016constraint}. They created a new sequence-based MIP and constraint programming (CP) model for solving the IORPS. 
Bargetto et al. \cite{bargetto2023branch} studied the IORPS in the context of explicitly incorporating the available time of nurses and surgeons in various surgical groups as constraints in the developed model.
Assignment of surgeons to surgeries was recently investigated by Wang et al. \cite{wang2020two} where the goal was to assign two or more surgeons with different skill levels to one surgery. Also, as described by Vijayakumar et al. \cite{vijayakumar2013dual}, the surgeon-surgery assignment approach is usually considered in training hospitals where the surgeries are assigned to different assistant surgeons as part of the training program.
Kamran et al. \cite{kamran2019adaptive} developed multi-stage stochastic ORPS with a modified block scheduling policy. They considered nine different terms in the objective function including waiting time and tardiness for both scheduled and postponed surgeries and total overtime of each surgery block.
Hashemi Doulabi et al. \cite{hashemi2022state} and Hashemi doulabi and Khalilpourazari \cite{hashemi2023stochastic} proposed state-variable models to formulate stochastic operating room planning scheduling problems. These state-varaible models are closed-form mixed integer programming models that incorporate uncertainty in surgical durations. Also, Hashemi doulabi et al. \cite{hashemi2020vehicle} demonstrated that operating room scheduling problem with stochastic surgical durations can be formulated a synchronized vehicle routing problem with stochastic service times. 

Ghandehari and Kianfar \cite{ghandehari2022mixed} explored a hybrid CP and CG model to address the OR planning problem. Their approach considered a block strategy at the tactical level and surgery sequencing at the operational level. They focused solely on elective surgeries and assumed that the available planning time was enough to perform all the surgeries. Consequently, they did not account for the possibility of surgery postponement in their model.\\
Beyond these foundational works, several recent studies have introduced novel formulations and features to further refine OR planning and scheduling. 
Zhu et al. \cite{zhu2024surgical} formulated the Surgical Case Assignment Problem (SCAP) as a multi-objective optimization problem, addressing both total operating cost minimization and maximum scheduled surgeries. Additionally, Hashemi Doulabi et al. \cite{hashemi2012effective} formulated an open-shop scheduling problem that closely resembles the operating room scheduling problem as a mixed-integer programming model with Big-M constraints, and developed a hybrid simulated annealing algorithm to solve it. Li et al. \cite{li2024scheduling} formulated the ORPS problem as a parallel machine scheduling problem, and considered bed availability constraints to jointly optimize surgery scheduling and hospital bed assignments. Li et al. \cite{li2024scheduling} similarly studied the bed planning and OR scheduling to optimize efficiency.
Majthoub Almoghrabi and Sagnol \cite{almoghrabi2025surgery} proposed a two-stage stochastic programming model for OR scheduling with emergency patient integration. Their approach introduced a convex surrogate model for second-stage costs, which lead to significantly computational efficiency and solution quality improvement. Boccia et al. \cite{boccia2024integrated} tackled a real-world IORPS problem at a Naples hospital. They developed an ILP model tailored for emergency responsiveness. Their study compared different OR management strategies (open, block, and modified block scheduling) and demonstrated how their model could significantly reduce emergency surgery delays.

\subsection{Advances in Solution Methodologies}\label{Developed_Solution_Methods}
Various methods and methodologies have been developed by researchers to solve IORPS problems effectively. These approaches primarily focus on hybrid heuristics, metaheuristic algorithms, and machine learning-based solutions to enhance efficiency in scheduling operating rooms. 

\subsubsection{Heuristics and Metaheuristic solution approaches}\label{Developed_Solution_Methods_Heuristics}
Marques et al. \cite{marques2014scheduling} proposed a new GA to address the same problem proposed in  \cite{marques2012integer} in real-world settings with large-scale instances. They considered maximizing total planned surgery duration and number of scheduled surgeries as two separate objective functions for the problem. Consequently, they reported the results separately for each objective function. 
Lin and Chou \cite{lin2020hybrid} introduce four heuristic approaches to quickly generate feasible solutions for OR scheduling problem, which involves assigning surgeries to multifunctional ORs while minimizing overtime costs, idle time, and maximizing utilization. Additionally, they developed a Hybrid GA (HGA) that integrates these heuristics with local search and elite search mechanisms to enhance performance. They compared its performance against four heuristic algorithms: Longest Processing Time first (LPT), Shortest Processing Time first (SPT), Earliest Due Date (EDD), and a heuristic based on the ratio of LPT to EDD. Their results demonstrated that the proposed hybrid genetic algorithm outperformed all four heuristics. Kamran et al. \cite{kamran2019adaptive} developed a two-phase heuristic algorithm to solve the proposed multi-stage model for the planning and sequencing problem.
Hashemi Doulabi et al. \cite{hashemi2016constraint} introduced a CG algorithm enhanced with a series of dominance rules and a fast infeasibility detection procedure to solve IORPS. They showed that their constraint-programming-based enhanced CG can achieve a feasible solution with an average optimality gap of just under 3\% for instances involving up to 120 surgeries for a five-day planning horizon. 
Roshanaei, and Naderi \cite{roshanaei2021solving} introduced a logic-based Benders' decomposition (LBBD) algorithm. The authors incorporated novel feasibility and optimality cuts to transfer the information from subproblems to the master problem. Furthermore, they developed a Branch-and-Check (B\&C) algorithm and demonstrated that the B\&C, coupled with CP subproblems, yielded the best results.
Akbarzadeh and Maenhout \cite{akbarzadeh2024dedicated} developed a Branch-Price-and-Cut algorithm with column generation and cutting planes for simultaneous planning of patients and scheduling of surgeons in ORs. They considered patient characteristics and surgeon preferences and introduced intelligent branching and valid inequalities to enhance efficiency.
Bargetto et al. \cite{bargetto2023branch} proposed a new Benders’ cutting procedure for tightening the LP relaxation problem used in the CG. They demonstrated that in certain cases, their method could outperform the developed solution methodology by Hashemi Doulabi et al. \cite{hashemi2016constraint}. 
Mashkani et al. \cite{mashkani2021operating} studied the IORPS considering the maximization of patients' priority scores. They introduced a heuristic algorithm for assigning specialties and patients to ORs on various days within the planning horizon. They applied a first-fit strategy which assigns the patient to the first available time block considering the capacity limitation. They used a sorted list of patients according to their level of priority, and surgery duration attributes, respectively.
Aringhieri et al. \cite{aringhieri2022combining} similarly focused on maximizing patient priority and levelling bed occupancy as hierarchical objective functions. They proposed a new ad hoc neighbourhood local search matheuristic for the bi-objective optimization problem. The authors presented three distinct local search methods, each involving the fixation of certain variable values followed by model re-optimization in pursuit of improved solutions.
Dodaro et al. \cite{dodaro2020solving} proposed a new formulation for OR scheduling problems considering a team of surgeons and anesthesiologists. The authors utilized Answer Set Programming (ASP) to identify the optimal solution that maximizes the number of scheduled surgeries. They evaluated their method on four types of randomly generated instances. Results indicate that increasing the time interval from 10 to 60 minutes will increase the overall time efficiency score for ORs, surgeons, and anesthesiologists. It is noteworthy that the developed algorithm could only handle instances up to a five-day planning horizon.

\subsubsection{Machine Learning-based Algorithms}\label{Developed_Solution_Methods_ML}
Recent developments in Machine Learning (ML) and data-gathering techniques have inspired researchers to leverage ML concepts and algorithms for creating hybrid ML-metaheuristic algorithms aimed at solving the IORPS. Wang et al. \cite{wang2022adaptive} developed an adaptive-learning-based GA designed to identify the optimal solution of the IORPS in a multi-hospital context. Their approach took into account the post-anesthesia care section and aimed to minimize the total average recovery completion time for patients. They introduced new surgery allocation and sequencing re-arrangement based on the random and social learning activity which were embedded within the GA. The results indicated that the proposed method enhanced the previously developed GAs targeting the same problem setting. Eshghali et al. \cite{eshghali2023machine} developed a machine-learning-based algorithm to address both the tactical and operational levels of IORPS. They benefited from a random forest classifier to predict the duration of the surgery. They proposed a hierarchical method to solve a tactical level and an initial daily schedule at the operational level. In the next step, they perform a rescheduling process triggered by the arrival of emergency patients and their predicted surgery duration. Results demonstrated that the developed algorithm could be utilized for real-world applications and the rescheduling process led to a 10.5\% enhancement in operating theatre performance. Akbarzadeh and Maenhout \cite{akbarzadeh2024two} introduced a two-layer machine learning-guided GA for sequencing surgeries across multiple resource phases to optimize operating room efficiency. They used Support Vector Machine (SVM) Regressor to predict the feasibility of the patient schedule generated in the GA crossover and mutation processes. Ala and Goli \cite{ala2024incorporating} presented a hybrid machine learning-tabu search algorithm to  assign patients to ORs with fairness and efficiency considerations. They used convolutional neural networks (CNNs) to predict the objective function value of the candidate solutions.
Wang et al. \cite{wang2021solving} proposed a hybrid Genetic Algorithm–Column Generation (GA-CG) approach for solving the multi-depot electric vehicle scheduling problem. Although their application domain differs from the IORPS problem, the underlying solution methodology shares some similarities with our approach. Their method follows a two-phase framework: the first phase employs CG to produce a pool of feasible columns, while the second phase utilizes a genetic algorithm to identify an optimal subset of columns that constructs the final solution.

Another example of hybrid ML-metaheuristic algorithms is the hybridization of RLA and GA. The reinforcement learning concept is a subset of ML and it is based on learning from the interaction with an environment and making decisions to maximize the long-term profit. The RLA in the field of computer science was initially introduced by Sutton and Barto in 1998 \cite{sutton2018reinforcement}. The advantage of RLA over other ML algorithms is that it can make decisions dynamically based on the current state and the accumulated results of previous decisions. In recent years researchers in the field of optimization studied the hybridization of the RLA with the state-of-the-art solution methods to improve the results and performance of the solution methodologies. For instance, from the latest studies, Song et al. \cite{song2023rl} used RLA in GA to develop an integrated RL-GA algorithm to schedule a satellite position based on the demands and requests. Different methods were introduced by Sutton and Barto \cite{sutton2018reinforcement} and other researchers to implement RLA. In all mentioned methods and algorithms, three mechanisms of recording previous results and decisions, reward calculation, and selection of the next action or decision, are common while they are implemented differently.

\subsection{Novelty and contributions of this study}\label{Novelty_contributions}
Based on the discussion in the Introduction and the preceding literature review, this study offers distinct contributions in both the formulation of the IORPS problem and the design of the solution methodology. These contributions highlight the originality and novelty of our work in comparison to existing studies.
The core contribution of this work lies in the development of a more efficient mathematical model for the IORPS problem. Unlike traditional models that employ a four-index binary decision variable structure (e.g, $x_{crtd}$), as seen in the works of Marques et al. \cite{marques2012integer}, Roland et al. \cite{roland2006operating}, Hashemi Doulabi et al. \cite{hashemi2016constraint}, Roshanaei and Naderi \cite{roshanaei2021solving}, and Boccia et al. \cite{boccia2024integrated}, our model introduces a decision variable with a reduced indexing scheme, thus significantly decreasing the problem dimensionality. In addition to simplifying the variable structure, our formulation also involves fewer constraints, making the model computationally more tractable and scalable for real-world applications. 

From a methodological standpoint, this study is—to the best of our knowledge—the first to use RL to enhance a GA as the initial column generator in a CG algorithm, where the RL adaptively determines the best crossover and mutation strategies based on instance characteristics. In fact, unlike existing hybrid GA-CG approaches, the GA is used as an initial column generator to warm start the CG procedure, thereby improving the convergence rate and solution quality of the CG algorithm.
Although both our method and the hybrid GA-CG approach proposed by Wang et al. \cite{wang2021solving} share a similar high-level structure, the roles of the components are fundamentally different. Wang et al. used CG to generate feasible columns and GA as a column selector to construct the final solution. In contrast, we reverse this structure by using GA to pre-generate high-quality columns and relying on a MIP model within CG to identify the optimal subset of columns for the final solution.


\section{Problem definition}\label{Problem_Definition}
In this section, we present assumptions and some notations that will then be used to formulate the IORPS. In this IORPS, we consider the decisions for the surgery allocation to ORs and sequencing the allocated surgeries in the same OR over a finite planning horizon. We formulate the problems as a MIP model where all parameters including surgery durations, are assumed to be deterministic. The ultimate goal of the model is to determine on which day, in which OR, and at what time, the surgeries should start.

The proposed IORPS consists of a set of surgeries $(i\in I)$ with specified due dates $(d_i)$. Also, the surgeries should be scheduled within a planning horizon which, typically spanning a week or a limited number of days $(d\in D)$ or can be postponed to the next planning period. Surgeries with due dates falling within the planning horizon are categorized as mandatory $(i\in I_1)$ and must be planned within the current planning time frame. In contrast, other surgeries are classified as elective $(i\in I_2)$ and could be postponed. The division of surgeries into mandatory and elective categories is a realistic approach, as not all surgeries are essential or life-threatening. Indeed, surgeries can be categorized based on multiple factors such as deadline, or the patient priority score  (\cite{aringhieri2022combining}, \cite{mashkani2021operating}).

Each day is subject to constraints regarding a limited number of working hours, represented as a regular time shift $(t\in T_1)$, and a limited overtime possibility $(t\in T_2)$. The planning horizon is defined here as a fixed number of days with eight-hour regular time shifts divided into five-minute time slots, yielding a total of 96 time slots daily. Additionally, two-hour overtime shifts divided into five-minute time slots (a total of 24 time slots per day) are considered in each OR. If the model requires more than two hours of overtime daily for planning or scheduling, a practical solution could be to either open an additional operating room (OR) or postpone certain surgeries to the next planning period.

Additionally, there are a limited number of identical operating rooms in terms of both available time slots and equipment on each day $(r\in R_d)$ during the planning period which limits the service supply capacity of the hospital. It is worth mentioning that, there are a few studies that considered non-identical ORs with different equipment such as Ronald et all. (\cite{roland2006operating}, \cite{roland2010scheduling}).

Surgeons typically perform surgeries on patients with whom they are familiar. Thus, we assume that surgeries are preassigned to the surgeons $(l\in L)$, and the total surgical time for each surgeon is limited on each day. This ensures that the final plan for each surgeon on a given day complies with rules and regulations governing daily working hours. Our approach involves an open scheduling strategy, permitting the sharing of operating rooms among various surgeons and surgeries. 
Finally, the objective function is to minimize the postponement cost of patients not scheduled within the planning horizon, along with the overtime cost and fixed cost of opening ORs.

\section{Existing mathematical model}\label{Existing_mathematical_model}
In this section, we introduce the existing IORPS model which was initially presented in Marques et al. \cite{marques2012integer}. The original objective function of the IORPS aimed to maximize the cumulative duration of scheduled surgical procedures. Nonetheless, recognizing the significance of cost management, as deliberated in the sections \ref{Introduction} and \ref{Literature_Review}, we modified the IORPS to be compatible with the new objective function. The modified version is referred to as MIORPS. A comprehensive listing of the notations employed within the IORPS is presented in Table \ref{tab_1}. 

\begin{table}[h]
\centering
\footnotesize
\caption{MIORPS notation} 
\begin{tabular}{R{1cm} L{13cm}}
\hline
               &                \\
\multicolumn{2}{l}{Sets:}       \\
  $I_1:$             &  Set of mandatory surgeries              \\
  $I_2:$             &  Set of elective surgeries              \\
  $I:$               &  Set of all surgeries where $I=I_1 \cup I_2$              \\
  $L:$               &  Set of surgeons              \\
  $I^{'}_{l}:$       &  Set of surgeries that can be performed by surgeon $l$              \\
  $D:$               &  Set of days in the planning horizon              \\
  $K:$               &  Set of operating rooms              \\
  $K_d:$             &  Set of available operating rooms on day d              \\
  $T_1:$             &  Set of time slots related to normal shift               \\
  $T_2:$             &  Set of time slots related to overtime shift               \\
  $T:$               &  Set of total time slots where $T=T_1 \cup T_2$              \\
  $T_i:$             &  Set of time slots for surgery i in such a way that if surgery i starts, it finishes before the end of the available time slots in the day (T)             \\
                     &                \\
\multicolumn{2}{l}{Parameters:} \\
  $t_i:$             &  Duration of surgery $i$              \\
  $d_i:$             &  Due date of surgery $i$              \\
  $a_{ld}:$          &  Total available time for surgeon $l$ on day $d$              \\
  $c^{Pos}:$         &  Cost of surgery postponement              \\
  $c^{OR}:$          &  Cost of opening of operating room              \\
  $c^{Ovt}:$         &  Overtime cost               \\
  $c_{it}:$          &  Cost of surgery $i$ if it starts at the beginning of time slot $t$ where:               \\
                     &  $c_{it}=c^{Ovt}*\begin{cases}
0 & t+t_{i}< \left| T_{1} \right| \\ 
t_{i} & t+t_{i}\ge  \left| T_{1} \right|\text { and }t\ge  \left| T_{1} \right| \\ 
t+t_{i}-\left| T_{1} \right| & t+t_{i}\ge  \left| T_{1} \right|\text { and }t< \left| T_{1} \right|  
\end{cases}$              \\
                     &                \\
\multicolumn{2}{l}{Variables:} \\
  $x_{idkt}:$        &  1 if surgery $i$ starts at the beginning of time slot $t$ on day $d$ and operating room $k$, 0 otherwise              \\
  $y_{dk}:$          &  1 if operating room k opened on day $d$, 0 otherwise              \\
  $z_i:$             &  1 if surgery $i$ postponed to the next planning horizon, 0 otherwise               \\
                     &                \\
\hline
\end{tabular}
\label{tab_1}
\end{table}

The MIORPS is as follows:

\vspace{-10 pt}

\begin{flalign}\label{Eq_1}
 \textbf{MIORPS: }\;\; \min \sum_{i \in I_2} c^{Pos} z_i + \sum_{d \in D}\sum_{k \in K_d} c^{OR} y_{dk} + \sum_{i \in I}\sum_{d \in D:d \leq d_i}\sum_{k \in K_d}\sum_{t \in T_i} c_{it} x_{idkt} &&
\end{flalign}

\noindent \mbox{Subject to:}

\vspace{-10 pt}

\begin{alignat}{2}
& \sum_{d \in D:d \leq d_i}\sum_{k \in K_d}\sum_{t \in T_i} x_{idkt} = 1                                                && \forall\, i \in I_1 \label{Eq_2}\\
& \sum_{d \in D:d \leq d_i}\sum_{k \in K_d}\sum_{t \in T_i} x_{idkt} + z_i = 1  \qquad \qquad \qquad \qquad      && \forall\, i \in I_2 \qquad \qquad \qquad \qquad \qquad \qquad \qquad \qquad \label{Eq_3}
\end{alignat}
\begin{alignat}{2}
& \sum_{i \in I:d \leq d_i}\sum_{t' \in T_i:max[t-t_{i'},0] < t^\prime \leq t} x_{idkt^\prime} \leq 1                            && \forall\, d \in D, \forall\, k \in K_d, \forall\, t \in T \label{Eq_4}\\
& \sum_{i \in I^{'}_{l}:d \leq d_i}\sum_{k \in K_d}\sum_{t^\prime \in T_i:max[t-t_{i'},0] < t^\prime \leq t} x_{idkt^\prime} \leq 1 \qquad \qquad  && \forall\, l \in L, \forall\, d \in D, \forall\, t \in T \label{Eq_5}\\
& \sum_{i \in I^{'}_{l}:d \leq d_i}\sum_{k \in K_d}\sum_{t \in T_i} t_{i} x_{idkt} \leq a_{ld}                   && \forall\, l \in L, \forall\, d \in D \label{Eq_6}\\
& x_{idkt} \leq y_{dk}                                                                                           && \forall\, i \in I, \forall\, d \in D:d \leq d_i,  \forall\, k \in K_d, \forall\, t \in T_i \label{Eq_7}\\
& x_{idkt}, y_{dk}, z_i \in \{0,1\}                                                                         && \forall\, i \in I, \forall\, d \in D:d \leq d_i,  \forall\, k \in K_d, \forall\, t \in T_i \nonumber
\end{alignat}

The objective function minimizes the aggregate costs of surgery postponement, the cost of opening the operating room and the costs incurred due to surgeries conducted within overtime. These components are inherently interdependent, and their interactions must be carefully captured to reflect realistic operational trade-offs. For instance, opening additional ORs increases surgical capacity and can reduce the number of postponed procedures. However, doing so incurs fixed or operational costs. If the cost of postponement is relatively high, the model will favor opening more ORs to accommodate demand. Conversely, if OR opening costs dominate, the model may opt to defer lower-priority cases.
Similarly, the number of open ORs directly influences the extent of overtime incurred. With limited OR availability, the model may be forced to allocate surgeries beyond regular working hours, which increases overtime costs. In contrast, opening more ORs can distribute the workload more evenly, reducing the need for overtime but at the expense of higher fixed costs. Finally, allowing overtime can prevent postponements, particularly for high-priority cases, but the associated costs may push the model to postpone elective surgeries. 

Constraint (\ref{Eq_2}) ensures that all mandatory surgeries must be scheduled during the current planning horizon. Constraint (\ref{Eq_3}) specifies that elective surgeries either could be operated in the current planning horizon or should be postponed to the next planning horizon. Constraint (\ref{Eq_4}) ensures that on each day, and in each operating room, at most one surgery can performed at each time slot. In this constraint, the inner summation is equal to one if surgery $i$ starts before time slot $t$ and finishes after this time slot, meaning that it is ongoing in time slot $t$.
Constraint (\ref{Eq_5}) enforces that each surgeon can perform at most one surgery at a time. Constraint (\ref{Eq_6}) ensures the respect of the availability time of each surgeon on each day. Finally, Constraint (\ref{Eq_7}) states that if a surgery is scheduled, the corresponding operating room must be opened to perform the surgery. 

Expanding the problem size, whether by including more surgeries, additional time slots, or extending the planning horizon, significantly increases the number of variables and constraints in the MIORPS. As a result, solving larger instances of the model becomes considerably more challenging. For instance, Lin and Chou \cite{lin2020hybrid} characterized the size of their studied instances based on the combination of available operating rooms on day $d$ denoted as  $M_d$, and the number of surgeries, $N$. They categorized problems as small when $M_d=4$ and $N=20,\; 30,\; 40$, and as large when $M_d \times N = 6\times50,\; 6\times80,\; 10\times100,\; 10\times200$. In our study, we go beyond these limitations by considering $M_d=5$ and $N = 40,\; 60,\; 80,\; 100,\; 120$ for small instances, and $M_d \times N=10\times200,\; 10\times240,\; 15\times300,\; 15\times360,\; 20\times400,\; 20\times480$ for large-scale instances. Further details about the instance definitions are provided in Sections \ref{Definition_Instances_Synthetic} and \ref{Definition_Instances_Real}.

\section{Proposed mathematical model}\label{Proposed_Model}
We propose a novel mathematical model to solve IORPS. The new model is achieved by reformulating the IORPS developed by Marques et al. (2012). The distinctive configuration of the MIORPS - where each operating room (OR) possesses dedicated resources - allows for a reduction in the MIORPS's scale, along with the elimination of the dimension associated with the operation rooms. Therefore, the proposed model, which is referred to as PMIORPS, has a smaller number of variables and constraints. From a complexity perspective, it is well-known that operating room planning and scheduling problems belong to the category of NP-Hard problems. This claim is supported by the complexity analysis of a simplified version of the problem, in which for a single day scheduling problem, only the fixed cost of operating rooms is considered. This simplified case is known as the bin packing problem, a problem that has been proven to be NP-hard in the literature  \cite{korte2008combinatorial}. Therefore, our original operating room planning and scheduling problem is also NP-hard. Table \ref{tab_2} shows the new variables of the PMIORPS.

\begin{table}[H]
\centering
\small
\caption{The new variables of the PMIORPS} 
\begin{tabular}{R{1cm} L{13cm}}
\hline
                     &                \\
\multicolumn{2}{l}{New variables:} \\
  $x_{idt}:$         &  1 if surgery $i$ starts at the beginning of time slot $t$ on day $d$, 0 otherwise              \\
  $y_d:$             &  Number of opened operating rooms on day $d$              \\
                     &                \\
\hline
\end{tabular}
\label{tab_2}
\end{table}

\begin{flalign}\label{Eq_27}
 \textbf{PMIORPS: }\;\; \min \sum_{i \in I_2} c^{Pos} z_i + \sum_{d \in D} c^{OR} y_d + \sum_{i \in I}\sum_{d \in D:d \leq d_i}\sum_{t \in T_i} c_{it} x_{idt} &&
\end{flalign}

\noindent \mbox{Subject to:}
\begin{alignat}{2}
& \sum_{d \in D:d \leq d_i}\sum_{t \in T_i} x_{idt} = 1                                                                                    && \forall\, i \in I_1 \label{Eq_28}\\
& \sum_{d \in D:d \leq d_i}\sum_{t \in T_i} x_{idt} + z_i = 1                                                                              && \forall\, i \in I_2 \label{Eq_29}\\
& \sum_{i \in I:d \leq d_i}\sum_{t' \in T_i:max[t-t_{i'},0] < t^\prime \leq t} x_{idt^\prime} \leq y_d                                   && \forall\, d \in D, \forall\, t \in T \label{Eq_30}\\
& \sum_{i \in I^{'}_{l}:d \leq d_i}\sum_{t^\prime \in T_i:max[t-t_{i'},0] < t^\prime \leq t} x_{idt^\prime} \leq 1 \qquad \qquad \qquad  && \forall\, l \in L, \forall\, d \in D, \forall\, t \in T \label{Eq_31}\\
& \sum_{i \in I^{'}_{l}:d \leq d_i}\sum_{t \in T_i} t_{i} x_{idt} \leq a_{ld}                                                       && \forall\, l \in L, \forall\, d \in D \label{Eq_32}\\
&  y_d \leq | K_d|                                                                                                                         && \forall\, d \in D \label{Eq_33}\\
& x_{idt}, z_i \in \{0,1\}                                                                                                                 && \forall\, i \in I, \forall\, d \in D:d \leq d_i, \forall\, t \in T_i \nonumber\\
& y_d \in \mathbb{Z}                                                                                                                       && \forall\, d \in D \nonumber
\end{alignat}

Generally, PMIORPS and MIORPS follow the same logic. However, using the new variable $x_{idt}$ empowered us to eliminate the constraint (\ref{Eq_7}). Constraints (\ref{Eq_28}), (\ref{Eq_29}), (\ref{Eq_31}), and (\ref{Eq_32}) follow the same concept of constraints (\ref{Eq_2}), (\ref{Eq_3}), (\ref{Eq_5}), and (\ref{Eq_6}) respectively. Constraint (\ref{Eq_30}) ensures that at any time slot $t$, the number of ongoing surgeries is less than or equal to the number of opened ORs. Finally, constraint (\ref{Eq_33}) enforces the upper bound on the number of ORs on each day. 

\begin{theorem}\label{Tr_1}
The PMIORPS model is a valid reformulation of the MIORPS Model. 
\end{theorem}

We provide the proof of this theorem in Appendix \ref{AppendixA}. In this proof, we first prove that the developed alternative model is a valid representation of MIORPS model. Then we show that for any feasible solution of the MIORPS, there always exists a feasible solution for the alternative model.

\section{Column generation algorithm}\label{CG_Algorithm}
In this section, we introduce the CG algorithm based on the model we proposed in the previous section, to address the challenges posed by large-scale instances when we deal with a large number of surgeries in an extended planning horizon. We can use commercial solvers to solve PMIORPS directly and this approach is powerful enough to handle relatively large problems. However, the capability of the solver to achieve zero-gap optimal solutions reduces by increasing the size of the problem. Therefore, we have developed a new pattern-based decomposition model based on the PMIORPS. In the following sections, first, we propose the master problem and subproblem of the CG approach. Then, we will discuss the hybrid reinforcement learning algorithm we developed to generate the initials columns for the CG procedure.

\subsection{Master problem and subproblem}\label{CG_mathematical_model}

We introduce the following model, denoted as M-CG, which serves as the master problem in our CG algorithm. The purpose of M-CG is to determine the optimal selection of patterns across the planning horizon, where each pattern corresponds to a feasible daily surgery schedule. The notation used in formulating this master problem is summarized in Table \ref{tab_3}.

\begin{figure}[h]
   \centering
   \includegraphics[width=0.8\textwidth]{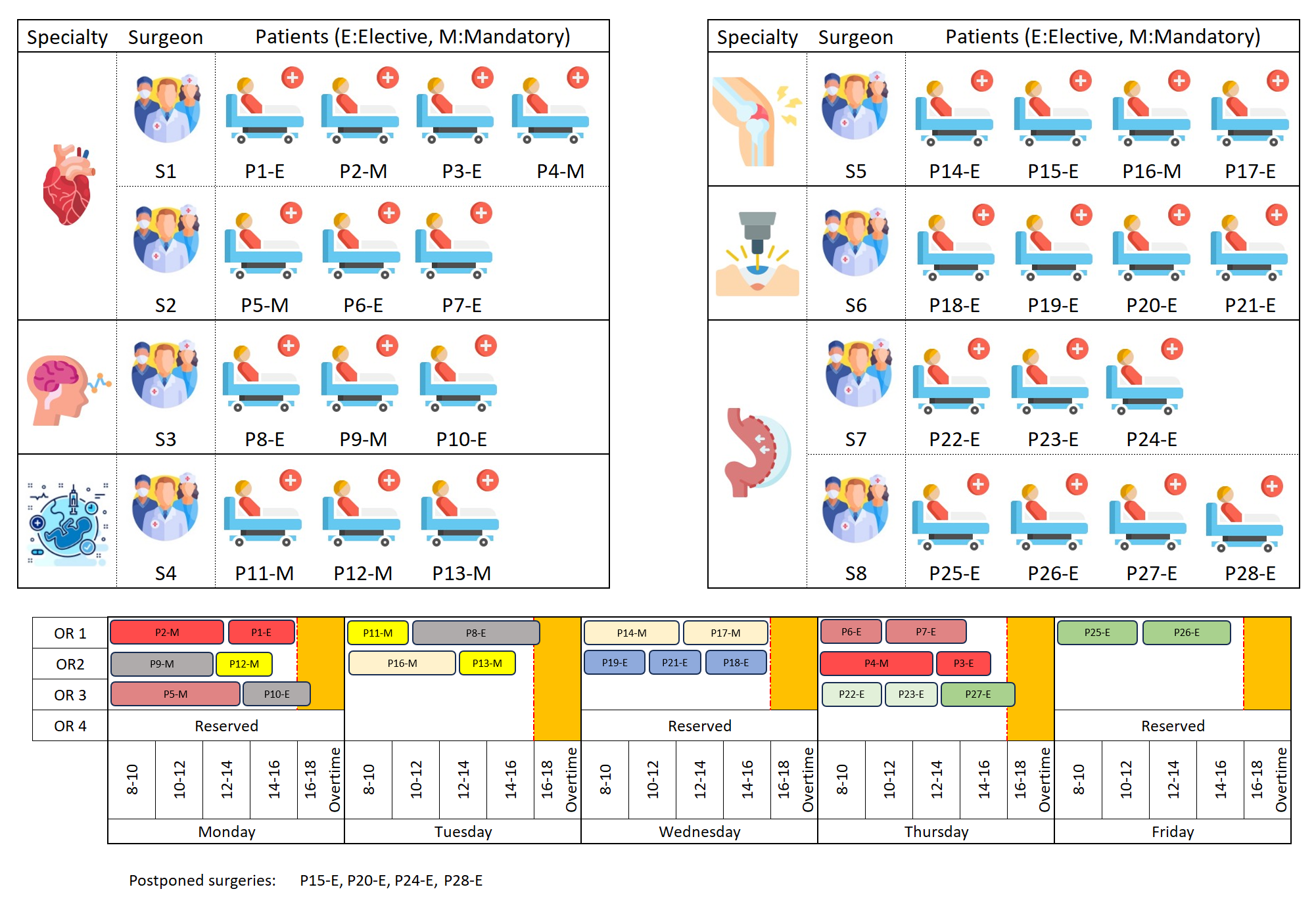}
   \caption{Illustration of the sample surgery plan and schedule output}
   \label{fig_1}
\end{figure}

To clarify the interpretation of the newly defined sets in Table \ref{tab_3}, we provide an illustrative example in Figure \ref{fig_1}, which depicts a representative final solution. In this example, there are 28 elective and mandatory patients assigned to 8 different surgeons over a 5-day planning horizon. Each day consists of 8 regular working hours and 2 hours of allowed overtime.
Consider the schedule for Friday: only two elective patients (Patient 25 and Patient 26) are assigned. There may exist alternative feasible ways to schedule these two patients. For instance, one pattern might schedule Patient 25 at 08:00 in OR1 and Patient 26 at 12:00 in OR3. Another pattern could reverse the OR assignments with Patient 25 at 08:00 in OR3 and Patient 26 at 12:00 in OR1. Each of these alternatives represents a distinct feasible pattern for Friday, potentially with identical or differing costs. If there are exactly three such valid combinations, the set $j_d$, $d=5$ (Friday) would contain these three patterns. 

Applying the same logic to the rest of the week, each day $d\in{1,2,3,4,5}$ (Monday to Friday) may have a different number of feasible daily patterns. For instance, the sizes of the sets $J_d$ could be 8, 5, 9, 6, and 3, respectively. Consequently, the total number of patterns in the planning horizon, represented by the union set $J$, would be $|J|=\sum_{d=1}^{5} |J_d|=31$. Moreover, it is important to note that patients can appear in multiple patterns across different days. For example, Patient 18 might be feasibly scheduled on Tuesday at 10:00 in OR4, or alternatively, on Thursday at 08:00 in OR4. Assuming these are the only patterns in which Patient 18 appears, the set $J'_i$ (representing all patterns involving patient $i=18$) would include these two patterns. This concept naturally extends to all other patients in the model, each associated with their corresponding pattern set $J'_i$.

Since both $x_j$ and $z_i$ are binary variables, the M-CG is an IP model and we cannot apply standard duality to calculate the $\pi^{1}_{i}$, $\pi^{2}_{i}$, and $\pi^{3}_{d}$ values directly. Therefore, we should relax the integrality of $x_j$ and $z_i$ variables and let them take any real values between 0 and 1 before solving the master problem. Then we can dynamically calculate $\pi^{1}_{i}$, $\pi^{2}_{i}$, and $\pi^{3}_{d}$ values and use them in the subproblem to produce and add new columns to the master problem.

\begin{table}[h]
\centering
\small
\caption{The new notation of the CG} 
\begin{tabular}{R{1cm} L{13cm}}
\hline
               &                \\
\multicolumn{2}{l}{New sets:}       \\
  $J:$               &  Set of all feasible daily surgical schedules (patterns or columns) for all operating rooms              \\
  $J_d:$             &  Set of patterns for day $d$              \\
  $J'_{i}:$       &  Set of patterns including surgery $i$              \\
                     &                \\
\multicolumn{2}{l}{New parameters:} \\
  $c_j:$             &  Calculated cost of column $j$              \\
  $\pi^{1}_{i}:$     &  Reduced costs of set of constraints (\ref{Eq_35})            \\
  $\pi^{2}_{i}:$     &  Reduced costs of set of constraints (\ref{Eq_36})            \\
  $\pi^{3}_{d}:$     &  Reduced costs of set of constraints (\ref{Eq_37})            \\
                     &                \\
\multicolumn{2}{l}{New variables:} \\
  $x_j:$             &  1 if column $j$ selected, 0 otherwise              \\
                     &                \\
\hline
\end{tabular}
\label{tab_3}
\end{table}

\begin{flalign}\label{Eq_34}
 \textbf{M-CG: }\;\; \min \sum_{j \in J} c_j x_j + \sum_{i \in I_2} c^{Pos} z_i &&
\end{flalign}

\noindent \mbox{Subject to:}
\begin{alignat}{2}
& \sum_{j \in J'_{i}} x_j \geq 1                                       && \forall\, i \in I_1 \label{Eq_35}\\
& \sum_{j \in J'_{i}} x_j + z_i \geq 1                                 && \forall\, i \in I_2 \qquad \qquad \qquad\label{Eq_36}\\
& \sum_{j \in J_d} x_j \leq 1 \qquad \qquad \qquad \qquad \qquad \qquad \qquad \qquad \qquad\qquad \qquad \qquad                        && \forall\, d \in D \label{Eq_37}\\
& x_j, z_i \in \{0,1\}                                                    && \forall\, i \in I, \forall\, j \in J \nonumber
\end{alignat}

The objective function of the M-CG is to minimize the total cost of selected columns as well as the cost of postponing surgeries to the next planning horizon. Constraint (\ref{Eq_35}) is guaranteeing the selection of at least one pattern for mandatory surgeries. Constraint (\ref{Eq_36}) implies that either at least one pattern should be selected for the elective surgeries, or they must postponed to the next planning horizon. Constraint (\ref{Eq_37}) enforces that for each day in the planning horizon, we are allowed to select at most one pattern. This restriction reflects a practical consideration where some problem instances may feature excessive resources relative to demand, it is not always necessary or efficient to schedule surgeries every day. In such cases, the model allows for the possibility of having no pattern selected on a given day that no surgeries are scheduled. This flexibility is essential to avoid unnecessary operating room usage or idle resource allocation when there is no surgical demand. Therefore, for each day, the model permits the selection of at most one pattern from the column pool, including the option of selecting none when appropriate. We have also proposed the following model S-CG as the subproblem of our CG algorithm. It is worth noting that the structure of the formulated subproblem closely resembles that of the classical bin packing problem, which has been proven to be NP-hard \cite{korte2008combinatorial, schepler2022solving}. In addition, the subproblem is strongly NP-hard, as it requires constructing a hospital-wide daily schedule across all operating rooms rather than scheduling a single room on a given day.

\begin{flalign}\label{Eq_38}
\begin{split}
 \textbf{S-CG: }\;\; &\min \left( \sum_{d \in D} c^{OR} y_d + \sum_{i \in I}\sum_{d \in D:d \leq d_i}\sum_{t \in T_i} c_{it} x_{idt} \right) \\
 & - \left( \sum_{i \in I_1} \pi^{1}_{i} \Bigl( \sum_{d \in D:d \leq d_i}\sum_{t \in T_i} x_{idt} \Bigl) -\sum_{i \in I_2} \pi^{2}_{i} \Bigl( \sum_{d \in D:d \leq d_i}\sum_{t \in T_i} x_{idt} \Bigl) -\sum_{d \in D} \pi^{3}_{d} \right)
 \end{split}
\end{flalign}

\noindent \mbox{Subject to:}
\begin{alignat}{2}
& \sum_{d \in D:d \leq d_i}\sum_{k \in K_d}\sum_{t \in T_i} x_{idt} \le 1                                                                  && \forall\, i \in I \label{Eq_39}\\
& \sum_{i \in I:d \leq d_i}\sum_{t' \in T_i:max[t-t_{i'},0] < t^\prime \leq t} x_{idt^\prime} \leq y_d                                   && \forall\, d \in D, \forall\, t \in T \label{Eq_40}\\
& \sum_{i \in I^{'}_{l}:d \leq d_i}\sum_{t^\prime \in T_i:max[t-t_{i'},0] < t^\prime \leq t} x_{idt^\prime} \leq 1  \qquad \qquad \qquad\qquad \qquad \qquad                     && \forall\, l \in L, \forall\, d \in D, \forall\, t \in T \label{Eq_41}\\
& \sum_{i \in I^{t^\prime}_{l}:d \leq d_i}\sum_{t \in T_i} t_{i} x_{idt} \leq a_{ld}                                                       && \forall\, l \in L, \forall\, d \in D \label{Eq_42}\\
\begin{split}
&\left( \sum_{d \in D} c^{OR} y_d + \sum_{i \in I}\sum_{d \in D:d \leq d_i}\sum_{t \in T_i} c_{it} x_{idt} \right) \\
& - \sum_{i \in I_1} \pi^{1}_{i} \Bigl( \sum_{d \in D:d \leq d_i}\sum_{t \in T_i} x_{idt} \Bigl) \\
&-\sum_{i \in I_2} \pi^{2}_{i} \Bigl( \sum_{d \in D:d \leq d_i}\sum_{t \in T_i} x_{idt} \Bigl) -\sum_{d \in D} \pi^{3}_{d} \le 0  
\end{split}                                                                                                                                && \label{Eq_43}\\
& x_{idt} \in \{0,1\}                                                                                                                      && \forall\, i \in I, \forall\, d \in D:d \leq d_i, \forall\, t \in T_i \nonumber
\end{alignat}

The objective function of the subproblem is to find the feasible pattern or column for each day $d$. The objective function of the S-CG is aimed to minimize the reduced cost of the new variables that are going to be generated by the subproblem and be added to the master problem.

In fact, the subproblem searches for the column with the highest negative objective function. Definition of constraints (\ref{Eq_39})-(\ref{Eq_42}) is similar to constraints (\ref{Eq_28})-(\ref{Eq_32}) respectively. Constraint (\ref{Eq_43}) enforces that any feasible solution should have a negative objective function. Therefore, if the S-CG fails to find a feasible solution, there is no more column to add to the columns pool. 

During the subproblem's execution, identifying the column featuring the most negative cost can sometimes be time-intensive \cite{abbaszadeh2026drone}. To boost the CG and increase the chance of its convergence, we implemented a fast feasible column detection technique. In this method, when the subproblem finds a feasible solution, we allow the searching process to continue for a couple of seconds hoping to find a better solution. When the mentioned time limit passes, the subproblem generates and adds a new column to the column pool based on the incumbent solution. When CG converges, the objective value of the solution obtained by solving linear programming relaxation of the master problem represents a lower bound of the problem. The upper bound will be obtained by solving the master problem, in which the integrality requirements of $x_j$ and $z_i$ variables are imposed to the model.

\subsection{Initial columns generator}\label{RGA_CG}
In this section, we discuss a reinforcement-learning-based heuristic algorithm specifically designed to generate the initial columns. CG needs at least one feasible solution for constructing the master problem model and calculating the reduced costs of the master problem's constraints. To find an initial feasible solution, we adopted the first-fit strategy described by Mashkani et al. \cite{mashkani2021operating}. We used the sorted list of surgeries based on their due date considering the lower surgery duration as a tie-breaking rule, and random selection in case of equality of surgery durations. We combined the RLA and GA and introduced Hybrid RGA-CG to increase the selection chance of operators which contributes more effectively to improving the solution. We have implemented three types of crossover operators: Partially Mapped Crossover (PMX), Order crossover (OX1) and ordered-based crossover (OX2). Additionally, we considered four types of mutation operators including swap, insertion, scramble, and inversion. As the GA searches for a solution, RLA records the selection of GA's operators along with their corresponding outputs. This information is then used to select the best operator based on the problem's structure.  Figure \ref{fig_3} illustrates the steps of the Hybrid RGA-CG. Solid blue lines show the procedure of the algorithm, while the dashed orange line shows the flow of data exchange between the dataset of columns and the algorithm's components. The Hybrid RGA-CG has two RGA phases, which are designed to search the solution space for better columns. Each phase also has its own enhancement procedure developed to perform a local search for a feasible solution. When the initial column generator algorithm meets the stopping criteria, the CG loop starts and solves the master problem and subproblem iteratively using MIP optimization solver. Upon finishing the CG process, either the total time limit passed or the CG algorithm converged, solving the master problem will provide the best solution based on the found columns. A detailed description of the Hybrid RGA-CG algorithm is presented in Appendix D.  

\begin{figure}[h]
	\centering
	\includegraphics[width=\textwidth]{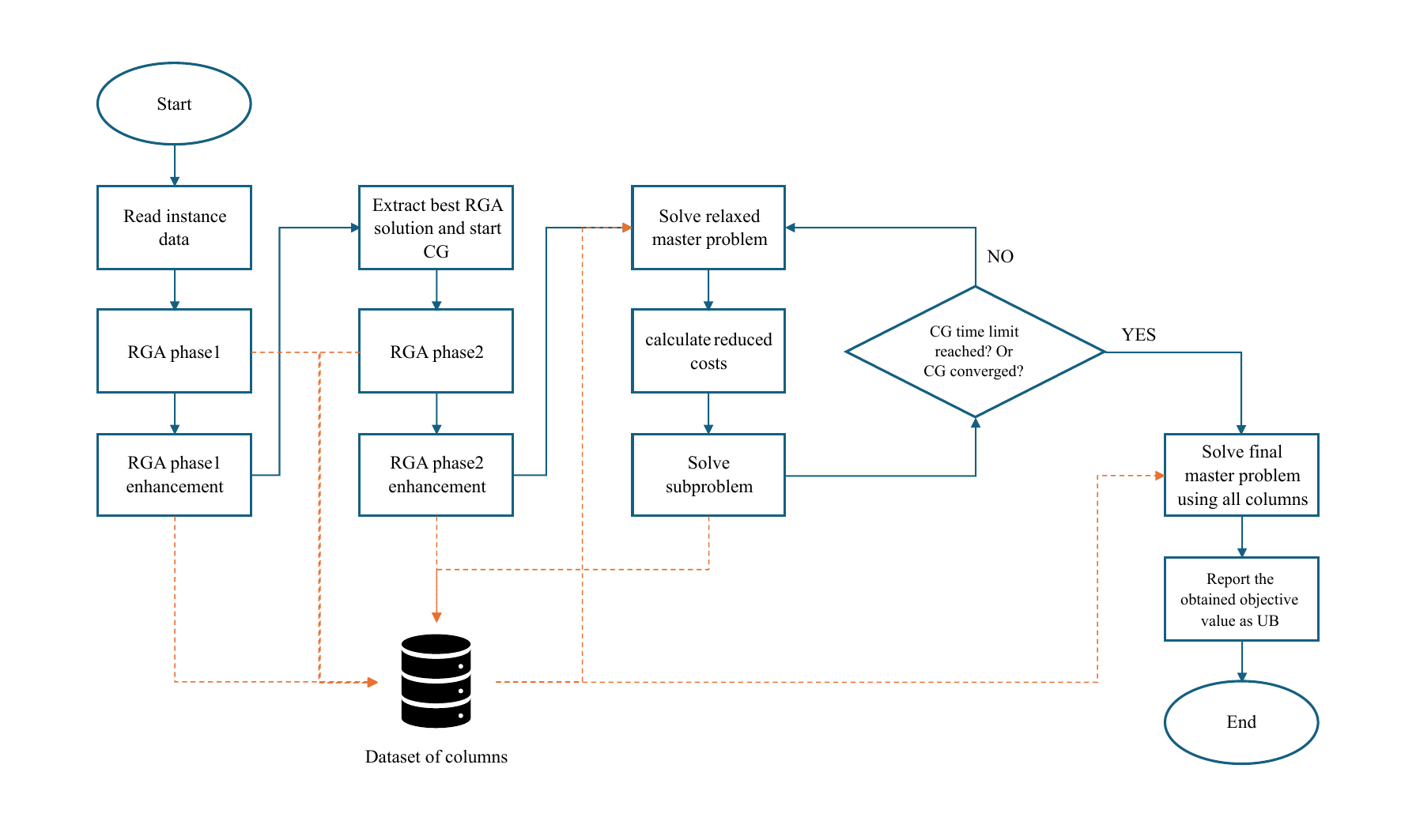}
	\caption{Flowchart of the column generation heuristic}
	\label{fig:cg_algorithm}
\end{figure}

\subsection{Integration of RLA and GA}\label{Hybrid_RL_GA}
As mentioned in the previous section, developed GA in this study, benefited from four mutation operators and three crossover operators to improve the exploration of the solution space. We cannot anticipate which operator will be the best in each category because each problem exhibits unique characteristics of it's solution space. In other words, we cannot say that for example, PMX is the best crossover operator for all cases and instances. Tailoring specific GA operators to each problem instance allows for better adaptation to the problem's characteristics. Therefore, We require a method that chooses the best operator according to the internal feedback developed based on the feasibility check and objective value. To achieve this purpose, we integrated RLA with GA as an adaptive operator selection (AOS) method in each iteration of the GA. In each iteration, RLA will determine which operator should be selected to perform crossover and mutation from the available operators for crossover and mutation.

In this study, we implement the Dynamic Multi-Armed Bandit (D-MAB), a variant of the Upper Confidence Bound (UCB) algorithm originally proposed by Da Costa et al. \cite{dacosta2008adaptive}. UCB is a well-known stochastic approach that has been widely applied in the development of recommendation systems \cite{afsar2022reinforcement}. This method was initially developed for the multi-armed bandit problem and it quickly gained researchers' attention due to its simplicity and effectiveness \cite{afsar2022reinforcement}. The classic multi-armed bandit strategy optimally balances exploration and exploitation under the assumption of independent, stationary rewards. However, many problems violate these assumptions. In other words, rewards can change over time or depend on the sequence of actions. In a reinforcement learning (RL) context, an agent’s actions influence not only immediate rewards but also the next state of the environment. This breaks the independent trial assumption. To tackle this issue, researchers have developed several approaches to justify and modify UCB. for example, Auer et al. \cite{auer2008near} developed an improved version of the UCB for RL called UCRL2. This algorithm maintains confidence intervals for the state-action values or model parameters and always choose the policy that is optimistic (upper-bound optimal) with respect to the confidence bounds. Furthermore, DaCostaet al. \cite{dacosta2008adaptive} introduced a method called Dynamic Multi-Armed Bandit (D-MAB) for operator selection, which augments the UCB1 algorithm with a change-detection mechanism. This effectively resets the accumulated reward estimates, allowing the method to forget outdated information when the environment shifts. Therefore, with the reset mechanism, D-MAB remains effective even though the reward probabilities are not independent over time. Additionally, in our problem setting, each arm corresponds to a distinct operator within either the crossover or mutation categories. These operators represent the available choices that an agent can select during the evolutionary search process. The agent's objective is to sequentially choose among these multiple operators or arms in a manner that maximizes the cumulative reward over time.

\begin{figure}[h]
    \centering
    \includegraphics[width=1\textwidth]{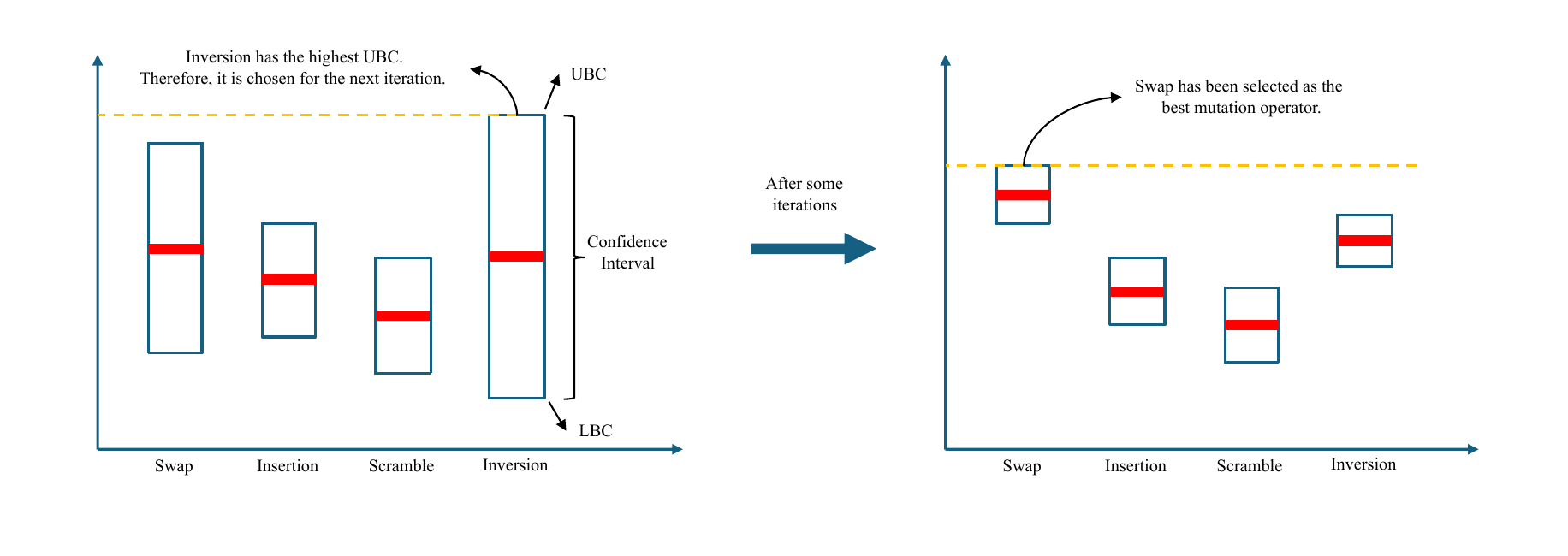}
    \caption{Procedure of mutation operator selection by UCB in RGA}
    \label{fig_5}
\end{figure}

The UCB agent maintains a record of both the total number of trials and the corresponding outcomes for each operator. After each iteration, it updates the cumulative reward values based on the success or failure of the trial. Through this iterative process, the agent continuously learns and refines its knowledge of the environment with which it interacts. In this study a successful trial is considered as obtaining a feasible solution or improvement in fitness value, while a failed trial refers to obtaining an infeasible solution. Then, the algorithm updates the UCB values for each operator based on the Equation (\ref{Eq_44}) \cite{chen2019context}. In Equation (\ref{Eq_44}), $TR_{op}$ is the total rewards gained by operator $op$, and $TOS_{op}$ is the total number of times operator $op$ selected. It is worth mentioning that, we used a failed trial penalty when the result becomes infeasible to make sure that the randomness of the GA's evolution will not result in the dominance of one operator. The failed trial penalty absorbs the random successful rewards and guarantees that over time and during the evolution process, proper mutation and crossover operators will be selected based on the structure of each problem. In cases where a single operator consistently outperforms the others by more than 50\% difference, we interpret this as an indication of an environmental shift. To adapt to the new conditions, we reset the accumulated rewards, and allowing the algorithm to discard outdated information and re-evaluate operator performance. At the end, the algorithm selects the operator with the maximum UCB value for the next iteration. As an example, Figure \ref{fig_5} shows how UCB could lead to selecting the best mutation operator in RGA.

\vspace{-0.25cm}

\vspace{-0.5cm}

\begin{flalign}\label{Eq_44}
 UCB_{op} = \frac{TR_{op}}{TOS_{op}} + 2\sqrt{\frac{2\ln(TOS_{op})}{TOS_{op}}}&&
\end{flalign}

\vspace{-0.25cm}

\section{Computational Experiments}\label{Computational_Experiments}
In this section, we conduct a comprehensive comparison of mathematical models using synthetic and real-world instances, examining their performance when solved directly through commercial optimization software such as CPLEX. Additionally, we aim to assess the efficacy of our developed solution approach by comparing it against existing state-of-the-art algorithms in addition to four heuristic algorithms.

We implemented the algorithms in a C++ environment using the IBM ILOG CPLEX Optimization Studio V20.10 as an optimizer engine. Experiments were run on a Compute Canada Béluga cluster where each compute node is equipped with two Intel Gold 6148 Skylake @ 2.4 GHz and a total of 40 cores for each node. We ran different instances using one core and 32GB of memory limit.

\subsection{Synthetic Instances}\label{Definition_Instances_Synthetic}
An analysis of real-world datasets with more than 8000 surgeries reveals that the average surgery duration is 130.45 minutes, with a standard deviation of 97.23 minutes \cite{costa2017assessment}. Moreover, Zhao et al. \cite{zhao2019machine} and Miller et al. \cite{miller2023using} applied machine learning techniques to predict the surgery duration in their studies based on different datasets. Comparing the average and standard deviation of surgery duration among the three above mentioned papers shows that the obtained results by Costa \cite{costa2017assessment} are valid and we can use them to generate instances for our computational experiments.

We categorized the instances based on the number of surgeries and generated five instances for each category. We generated instances ranging from 40 to 120 surgeries, referred to as small and medium-sized instances, with a fixed planning horizon of five days. In addition, some health centers may have to deal with more surgeries or they want to consider an extended planning horizon. For this case, we generated instances with 200, 240, 300, 360, 400, and 480 surgeries with a planning horizon of 10, 15, and 20 days. We assume that all surgeries will fall uniformly into the interval $130.45 \pm 97.23 \approx U\sim[30,230]$. We considered a five-minute interval units in generating the surgical durations. We also considered the total available time for surgeons on each day ($a_{ld}$) as in Table \ref{tab_4}; these values are the same as those reported in Table 1 of Fei et al. \cite{fei2009solving}. Moreover, we supposed 8 hours of regular time and 2 hours of overtime available for each OR on each day. We generated the due date of each surgery uniformly within the interval [1, 14] days as recommended by Fei et al. \cite{fei2009solving}. Also, according to Denton et al. \cite{denton2010optimal}, we set the overtime cost as $c^{Ovt}= c^{OR} / |T_2| = 1000 / 24$. This means the two hours of overtime is equivalent to opening a new operating room. 

Extending the horizon beyond 20 working days (equivalent to four weeks) is, however, rarely practical. Detailed planning and scheduling for such long horizons are of limited use in real settings, as several factors such as surgeon availability, patient health condition and priority, and operating room availability, are subject to frequent change. In practice, for longer horizons, hospitals typically adopt block scheduling strategies, while detailed day-to-day schedules are reserved for shorter horizons. In addition, most hospitals operate only a small number of ORs and perform fewer than 200 surgeries per week, as observed in the real-world instances from Naples. Thus, scenarios involving substantially larger numbers of ORs are of limited practical relevance. Nevertheless, to benchmark the scalability of our approach, we also generated very large instances with 800 and 960 surgeries over a 40-day horizon. On these instances, all executions failed due to memory exhaustion under the computational resources described in Section \ref{Computational_Experiments}. These results indicate that while the proposed algorithm is highly efficient for instances up to 480 surgeries, tackling problems with more than 800 surgeries requires additional computational strategies. Potential directions include parallel or distributed computing, as well as memory-optimized data structures. Such enhancements would be necessary to extend the applicability of the method to extremely large-scale cases, although these remain uncommon in most hospital environments.

\begin{table}[]
\centering
\caption{Value of daily availability of surgeons}
\begin{tabular}{lccccc}
\hline
 & \multicolumn{5}{c}{Day} \\ \cline{2-6} 
 & \multicolumn{1}{l}{Monday} & \multicolumn{1}{l}{Tuesday} & \multicolumn{1}{l}{Wednesday} & \multicolumn{1}{l}{Thursday} & \multicolumn{1}{l}{Friday} \\ \hline
Surgeon 1 & 8 & 0 & 7 & 0 & 6 \\
Surgeon 2 & 8 & 4 & 5 & 6 & 5 \\
Surgeon 3 & 8 & 3 & 6 & 7 & 8 \\
Surgeon 4 & 5 & 3 & 4 & 8 & 8 \\
Surgeon 5 & 6 & 5 & 0 & 6 & 8 \\
Surgeon 6 & 6 & 0 & 5 & 7 & 8 \\
Surgeon 7 & 0 & 6 & 6 & 6 & 8 \\
Surgeon 8 & 0 & 6 & 6 & 8 & 8 \\ \hline
\end{tabular}
\label{tab_4}
\end{table}

\subsection{Real-world Instances}\label{Definition_Instances_Real}
In addition to the synthetic data, we also used real-world data from a local hospital in Naples. In this setting, four operating rooms are available each day, from Monday to Saturday, between 8:00 AM and 8:00 PM. However, on Saturdays, Operating Room 2 is reserved exclusively for emergency cases. The dataset comprises 30 different instances, with patient sizes ranging from 110 to 164. Based on the synthetic data setting described earlier, this range can be classified as small- to medium-sized instances. To evaluate the performance of the proposed model and algorithms on large-scale instances, we generated an extended dataset by combining multiple instances. For example, instance 31 is constructed by merging instances 1 and 2; instance 32 by combining instances 3 and 4; and so on. Additionally, instance 46 is created by aggregating instances 1 to 3, while instance 55 is generated by combining instances 28 to 30. For more details regarding the original 30 instances, please refer to the work by Boccia et al. \cite{boccia2024integrated}

\subsection{Hyperparameters fine tuning}\label{GA_Parameters_Fine_Tuning}
This section details the hyperparameter fine tuning for the GA, the RL components, and the Bayesian Optimization (BO) algorithm. GA parameters are tuned via BO, RL parameters are examined through a targeted sensitivity analysis, and BO’s own hyperparameters are selected by a full factorial study.

\subsubsection{GA Hyperparameters fine tuning}\label{GA_Parameters_Fine_Tuning}
The performance of a GA heavily depends on the careful tuning of its hyperparameters. These hyperparameters, such as population size, crossover and mutation rates, and selection strategies, determine the balance between exploration (diversity in the solution space) and exploitation (intensifying the search around promising areas). The selected values of the mentioned hyperparameters significantly influence solution quality and convergence speed \cite{mao2021boosted}. There are several methods developed by researchers for finding the optimal values of the parameters involved in the optimization problem including but not limited to Manual Search, Grid search, Random search, Bayesian Optimization, Meta-heuristic approaches, and Adaptive Tuning \cite{feurer2019hyperparameter, mao2021boosted, ippolito2022hyperparameter}. Each of the above-mentioned algorithms has their own benefits and disadvantages. For example, Grid Search systematically evaluates all combinations of parameter values from predefined ranges. Therefore, it is exhaustive and straightforward and at the same time, computationally expensive, especially for large parameter spaces. However, Random Search is often faster than grid search but it may require a large number of samples for accurate results \cite{yang2020hyperparameter}.

Bayesian Optimization, on the other hand, has emerged as a powerful alternative due to its ability to balance exploration and exploitation systematically through probabilistic modeling. BO is particularly advantageous when dealing with expensive black-box functions, as it constructs a surrogate model to approximate the objective function and minimizes the number of evaluations required \cite{garrido2020dealing, pourmohamad2021bayesian}. This approach is especially effective when the number of evaluations is constrained, providing significant computational savings while achieving near-optimal results. Moreover, BO incorporates uncertainty explicitly into its decision-making, making it robust to noisy evaluations and well-suited to scenarios where the search space includes categorical or integer-valued variables \cite{garrido2020dealing}. Its application to hyperparameter tuning has demonstrated substantial improvements over traditional methods, solidifying its role as a state-of-the-art optimization technique \cite{liu2023bayesian}.

Performing a grid search to tune Genetic Algorithm (GA) parameters involves testing all possible parameter combinations. Given the ranges for the parameters: $3\times 3\times 3\times 3\times 2\times 3\times 4 = 1944$ combinations need to be tested for each instance set. Assuming each test takes approximately one minute to complete (although empirical results suggest that this could extend to an hour or more depending on the parameter values), the total computational time required exceeds 32 hours per instance set. This makes an exhaustive grid search computationally prohibitive, particularly for large problem instances or multiple instance sets. To address this limitation, we employed the Bayesian Optimization Python package to fine-tune the hyperparameters of the GA \cite{Fernando2014}. The following parameters and their respective value ranges were tested during the optimization process:

\begin{itemize}
    \item Population sizes: Integer values ranging from 100 to 300 
    \item Crossover rates: Real values between 0.5 and 0.9
    \item Mutation rates: Real values between 0.01 and 0.1
    \item Generations: Integer values ranging from 100 to 300 
    \item Selection strategies: Tournament selection, Roulette wheel
    \item Crossover types: PMX, Order, Order Based
    \item Mutation strategies: Swap, Insertion, Scramble, Inversion
\end{itemize}

This parameter space design ensured a comprehensive exploration of both continuous and discrete variables critical to GA performance. Bayesian Optimization was chosen due to its ability to efficiently navigate high-dimensional and potentially non-convex search spaces, offering a robust mechanism to optimize hyperparameters while minimizing the number of evaluations required \cite{garrido2020dealing}. Table \ref{tab_GA_Tuning} presents the results of the Bayesian Optimization-based fine-tuning process for instance number 1, across different patient sets. As observed, the total computational time required to obtain these results poses a challenge for extending the same method to other instances. Consequently, we will adopt the final average and dominant hyperparameter values in the solution algorithm to enhance computational efficiency. Therefore, for all problem instances, the following GA parameters are fixed: the population size is set to 254, the crossover rate to 0.64, the mutation rate to 0.07, and the number of generations to 185. Tournament selection is used as the parent selection method throughout. Unlike these static parameters, the crossover operator and mutation strategy are dynamically selected at each iteration by the RL module, enabling adaptive control based on ongoing performance feedback.

The Table \ref{tab_GA_Tuning} summarizes the best-performing hyperparameter configurations identified from 100 evaluated parameter sets. Each row corresponds to a configuration that achieved the best objective function value, with execution time serving as a tie-breaker when multiple configurations yielded equivalent objective values. As expected, the computational effort required for hyperparameter tuning increases significantly with problem size. While Bayesian Optimization effectively reduces the number of required evaluations compared to exhaustive search methods, the total convergence time remains substantial for large-scale instances. To further accelerate the solution methodology, additional techniques such as adaptive exploration strategies may be necessary.

\begin{table}[h]
\caption{Results for GA hyperparameter fine tuning}
\resizebox{\textwidth}{!}{%
\begin{tabular}{cccccccccc}
\hline
P & Population Size & Crossover Rate & Mutation Rate & Generations & Selection Strategy & Crossover Type & Mutation Strategy & Best Objective Value & Total time (Hrs) \\ \hline
40 & 208 & 0.68 & 0.01 & 151 & Roulette wheel & PMX & Swap & 10166.7 & 5.34 \\
60 & 209 & 0.50 & 0.06 & 100 & Tournament selection & Order & Inversion & 18125 & 7.82 \\
80 & 263 & 0.54 & 0.09 & 122 & Tournament selection & Order\_Based & Scramble & 26125 & 8.76 \\
100 & 218 & 0.50 & 0.07 & 175 & Tournament selection & Order & Inversion & 33708.3 & 12.41 \\
120 & 300 & 0.50 & 0.10 & 100 & Roulette wheel & Order\_Based & Inversion & 42250 & 14.18 \\
200 & 300 & 0.50 & 0.10 & 300 & Tournament selection & PMX & Inversion & 70875 & 15.76 \\
240 & 300 & 0.50 & 0.10 & 100 & Tournament selection & Order & Inversion & 87083.3 & 28.72 \\
300 & 299 & 0.81 & 0.08 & 283 & Tournament selection & PMX & Inversion & 107875 & 35.89 \\
360 & 300 & 0.88 & 0.10 & 300 & Tournament selection & Order & Inversion & 128000 & 45.07 \\
400 & 100 & 0.90 & 0.10 & 300 & Roulette wheel & PMX & Scramble & 139042 & 39.91 \\
480 & 300 & 0.76 & 0.02 & 100 & Roulette wheel & Order & Scramble & 168042 & 39.31 \\ \hline
Average & 254 & 0.64 & 0.07 & 185 &  &  &  &  & 23.02 \\
Dominant Parameters &  &  &  &  & Tournament   selection & Order & Inversion &  &  \\ \hline
\end{tabular}%
}
\label{tab_GA_Tuning}
\end{table}

\subsubsection{RL parameters selection and fine tuning}\label{RL_Parameters_Fine_Tuning}
The values for the reward and penalty parameters used in the RLA were selected based on both intuitive reasoning and fine tuning. Specifically, the $\text{L\_Reward}$ parameter, set to 1.0, provides a substantial incentive when a crossover or mutation operation results in an improvement in the cost function. This aligns with the primary objective of the algorithm and should be strongly reinforced. To determine the best value of the $\text{F\_Reward}$ parameter, we conducted an extensive sensitivity analysis of the $\text{F\_Reward}$ parameter over the set $\{0.05, 0.1, 0.25, 0.5\}$. For this analysis, all GA parameters were fixed at the values obtained from the GA hyperparameter fine-tuning stage. Each instance, both synthetic and real world (Naples Hospital), was executed three times to ensure robustness against stochastic effects. The results, summarized in Tables \ref{RL_Fine_tuning_Naple_Instances} and \ref{RL_Fine_tuning_Synthetic_Instan}, show that $\text{F\_Reward} = 0.25$ consistently yielded the best performance, with this setting achieving the lowest values in the majority of cases. These findings confirm that the chosen value of $\text{F\_Reward} = 0.25$ effectively balances the exploration of feasible regions with cost minimization, and that the proposed approach remains robust across different classes of problem instances. The $\text{R\_Penalty}$ parameter, assigned a modest value of 0.1, was chosen to avoid prematurely discouraging exploration, especially in early iterations where infeasible solutions might provide useful search directions. This reward–penalty balance encourages diversification while guiding the search toward high-quality feasible solutions. These values were validated through empirical observation and are consistent with practices in related adaptive evolutionary algorithms, where heavier rewards are typically used for primary objectives and smaller penalties are applied to maintain a balance between intensification and diversification \cite{sonucc2023adaptive, li2023differential}.

\begin{table}[h]
\caption{Results for F\_Reward fine tuning for real-world (Naples hospital) instances}
\centering
\tiny
\begin{tabular}{cccccccc}
\hline
\multirow{2}{*}{\begin{tabular}[c]{@{}c@{}}Instance\\ name\end{tabular}} & \multirow{2}{*}{$\lvert P \lvert$} & \multirow{2}{*}{$\lvert D \lvert$} &  & \multicolumn{4}{c}{$\text{F\_Reward}$} \\ \cline{5-8} 
 &  &  &  & 0.05 & 0.1 & 0.25 & 0.5 \\ \hline
INST31 & 282 & 12 &  & 53986.1 & 53625 & 54069.43 & \textbf{53611.13} \\
INST32 & 289 & 12 &  & 54194.43 & 54527.77 & \textbf{54013.87} & 54222.23 \\
INST33 & 306 & 12 &  & 56041.67 & \textbf{55916.7} & 56597.23 & 56500 \\
INST34 & 298 & 12 &  & \textbf{57208.33} & 57791.7 & 57444.47 & 58291.67 \\
INST35 & 319 & 12 &  & \textbf{61611.1} & 62694.43 & 62444.47 & 61638.9 \\
INST36 & 285 & 12 &  & 49958.37 & 50000 & 50277.8 & \textbf{49722.2} \\
INST37 & 262 & 12 &  & 48083.33 & \textbf{47888.9} & 48263.87 & 48166.67 \\
INST38 & 278 & 12 &  & 51250.03 & 51125 & \textbf{50847.23} & 51138.87 \\
INST39 & 319 & 12 &  & 64861.1 & 65069.43 & 64569.43 & \textbf{63972.2} \\
INST40 & 296 & 12 &  & 53916.67 & \textbf{53736.1} & 54111.13 & 54847.2 \\
INST41 & 259 & 12 &  & 48958.3 & 49444.43 & \textbf{48583.33} & 49027.8 \\
INST42 & 258 & 12 &  & 45638.9 & \textbf{45166.67} & 45458.33 & 45375 \\
INST43 & 261 & 12 &  & - & - & - & - \\
INST44 & 274 & 12 &  & 47986.13 & 47736.1 & \textbf{47541.67} & 47888.87 \\
INST45 & 280 & 12 &  & 53111.1 & 53277.77 & 53819.43 & \textbf{52972.2} \\
INST46 & 446 & 18 &  & 85083.33 & 85013.9 & \textbf{84986.1} & 86027.77 \\
INST47 & 431 & 18 &  & 83333.37 & \textbf{83097.23} & 83194.43 & 83111.1 \\
INST48 & 462 & 18 &  & 89305.57 & 89444.43 & \textbf{88791.63} & 89236.13 \\
INST49 & 440 & 18 &  & 92138.9 & 92444.43 & 91694.43 & \textbf{91000} \\
INST50 & 399 & 18 &  & 84513.9 & \textbf{82875.03} & 83611.1 & 83194.43 \\
INST51 & 460 & 18 &  & 73861.1 & 73333.33 & \textbf{72847.2} & 72916.67 \\
INST52 & 406 & 18 &  & 74041.67 & 73513.9 & \textbf{73402.77} & 73555.57 \\
INST53 & 407 & 18 &  & 84444.43 & 83375 & 83666.67 & \textbf{83111.1} \\
INST54 & 403 & 18 &  & 95305.57 & 95125 & 94861.13 & \textbf{94305.57} \\
INST55 & 412 & 18 &  & 75777.77 & 75625 & \textbf{75000} & 75208.33 \\ \hline
\multicolumn{2}{c}{\# Times selected as minimum} &  &  & 2 & 6 & \textbf{9} & 7 \\ \hline
\end{tabular}%
\label{RL_Fine_tuning_Naple_Instances}
\end{table}

\begin{table}[h]
\caption{Results for F\_Reward fine tuning for synthetic instances}
\centering
\tiny
\begin{tabular}{cccccc}
\hline
\multirow{2}{*}{$\lvert P \lvert$} & \multirow{2}{*}{\# Instance} & \multicolumn{4}{c}{$\text{F\_Reward}$} \\ \cline{3-6} 
 &  & 0.05 & 0.1 & 0.25 & 0.5 \\ \hline
\multirow{5}{*}{200} & 1 & 65611.13 & 66250 & \textbf{65458.33} & 66513.9 \\
 & 2 & \textbf{66486.1} & 66625 & 66930.57 & 66888.87 \\
 & 3 & 61791.7 & 61833.33 & \textbf{61097.2} & 61500 \\
 & 4 & \textbf{58194.47} & 58486.1 & \textbf{58305.53} & 58750 \\
 & 5 & 63000.03 & \textbf{62555.57} & 62750 & 62583.37 \\
 \hline
\multirow{5}{*}{240} & 1 & 81430.53 & 82361.1 & 81097.23 & \textbf{80625.03} \\
 & 2 & 79319.47 & 79305.53 & \textbf{79152.77} & 79333.3 \\
 & 3 & \textbf{76500} & 76861.1 & 77847.23 & 77333.37 \\
 & 4 & 78430.57 & 78694.43 & \textbf{78388.87} & 78458.33 \\
 & 5 & 77625 & 77333.37 & 77513.9 & \textbf{77319.43} \\
 \hline
\multirow{5}{*}{300} & 1 & 100764 & 101250 & \textbf{100472.3} & 101458.3 \\
 & 2 & \textbf{100708.2} & 101013.8 & 101180.7 & 101819.3 \\
 & 3 & 101972 & \textbf{101541.7} & 101972 & 102250 \\
 & 4 & 97375 & 97819.43 & \textbf{96986.13} & 97736.1 \\
 & 5 & 101389 & 101097.3 & \textbf{100194.3} & 101597.3 \\
 \hline
\multirow{5}{*}{360} & 1 & 118458.3 & 118500 & \textbf{117917} & 118875 \\
 & 2 & 123055.3 & 121916.7 & 122111.3 & \textbf{121861.3} \\
 & 3 & 128402.7 & \textbf{128180.7} & 129514 & 128361.3 \\
 & 4 & 117486 & 116708 & 117194.3 & \textbf{116194.3} \\
 & 5 & 117069.7 & 117194.7 & 118777.7 & \textbf{115472.3} \\
 \hline
\multirow{5}{*}{400} & 1 & \textbf{129680.7} & 130416.7 & 130805.7 & 131944.3 \\
 & 2 & 135486.3 & \textbf{133916.7} & 134750 & 134833.3 \\
 & 3 & \textbf{130736} & 131680.7 & 131680.7 & 131583 \\
 & 4 & 132736.3 & \textbf{132139} & 134416.7 & 133416.7 \\
 & 5 & 136625 & 138208.3 & 136361 & \textbf{136250} \\
 \hline
\multirow{5}{*}{480} & 1 & 159486 & 159375 & \textbf{158264} & 160736.3 \\
 & 2 & 159541.7 & \textbf{159305.7} & 160944.3 & 162000 \\
 & 3 & 162347.3 & 161416.7 & \textbf{161055.3} & 162250 \\
 & 4 & 162417 & 159916.5 & 159722.3 & \textbf{158083.3} \\
 & 5 & 161875 & 160611 & \textbf{157902.7} & 159694.7 \\ \hline
\multicolumn{2}{c}{\# Times selected as minimum} & 6 & 6 & \textbf{11} & 7 \\ \hline
\end{tabular}%
\label{RL_Fine_tuning_Synthetic_Instan}
\end{table}

\subsubsection{BO hyperparameters fine tuning}\label{BO_Parameters_Fine_Tuning}
We performed a full factorial study over three BO hyperparameters following Fernando \cite{Fernando2014}: acquisition function (acq), number of random initial points ($init\_points$) and number of BO iterations ($n\_iter$), yielding 27 configurations per instance. Generalizability was assessed on twelve instances spanning three scales (small/medium/large) and two datasets (synthetic and real-world), with two distinct instances per scale. As summarized in Table \ref{BO_hyperparameter_tuning}, UCB attained the best mean performance in 6/12 instances, $init\_points = 5$ was most frequently selected (5/12), and $n\_iter$ exhibited a near tie between 50 and 100 (5/12 each). Therefore, the acquisition choice had a larger impact on performance than $init\_points$ or $n\_iter$. We verified that the qualitative conclusions of the GA hyperparameter study remain unchanged under these BO settings since the optimal values are the same as their default values, except for $n\_iter$ where the default value is 50. Furthermore, repeating the key experiments with $n\_iter = 100$ did not alter any conclusions related to Table \ref{tab_GA_Tuning}. Overall, this sensitivity study supports our original choice of $acq = UCB$ and $init\_points = 5$, and shows that the method is robust across small, large, and real-world instances.

\begin{table}[h]
\caption{Results for Bayesian Optimization hyperparameter fine tuning}
\centering
\tiny
\resizebox{\textwidth}{!}{%
\begin{tabular}{ccclccclccclccclc}
\hline
\multirow{2}{*}{Type of instance} & \multirow{2}{*}{Size\_Class} & \multirow{2}{*}{$\lvert P \lvert$ (name)} &  & \multicolumn{3}{c}{$acq$} &  & \multicolumn{3}{c}{$init\_points$} &  & \multicolumn{3}{c}{$n\_iter$} &  & \multirow{2}{*}{Best Objective Value} \\ \cline{5-7} \cline{9-11} \cline{13-15}
 &  &  &  & UCB & EI & POI &  & 5 & 10 & 20 &  & 25 & 50 & 100 &  &  \\ \hline
\multirow{6}{*}{Synthetic Instances} & \multirow{2}{*}{small} & 60 &  &  & \checkmark &  &  &  & \checkmark &  &  & \checkmark &  &  &  & 18542 \\
 &  & 120 &  &  &  & \checkmark &  & \checkmark &  &  &  &  & \checkmark &  &  & 43011 \\
 & \multirow{2}{*}{Medium} & 200 &  & \checkmark &  &  &  &  & \checkmark &  &  &  &  & \checkmark &  & 72434 \\
 &  & 240 &  &  &  & \checkmark &  &  &  & \checkmark &  &  &  & \checkmark &  & 88564 \\
 & \multirow{2}{*}{Large} & 400 &  & \checkmark &  &  &  & \checkmark &  &  &  &  & \checkmark &  &  & 142240 \\
 &  & 480 &  & \checkmark &  &  &  & \checkmark &  &  &  &  &  & \checkmark &  & 172411 \\ \hline
\multirow{6}{*}{Naples   Instances} & \multirow{2}{*}{small} & 110 (INST21) &  &  &  & \checkmark &  &  &  & \checkmark &  & \checkmark &  &  &  & 14817 \\
 &  & 162 (INST18) &  &  & \checkmark &  &  & \checkmark &  &  &  &  & \checkmark &  &  & 23445 \\
 & \multirow{2}{*}{Medium} & 259 (INST41) &  & \checkmark &  &  &  &  & \checkmark &  &  &  &  & \checkmark &  & 38836 \\
 &  & 258 (INST42) &  &  &  & \checkmark &  & \checkmark &  &  &  &  & \checkmark &  &  & 45494 \\
 & \multirow{2}{*}{Large} & 460 (INST51) &  & \checkmark &  &  &  &  &  & \checkmark &  &  & \checkmark &  &  & 55747 \\
 &  & 412 (INST55) &  & \checkmark &  &  &  &  & \checkmark &  &  &  &  & \checkmark &  & 59870 \\ \hline
\multicolumn{3}{c}{\# times picked as   the best value} &  & \textbf{6} & 2 & 4 &  & \textbf{5} & 4 & 3 &  & 2 & \textbf{5} & \textbf{5} &  &  \\ \hline
\end{tabular}%
}
\label{BO_hyperparameter_tuning}
\end{table}

\subsection{Numerical Results}\label{Results}
In this section, we present a detailed analysis of the two sets of synthetic instances and the real-world instances, as described in Sections \ref{Definition_Instances_Synthetic} and \ref{Definition_Instances_Real}, respectively. We set all algorithms' time limits to one hour for small and medium-sized instances and two hours for large instances to obtain the results in a reasonable time. For CG, we set the time limit for the fast feasible column detection technique to 10 seconds. When CG stopping criteria are met, we solve the M-CG for a 30-minute time limit and report the upper bound found by CG as a final solution to the problem.
The time limit for the RGA is set to one hour. Therefore, from the total two hours of time limitation of the Hybrid RGA-CG, at most one hour is dedicated to the RGA. At the end of this time limit, we use the best-found solution of RGA as an initial column to start the CG algorithm.

\subsubsection{Computational Results for Synthetic Instances}\label{Results_Synthetic}
In Tables \ref{tab_6} to \ref{tab_8}, each row shows the average of five instances generated for each category. The "IS\%" column shows the percentage of instances for which a feasible solution is found. The "OSI\%" column that stands for Optimally Solved Instances, gives the percentage of instances solved to optimality. Table \ref{tab_6} shows the computational results of the small and medium-sized instances for the MIORPS, and PMIORPS. In the Table \ref{tab_6}, according to the "IS\%" column, only five instances achieved $Gap\% < 100$ in the MIORPS model.

\begin{sidewaystable}
\doublespacing
\caption{Evaluation of the MIORPS model, and PMIORPS model}
\resizebox{\textwidth}{!}{%
\tiny
\begin{tabular}{ccccccccccccccc}
\hline
\multirow{2}{*}{\begin{tabular}[c]{@{}c@{}}$\lvert P \lvert$\end{tabular}} &  & \multicolumn{6}{c}{MIORPS} &  & \multicolumn{6}{c}{PMIORPS} \\ \cline{3-8} \cline{10-15} 
 &  & IS\% & OSI\% & LB & UB & Gap\% & Time &  & IS\% & OSI\% & LB & UB & Gap\% & Time \\ \hline
40 &  & 80 & 0 & 7086.63 & 16308.34 & 6.53 & 3600 &  & 100 & 40 & 8738.91 & 8800 & \textbf{0.67} & 2369.42 \\
60 &  & 20 & 0 & 6507.56 & 44150 & 28.84 & 3600 &  & 100 & 20 & 14206.52 & 14400 & \textbf{1.37} & 2900 \\
80 &  & 0 & 0 & 1317.67 & 61400.02 & - & 3600 &  & 100 & 0 & 19341.66 & 19666.66 & \textbf{1.69} & 3600 \\
100 &  & 0 & 0 & 0 & 69783.34 & - & 3600 &  & 100 & 0 & 23533.34 & 23700 & \textbf{0.7} & 3600 \\
120 &  & 0 & 0 & 0 & 75433.32 & - & 3600 &  & 100 & 0 & 29737.5 & 30250 & \textbf{1.71} & 3600 \\ \hline
Ave. &  & 20 & 0 & - & - & 17.69 & 3600 &  & \textbf{100} & \textbf{12} & - & - & \textbf{1.23} & \textbf{3213.89} \\ \hline
\end{tabular}%
}
\label{tab_6}

\vspace{2cm}
\hspace{1cm}

\caption{Evaluation of the RGA-CG and PMIORPS model}
\resizebox{\textwidth}{!}{%
\begin{tabular}{cccccccccccccccccc}
\hline
\multirow{2}{*}{\begin{tabular}[c]{@{}c@{}}$\lvert P \lvert$\end{tabular}} & \multirow{2}{*}{$\lvert D \lvert$} &  & \multicolumn{6}{c}{RGA\_CG} &  & \multicolumn{6}{c}{BCIORPS} &  & \multirow{2}{*}{\begin{tabular}[c]{@{}c@{}}\# Times CG \\ dominated BC\end{tabular}} \\ \cline{3-9} \cline{11-16}
 &  &  & \# Gen-Col & CG-Time & IP-Time & Time & UB & Time &  & \# Nodes & \begin{tabular}[c]{@{}c@{}}\# Feasibility\\ Cuts\end{tabular} & \begin{tabular}[c]{@{}c@{}}\# Optimality\\ Cuts\end{tabular} & \# Solved & UB & \# Solved &  &  \\ \hline
200 & 10 &  & 2965.6 & 7200 & 45.324 & 7245.324 & 63125 & 5 &  & 0 & 0 & 5.2 & 7200 & 71166.7 & 1 &  & \textbf{5} \\
240 & 10 &  & 2717.6 & 7200 & 21.156 & 7221.156 & 79383.34 & 5 &  & 0 & 0 & 2.8 & 7200 & INF & 0 &  & \textbf{5} \\
300 & 15 &  & 3708.2 & 7200 & 82.52 & 7282.52 & 101541.46 & 5 &  & 0 & 0 & 5.8 & 7200 & INF & 0 &  & \textbf{5} \\
360 & 15 &  & 3749 & 7200 & 56.712 & 7256.712 & 121749.8 & 5 &  & 0 & 0 & 7.8 & 7200 & INF & 0 &  & \textbf{5} \\
400 & 20 &  & 5710 & 7200 & 176.152 & 7376.152 & 134266.8 & 5 &  & 0 & 0 & 5.2 & 7200 & INF & 0 &  & \textbf{5} \\
480 & 20 &  & 7050.6 & 7200 & 43.418 & 7243.418 & 158841.8 & 5 &  & 0 & 0 & 10.6 & 7200 & INF & 0 &  & \textbf{5} \\ \hline
Ave. & - &  & 4316.8 & 7200 & 70.88 & 7270.88 & - & \textbf{30} &  & 0 & 0 & - & 7200 & - & 1 &  & \textbf{30} \\ \hline
\end{tabular}%
}
\label{tab_7}

\end{sidewaystable}

The key observations based on the data presented in Table \ref{tab_6} are as follows:
\begin{itemize}
    \item MIORPS demonstrates inadequacy as a solution strategy for our integrated operating room planning and scheduling problem. It consistently fails to achieve reasonable solution gaps (i.e., $<100\%$) within the designated computational timeframe for most instances.
    \item We observed the best remarkable performance with our proposed PMIORPS model. This model not only maintains acceptable solution gaps across all instances but also drastically reduces the average gap to 1.23\%. Such a reduction highlights a significant enhancement in solution quality compared to MIORPS.
\end{itemize}

So far, all the results have proven the superiority of our proposed model. However, we went a step further and investigated the performance of the PMIORPS against problems with a large number of surgeries and extended planning horizons. In fact, by extending the planning horizon, we are testing the applicability of our model at the tactical level of planning and scheduling. For this purpose, we implemented the Hybrid RGA-CG described in Sections \ref{RGA_CG} and \ref{Hybrid_RL_GA}. Table \ref{tab_7} present the computational results of our proposed RGA-CG and an existing branch-and-cut algorithm from the literature \cite{roshanaei2021solving} that was proposed for a fairly similar problem. In Appendix \ref{AppendixC}, we have explained how we have modified the original branch-and-cut algorithm for solving our problem. According to Table \ref{tab_7}, the branch-and-cut algorithm from the literature fails to find a valid upper bound for all instances except one, while RGA-CG shows its capability to find a feasible solution for all of the instances. For the only instance solved by the branch-and-cut algorithm, RGA-CG also provided a better upper bound, improving it by 8.03\%.

The "$\# Gen-Col$" column shows the average number of columns generated by the Hybrid RGA-CG. Also, the "$\# Solved$" column shows the number of instances solved with a valid UB as a feasible solution to the problem. While Hybrid RGA-CG can find the valid UB or a feasible solution in all instances, BCIORPS only found a valid feasible solution in just one of the cases. In fact, due to the size of instances, BCIORPS could not get out of the root node as indicated by the "$\# Nodes$" column. In all instances, Hybrid RGA-CG outperformed the BCIORPS and proved that it could obtain a feasible solution for large instances while the others failed. 

\begin{figure}[H]
    \centering
    \begin{subfigure}{0.4\textwidth}
        \includegraphics[width=0.95\textwidth, height=1.5in]{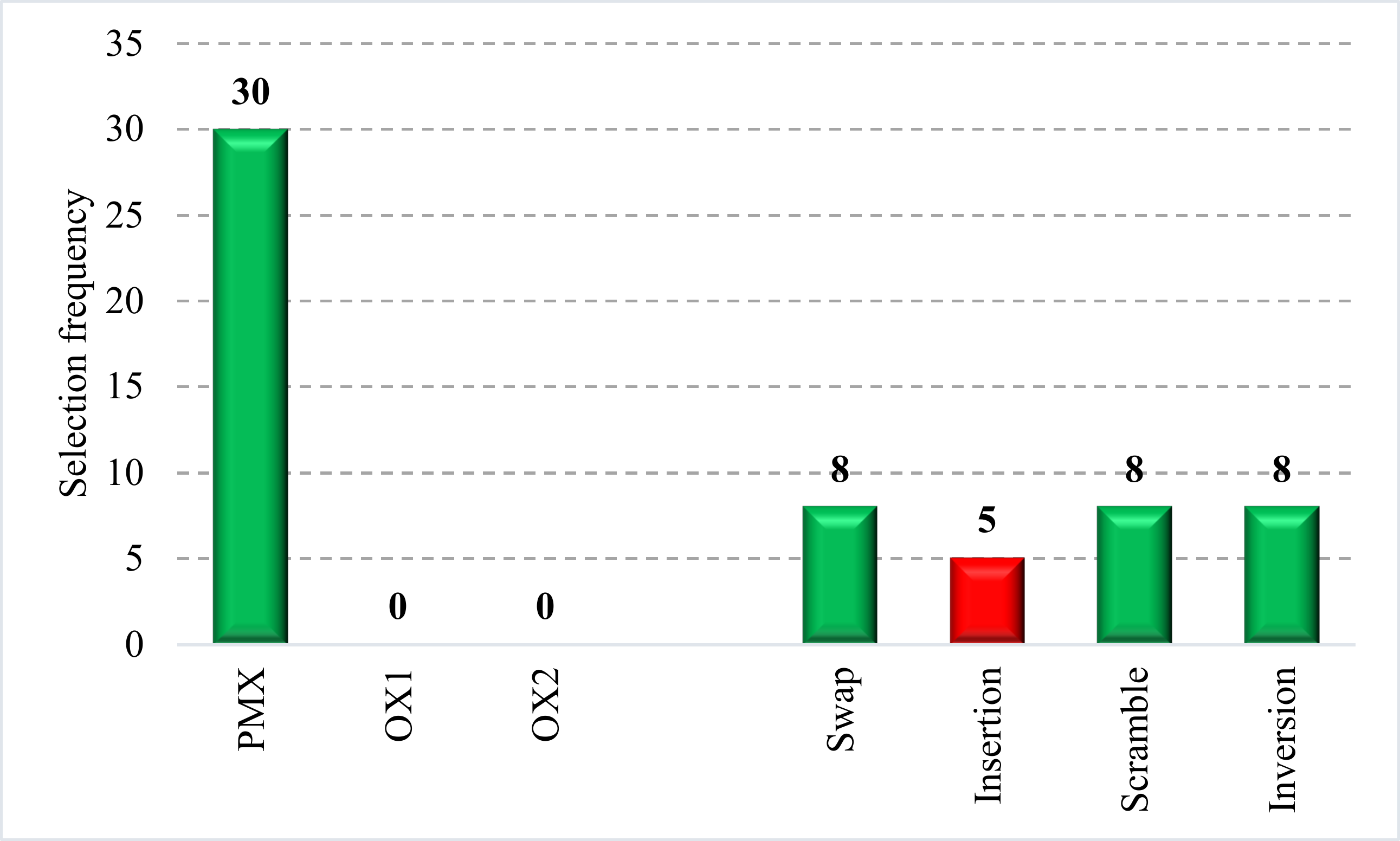}
        \caption{\label{fig_7a}}
    \end{subfigure}
    \begin{subfigure}{0.4\textwidth}
        \includegraphics[width=0.95\textwidth, height=1.5in]{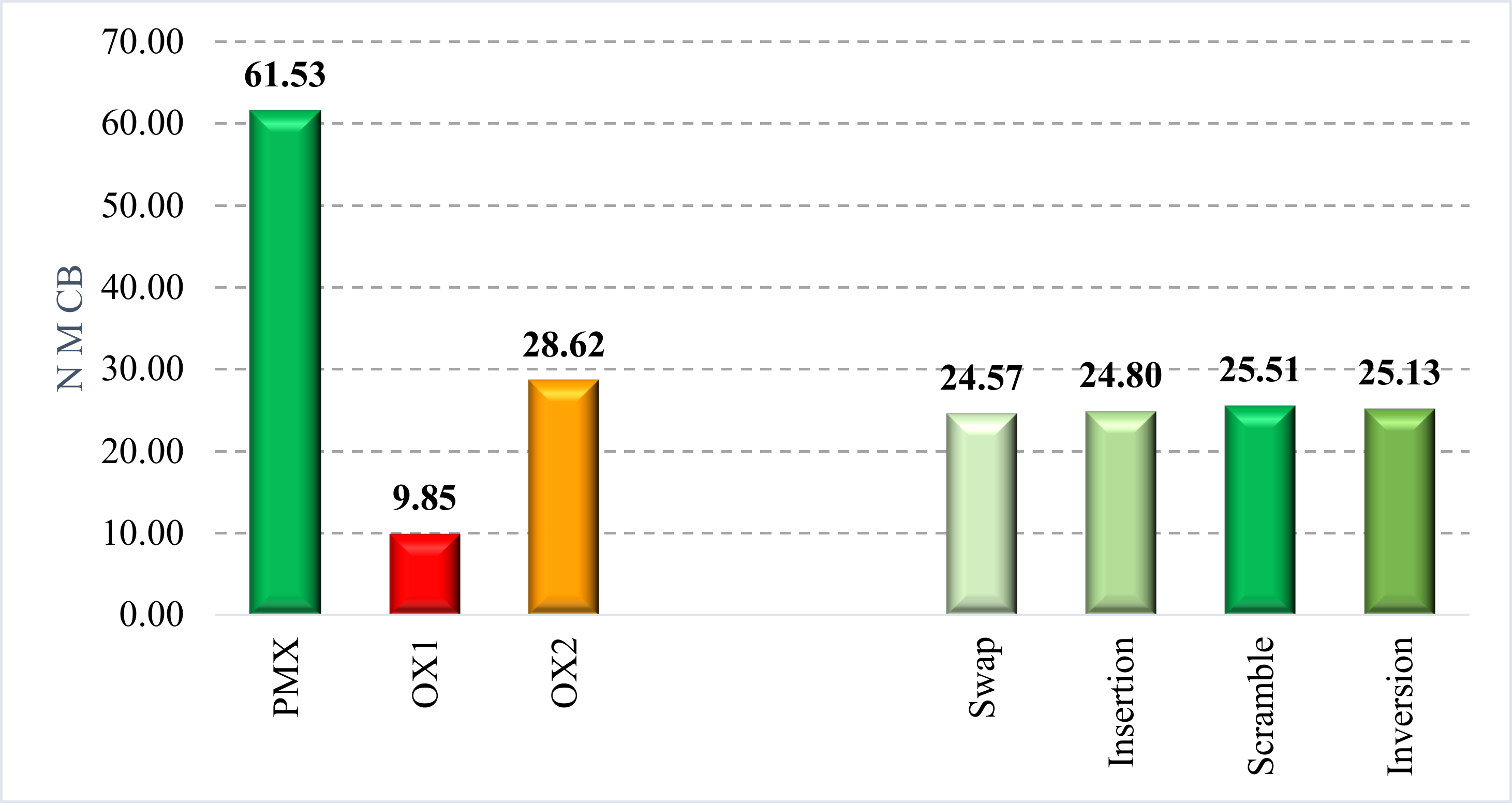}
        \caption{\label{fig_7b}}
    \end{subfigure}
    \caption{(\subref{fig_7a}) GA operator selection's frequency. (\subref{fig_7b}) Normalized $MCB_{Op}$}
    \label{fig_7}
\end{figure}

One of the contributions of this study is to integrate RLA and GA. This integration enables dynamic selection of genetic operators during the search process. Figure \ref{fig_7a} shows PMX crossover was consistently selected as the best crossover method in all solved instances reported in Table \ref{tab_7}. This outcome contrasts with the results from Table \ref{tab_GA_Tuning}, where the Order or OX1 Crossover operator was previously identified as the most effective strategy. This discrepancy highlights the value of online operator selection, which is made possible through the use of the RGA. Regarding mutation strategies, Swap, Scramble, and Inversion operators demonstrated similar performance, each being selected eight times. In contrast, the Insertion operator was chosen only five times, suggesting that it was less effective in the evaluated scenarios. These results reinforce the importance of adaptive operator selection to enhance the performance of GA across different instances.

As previously discussed, Figure \ref{fig_7a} presents the final selection outcomes for each instance. In order to have a clearer understanding of the details of the operators' selection procedure, we can use the normalized mean confidence bound ($NMCB_{Op}$) value that is computed as follows:

\begin{flalign}\label{Eq_45}
 NMCB_{Op} = \frac{MCB_{Op}}{\sum_{Op' \in Each\; Category \;of \;Operators} MCB_{Op'}} &&
\end{flalign}

The $NMCB$ indicates the probability of selecting each operator across all iterations and instances. It offers more detailed insight by capturing the selection behavior of the GA at each iteration for every instance. To compute $NMCB$,  we first calculate the mean confidence bound for each operator, denoted as $MCB_{Op} = \frac{TR_{Op}}{TOS_{Op}}$. This quantity serves as the initial component in the computation of the $UCB_{op}$ as defined in Equation \ref{Eq_44}. Since the total summation of $NMCB_{Op}$ for all operators belong to crossover and mutation categories are not necessarily equal to 1,we normalize these values using Equation \ref{Eq_45} to make sure that they have common scale so they can be comparable to each other. 

Figure \ref{fig_7b} shows that the average confidence level for PMX is clearly more than that of OX1 and OX2. Therefore, with high confidence, we could say that PMX is the best crossover strategy for the developed GA in this study. However, comparing $NMCB$ for mutation operators reveals that all four methods could be considered equal. 

To evaluate the performance of the RGA-CG algorithm, we implemented the four heuristic methods described by Lin and Chou \cite{lin2020hybrid}, namely LPT+, SPT+, EDD+, and LPT/EDD+. LPT rule schedules surgeries in decreasing order of their expected durations. In contrast, the SPT rule prioritizes surgeries with the shortest durations. The EDD rule orders surgeries based on their respective deadlines, minimizing maximum lateness and reducing the risk of delaying time-sensitive procedures. This is particularly relevant when clinical urgency dictates strict completion windows. A hybrid approach, the LPT/EDD priority rule, combines elements of both LPT and EDD by selecting longer surgeries among those with the earliest due dates. This integrated strategy balances the goals of efficient room utilization and timely treatment, making it suitable for overloaded systems that must manage both extended procedures and urgent clinical deadlines.

As shown in Table \ref{tab_8}, the RGA-CG consistently yields superior upper bounds across all problem sizes. Specifically, for each instance size, the upper bound obtained by RGA-CG is lower than those generated by the four heuristics. On average, RGA-CG achieves an upper bound of 109,918, outperforming the next-best heuristic, SPT+, by a significant margin of 10,562.5. This result clearly demonstrates the effectiveness of the RGA-CG algorithm in providing high-quality solutions for all tested instance sizes.

\begin{table}[h]
\caption{Comparison between various heuristic algorithms and RGA-CG}
\centering
\small
\begin{tabular}{cccccc}
\hline
\multirow{2}{*}{$\lvert P \lvert$} & \multicolumn{5}{c}{{UB}} \\ \cline{2-6} 
 & LPT+ & SPT+ & EDD+ & LPT/EDD+ & RGA CG \\ \hline
200 & 70525 & 70250 & 70383 & 70425.00 & 63125 \\
240 & 88525 & 88433 & 88408 & 89108 & 79383 \\
300 & 111075 & 111200 & 111075 & 110600 & 101541 \\
360 & 134175 & 133717 & 133542 & 134167 & 121750 \\
400 & 144517 & 144275 & 144225 & 144500 & 134267 \\
480 & 175142 & 174408 & 174883 & 174692 & 158842 \\ \hline
Ave. & 120660 & 120381 & 120420 & 120582 & 109818 \\ \hline
\end{tabular}%
\label{tab_8}
\end{table}

\vspace{-0.25cm}

\subsubsection{Computational Results for Real-world Instances}\label{Results_Real}
Tables \ref{tab_real_1} and \ref{tab_real_2} present the computational results for the real-world instances. As shown in Table \ref{tab_real_1}, the baseline MIORPS model was only able to find solutions with a GAP\% below 100\% for two instances (Instance 4 and Instance 8). Among the remaining 28 instances, half could not be solved at all, while the other half yielded solutions with GAP values exceeding 100\%. In contrast, the proposed PMIORPS model successfully solved all 30 instances, achieving an average optimality gap of  1.49\%. Moreover, PMIORPS reached the optimal solution in five instances, clearly demonstrating its superiority over MIORPS in terms of both solution quality and robustness.
A comparison between the results in Table \ref{tab_real_1} and those in Table \ref{tab_6} also confirms the consistent performance of PMIORPS across different instance types. While the average optimality gap for real-world instances is slightly higher than that for synthetic ones, this is expected given the greater complexity of real-world cases with an average instance size of 142.2, compared to 80 in synthetic scenarios.

Table~\ref{tab_real_2} presents a comparative analysis between the RGA-CG and BCIORPS algorithms on real-world instances with an extended planning horizon. Out of the 25 instances evaluated, RGA-CG failed to solve only one instance (Instance 43) due to memory limitations. In contrast, BCIORPS was unable to find feasible solutions for three instances.
When comparing the upper bound values obtained by both algorithms, RGA-CG outperformed BCIORPS in 17 out of the 25 instances, with lower upper bounds in the majority of cases. These results demonstrate that the proposed hybrid RGA-CG algorithm not only exhibits greater robustness but also consistently delivers higher-quality solutions compared to BCIORPS, particularly under more complex and extended planning scenarios.

\begin{figure}[h]
    \centering
    \begin{subfigure}{0.4\textwidth}
        \includegraphics[width=0.95\textwidth, height=1.5in]{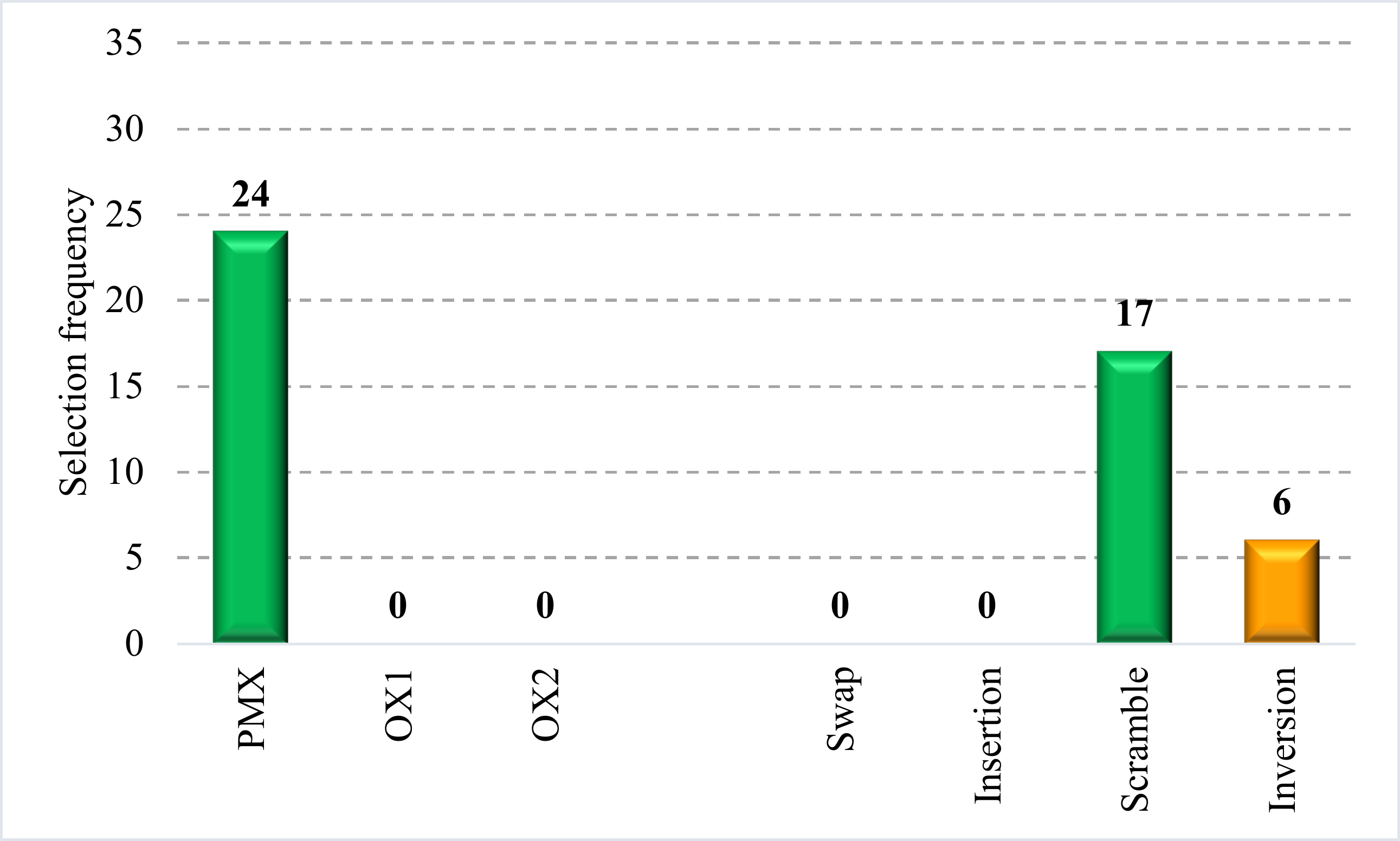}
        \caption{\label{fig_RL_reala}}
    \end{subfigure}
    \begin{subfigure}{0.4\textwidth}
        \includegraphics[width=0.95\textwidth, height=1.5in]{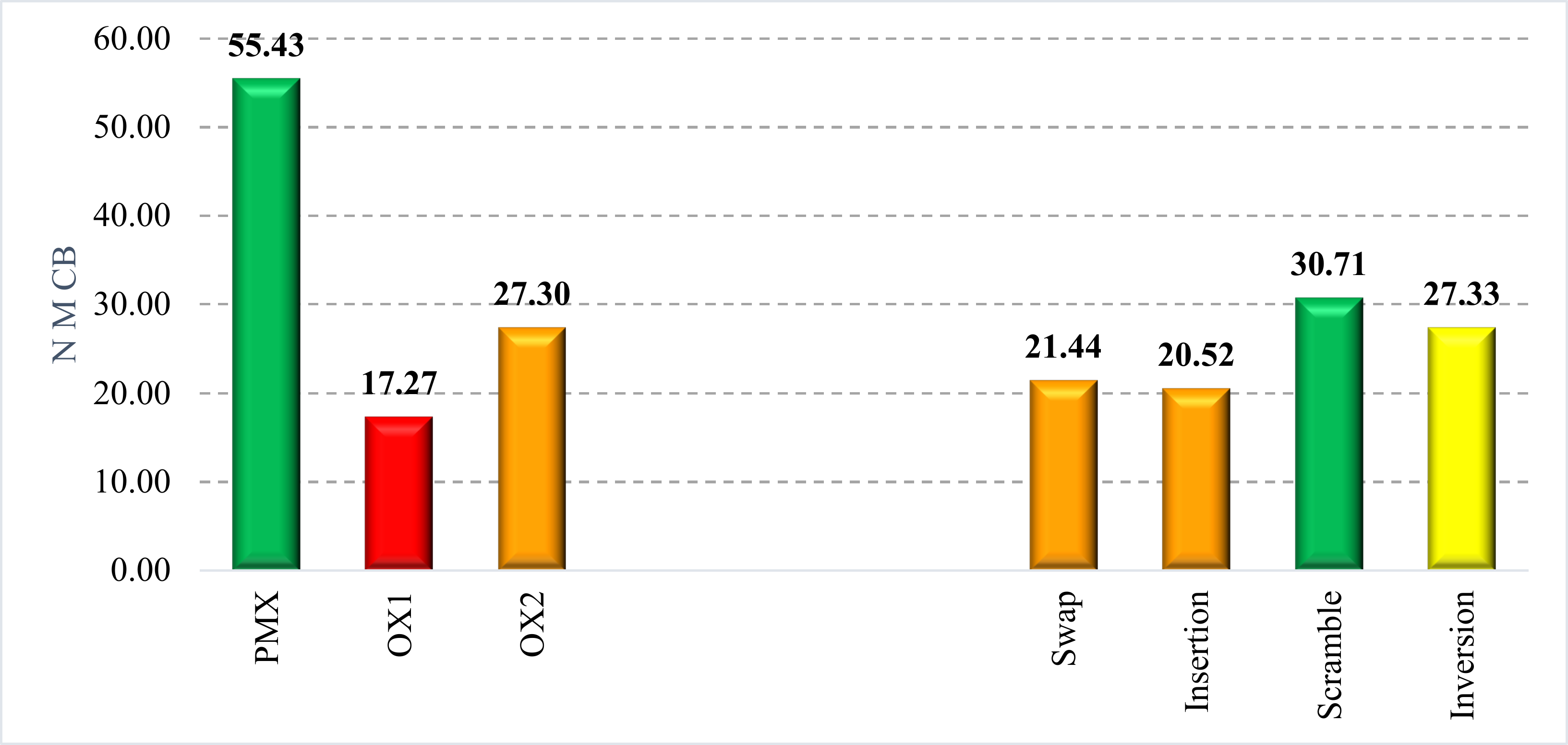}
        \caption{\label{fig_RL_realb}}
    \end{subfigure}
    \caption{(\subref{fig_7a}) GA operator selection's frequency. (\subref{fig_7b}) Normalized $MCB_{Op}$}
    \label{fig_RL_real}
\end{figure}

An analysis of the detailed RL results reveals that PMX was consistently selected as the most effective crossover operator across the real-world instances. For mutation operators, the Scramble strategy was chosen as the best-performing operator in 17 instances, while Inversion was selected in the remaining 6 instances. These findings are illustrated in Figure~\ref{fig_RL_real}.
When comparing Figure~\ref{fig_RL_real} with the earlier results shown in Figure~\ref{fig_7}, it is evident that PMX consistently outperforms the other crossover strategies, confirming its robustness across different types of instances. Among the mutation strategies, Scramble remains the top choice, followed by Inversion, which indicates the effectiveness of these two operators within the RL-guided GA framework.

\begin{sidewaystable}
\doublespacing
\caption{Evaluation of the MIORPS model, and PMIORPS model}
\resizebox{\textwidth}{!}{%
\tiny
\begin{tabular}{cccccccccccccccc}
\hline
\multirow{2}{*}{\begin{tabular}[c]{@{}c@{}}Instance\\ name\end{tabular}} & \multirow{2}{*}{$\lvert P \lvert$} &  & \multicolumn{6}{c}{MIORPS} &  & \multicolumn{6}{c}{PMIORPS} \\ \cline{4-9} \cline{11-16} 
 &  &  & IS & OSI & LB & UB & Gap\% & Time &  & IS & OSI & LB & UB & Gap\% & Time \\ \hline
INST1 & 143 &  & 0 & 0 & 0 & 93875 & INF & 7200 &  & 1 & 0 & 17775 & 18000 & 1.27 & 7200 \\
INST2 & 139 &  & 0 & 0 & 0 & 92291.7 & INF & 7200 &  & 1 & 1 & 18000 & 18000 & 0 & 2391 \\
INST3 & 164 &  & 0 & 0 & 0 & 104000 & INF & 7200 &  & 1 & 0 & 19708.3 & 20125 & 2.11 & 7200 \\
INST4 & 125 &  & 1 & 0 & 11696.9 & 22541.7 & 92.715 & 7200 &  & 1 & 1 & 16000 & 16000 & 0 & 1919.27 \\
INST5 & 160 &  & 1 & 0 & 500 & 103417 & 20583.3 & 7200 &  & 1 & 0 & 18500 & 18875 & 2.03 & 7201.9 \\
INST6 & 146 &  & 1 & 0 & 500 & 97541.7 & 19408.3 & 7200 &  & 1 & 0 & 18675 & 19000 & 1.74 & 7200 \\
INST7 & 162 &  & 0 & 0 & 0 & 103500 & INF & 7200 &  & 1 & 0 & 20941.7 & 21583.3 & 3.06 & 7200 \\
INST8 & 136 &  & 1 & 0 & 15359.2 & 27541.7 & 79.3166 & 7200 &  & 1 & 0 & 17991.7 & 18500 & 2.75 & 7200 \\
INST9 & 164 &  & 1 & 0 & 1000 & 103417 & 10241.7 & 7200 &  & 1 & 0 & 23033.3 & 23500 & 2.03 & 7200 \\
INST10 & 155 &  & 0 & 0 & 0 & 100083 & INF & 7200 &  & 1 & 0 & 19316.7 & 19500 & 0.95 & 7200 \\
INST11 & 136 &  & 1 & 0 & 7963.5 & 27833.3 & 249.512 & 7200 &  & 1 & 0 & 14758.3 & 15000 & 1.64 & 7200 \\
INST12 & 149 &  & 0 & 0 & 0 & 96958.3 & INF & 7200 &  & 1 & 0 & 18091.7 & 18500 & 2.26 & 7200 \\
INST13 & 129 &  & 1 & 0 & 7417.32 & 88000 & 1086.41 & 7200 &  & 1 & 0 & 14408.3 & 14500 & 0.64 & 7200 \\
INST14 & 133 &  & 1 & 0 & 1000 & 88708.3 & 8770.83 & 7200 &  & 1 & 0 & 17716.7 & 18000 & 1.6 & 7200 \\
INST15 & 137 &  & 0 & 0 & 0 & 91000 & INF & 7200 &  & 1 & 0 & 16308.3 & 16500 & 1.18 & 7200 \\
INST16 & 141 &  & 1 & 0 & 1000 & 92125 & 9112.5 & 7200 &  & 1 & 1 & 18541.7 & 18541.7 & 0 & 6143.16 \\
INST17 & 157 &  & 0 & 0 & 0 & 101167 & INF & 7200 &  & 1 & 0 & 22116.7 & 22500 & 1.73 & 7200 \\
INST18 & 162 &  & 0 & 0 & 0 & 103375 & INF & 7200 &  & 1 & 0 & 21566.7 & 22000 & 2.01 & 7200 \\
INST19 & 159 &  & 0 & 0 & 0 & 102875 & INF & 7200 &  & 1 & 0 & 19333.3 & 19500 & 0.86 & 7200 \\
INST20 & 137 &  & 1 & 0 & 500 & 91375 & 18175 & 7200 &  & 1 & 0 & 16325 & 16500 & 1.07 & 7200 \\
INST21 & 110 &  & 1 & 0 & 8155.81 & 78125 & 857.906 & 7200 &  & 1 & 0 & 14025 & 14125 & 0.71 & 7200 \\
INST22 & 149 &  & 1 & 0 & 1000 & 97500 & 9650 & 7200 &  & 1 & 0 & 19675 & 20000 & 1.65 & 7200 \\
INST23 & 128 &  & 1 & 0 & 5009.87 & 27375 & 446.421 & 7200 &  & 1 & 1 & 14500 & 14500 & 0 & 605.17 \\
INST24 & 130 &  & 0 & 0 & 0 & 87583.3 & INF & 7200 &  & 1 & 0 & 15600 & 16000 & 2.56 & 7200 \\
INST25 & 122 &  & 1 & 0 & 8309.35 & 27791.7 & 234.463 & 7200 &  & 1 & 0 & 14316.7 & 14500 & 1.28 & 7200 \\
INST26 & 139 &  & 0 & 0 & 0 & 93291.7 & INF & 7200 &  & 1 & 1 & 16500 & 16500 & 0 & 3024.9 \\
INST27 & 142 &  & 0 & 0 & 0 & 93750 & INF & 7200 &  & 1 & 0 & 16258.3 & 16791.7 & 3.28 & 7200 \\
INST28 & 132 &  & 1 & 0 & 500 & 88916.7 & 17683.3 & 7200 &  & 1 & 0 & 14808.3 & 15000 & 1.29 & 7200 \\
INST29 & 143 &  & 1 & 0 & 8912.77 & 94875 & 964.484 & 7200 &  & 1 & 0 & 18366.7 & 19000 & 3.45 & 7200 \\
INST30 & 137 &  & 0 & 0 & 0 & 91000 & INF & 7200 &  & 1 & 0 & 16227.6 & 16500 & 1.68 & 7200 \\ \hline
\multicolumn{2}{c}{Ave.} &  & 0.53 & 0 & 2627.49 & 83727.8 & 7352.26 & 7200 &  & \textbf{1} & \textbf{0.17} & 17646.2 & 17918.1 & \textbf{1.49} & 6469.51 \\ \hline
\end{tabular}%
}
\label{tab_real_1}

\end{sidewaystable}

\begin{sidewaystable}
\doublespacing
\caption{Evaluation of the RGA\_CG, and BCIORPS algorithms}
\resizebox{\textwidth}{!}{%
\small
\begin{tabular}{ccccccccccccccccccc}
\hline
\multirow{2}{*}{\begin{tabular}[c]{@{}c@{}}Instance\\ Name\end{tabular}} & \multirow{2}{*}{\begin{tabular}[c]{@{}c@{}}$\lvert P \lvert$\end{tabular}} & \multirow{2}{*}{$\lvert D \lvert$} &  & \multicolumn{6}{c}{RGA\_CG} &  & \multicolumn{6}{c}{BCIORPS} &  & \multirow{2}{*}{\begin{tabular}[c]{@{}c@{}}\# Times CG \\ dominated BC\end{tabular}} \\ \cline{5-10} \cline{12-17}
 &  &  &  & \# Gen-Col & CG-Time & IP-Time & Time & UB & \# Solved &  & \# Nodes & \begin{tabular}[c]{@{}c@{}}\# Feasibility\\ Cuts\end{tabular} & \begin{tabular}[c]{@{}c@{}}\# Optimality\\ Cuts\end{tabular} & \# Solved & UB & \# Solved &  &  \\ \hline
INST31 & 282 & 12 &  & 7074 & 7200 & 6.54 & 7206.54 & 50708.3 & 1 &  & 0 & 0 & 6 & 7200 & \textbf{39625} & 1 &  & 0 \\
INST32 & 289 & 12 &  & 6013 & 7200 & 0.54 & 7200.54 & \textbf{40541.7} & 1 &  & 0 & 0 & 0 & 7200 & 42000 & 1 &  & 1 \\
INST33 & 306 & 12 &  & 5988 & 7200 & 0.63 & 7200.63 & \textbf{41500} & 1 &  & 0 & 0 & 0 & 7200 & 45000 & 1 &  & 1 \\
INST34 & 298 & 12 &  & 7106 & 7200 & 0.71 & 7200.71 & 47166.7 & 1 &  & 0 & 0 & 5 & 7200 & \textbf{44250} & 1 &  & 0 \\
INST35 & 319 & 12 &  & 7126 & 7200 & 146.19 & 7346.19 & \textbf{54416.7} & 1 &  & 0 & 0 & 1 & 7200 & INF & 0 &  & 1 \\
INST36 & 285 & 12 &  & 6542 & 7200 & 0.65 & 7200.65 & \textbf{39875} & 1 &  & 0 & 0 & 1 & 7200 & 59500 & 1 &  & 1 \\
INST37 & 262 & 12 &  & 5327 & 7200 & 313.84 & 7513.84 & 43291.7 & 1 &  & 20 & 0 & 0 & 7200 & \textbf{37625} & 1 &  & 0 \\
INST38 & 278 & 12 &  & 6300 & 7200 & 0.82 & 7200.82 & \textbf{39500} & 1 &  & 0 & 0 & 0 & 7200 & 40625 & 1 &  & 1 \\
INST39 & 319 & 12 &  & 6950 & 7200 & 3.02 & 7203.02 & 56500 & 1 &  & 50 & 0 & 0 & 7200 & \textbf{49541.7} & 1 &  & 0 \\
INST40 & 296 & 12 &  & 6035 & 7200 & 0.59 & 7200.59 & \textbf{41041.7} & 1 &  & 0 & 0 & 1 & 7200 & 66500 & 1 &  & 1 \\
INST41 & 259 & 12 &  & 6174 & 7200 & 0.56 & 7200.56 & \textbf{36500} & 1 &  & 0 & 0 & 0 & 7200 & 38625 & 1 &  & 1 \\
INST42 & 258 & 12 &  & 4987 & 7200 & 150.9 & 7350.9 & 43000 & 1 &  & 0 & 0 & 0 & 7200 & \textbf{37583.3} & 1 &  & 0 \\
INST43 & 261 & 12 &  & - & - & - & - & - & 0 &  & 0 & 0 & 0 & 7200 & \textbf{35583.3} & 1 &  & 0 \\
INST44 & 274 & 12 &  & 5217 & 7200 & 0.5 & 7200.5 & \textbf{35250} & 1 &  & 0 & 0 & 0 & 7200 & 40375 & 1 &  & 1 \\
INST45 & 280 & 12 &  & 6095 & 7200 & 3 & 7203 & \textbf{42000} & 1 &  & 0 & 0 & 0 & 7200 & 43291.7 & 1 &  & 1 \\
INST46 & 446 & 18 &  & 9384 & 7200 & 1523.68 & 8723.68 & \textbf{86208.3} & 1 &  & 0 & 0 & 1 & 7200 & 102500 & 1 &  & 1 \\
INST47 & 431 & 18 &  & 7848 & 7200 & 0.89 & 7200.89 & \textbf{62041.7} & 1 &  & 0 & 0 & 1 & 7200 & 94000 & 1 &  & 1 \\
INST48 & 462 & 18 &  & 10498 & 7200 & 1.39 & 7201.39 & 67541.7 & 1 &  & 0 & 0 & 21 & 7200 & \textbf{63083.3} & 1 &  & 0 \\
INST49 & 440 & 18 &  & 9088 & 7200 & 0.97 & 7200.97 & \textbf{74500} & 1 &  & 0 & 0 & 1 & 7200 & 79500 & 1 &  & 1 \\
INST50 & 399 & 18 &  & 11374 & 7200 & 40.18 & 7240.18 & \textbf{83208.3} & 1 &  & 0 & 0 & 1 & 7200 & INF & 0 &  & 1 \\
INST51 & 460 & 18 &  & 7672 & 7200 & 0.87 & 7200.87 & \textbf{53041.7} & 1 &  & 0 & 0 & 1 & 7200 & 63000 & 1 &  & 1 \\
INST52 & 406 & 18 &  & 9129 & 7200 & 1.03 & 7201.03 & \textbf{55125} & 1 &  & 17 & 0 & 0 & 7200 & 55458.3 & 1 &  & 1 \\
INST53 & 407 & 18 &  & 11136 & 7200 & 1.24 & 7201.24 & \textbf{66458.3} & 1 &  & 0 & 0 & 3 & 7200 & INF & 0 &  & 1 \\
INST54 & 403 & 18 &  & 10497 & 7200 & 1.52 & 7201.52 & \textbf{74000} & 1 &  & 0 & 0 & 1 & 7200 & 86500 & 1 &  & 1 \\
INST55 & 412 & 18 &  & 9829 & 7200 & 1.07 & 7201.07 & 57291.7 & 1 &  & 0 & 0 & 10 & 7200 & \textbf{57250} & 1 &  & 0 \\ \hline
Ave. &  &  &  & 7641.208 & 7200 & 91.72208 & 7291.722 &  & \textbf{24} &  & 3.48 & 0 & 2.16 & 7200 &  & 22 &  & \textbf{17} \\ \hline
\end{tabular}%
}
\label{tab_real_2}
\end{sidewaystable}

\newpage
\subsubsection{Disruption Analysis}\label{Emergency_Surgery_Disruption_Analysis}
Unplanned (non-elective) surgeries routinely perturb elective OR schedules and can trigger overtime, or surgery deferrals. Quantifying these effects and evaluating mitigation policies are therefore essential for decision support aligned with hospital KPIs (e.g., surgery postponement rate, and total cost).

Emergency-related disruptions manifest in several operationally distinct forms, including: (i) true emergency arrivals that require same-day insertion; (ii) priority escalation of already-scheduled elective cases due to clinical deterioration or external constraints; and (iii) demand bursts around weekends, holidays or special occasions, which are usually handled by reserved operating rooms and on-call surgeons.
This study focuses on the first two clinically prevalent scenarios:
\vspace{-0.10cm}
\begin{itemize}
\item True Emergency Arrivals: Exogenous, same-day cases are injected into the schedule, stressing insertion and capacity-management policies which are analyzed on synthetic instances.
\vspace{-0.10cm}
\item Priority Escalation via Due-Date Tightening: A subset of elective surgeries experiences a tightened due date (higher priority); this is analyzed on the Naples dataset.
\end{itemize}

Although the proposed model is developed for elective surgery planning and scheduling, it accommodates emergency disruptions through a rolling-horizon rescheduling methodology where near-term decisions are frozen, and the remaining horizon is re-optimized under the realized disruption. The two disruption scenarios were analyzed separately on synthetic and real-world (Naples) instances. This design demonstrates that our proposed model, coupled with a rolling-horizon rescheduling methodology, serves as an effective operational tool for day-by-day incorporation of emergency surgeries. Furthermore, the scenario–dataset pairing is not restrictive, and each scenario can be applied to either instance set without modification of the model. Moreover, combined disruption (e.g., simultaneous same-day emergency insertions and due-date tightening) is also could be naturally accommodated by the same framework.

\paragraph{Emergency Arrivals}\label{True_Emergency_Arrivals}
True emergency arrivals, which require same‐day insertion, are one of the disruption sources in IORPS. They compress capacity, trigger overtime and surgery postponement, and directly affect hospital KPIs. Evaluating how a baseline elective plan absorbs such shocks is therefore critical for operational robustness. On synthetic instances, we generated a baseline schedule and then introduced a same-day emergency shock on day 3 by removing three surgeries from the planning set, solving the deterministic problem to optimality (baseline), and subsequently re-inserting these three cases as emergency arrivals. Rescheduling was carried out in a rolling-horizon manner where decisions for days 1–2 were frozen at their baseline values, and the remaining horizon was re-optimized. Table \ref{true_emergency_arrivals_table} reports hospital KPIs including total cost as UB, postponements cost, OR-opening cost, overtime cost, and $\Delta\%$ which defines the percentage of difference between obtained UB of rescheduling and the baseline schedule. 

The average $\Delta\%$ over the synthetic instances is 4.13\% indicating that the proposed model, coupled with rolling-horizon rescheduling, accommodates true emergency insertions with modest deviation from the baseline. The absorption mechanism varies with problem size. For 40–60 surgery instances, the shock is primarily absorbed by opening an additional OR and limited overtime. For 80–100 surgery instances, the effect is absorbed through overtime and targeted postponements of elective cases. Finally, for 120 surgery instances, both levers are activated (extra ORs and selective postponements). These patterns suggest that the model flexibly trades off OR opening, overtime, and deferrals to preserve feasibility and control costs under acute, same-day demand surges.
In this experiment, we introduced three same-day emergency cases on Day 3 to demonstrate the model’s ability to accommodate disruptions imposed by new emergency arrivals. The same rolling-horizon methodology applies to any day of the planning horizon and to differing emergency loads (fewer or more cases) without modifying the formulation.

\begin{table}[H]
\caption{Hospital KPIs for true emergency arrivals}
\resizebox{\textwidth}{!}{%
\small
\begin{tabular}{ccccccccccccccccccccc}
\hline
\multirow{2}{*}{$\lvert P \lvert$} & \multirow{2}{*}{\# Instance} &  & \multicolumn{7}{c}{Initial Schedule} &  & \multicolumn{8}{c}{Rescheduling to consider emergency surgeries} \\ \cline{4-10} \cline{12-19} 
 &  &  & LB & UB & GAP\% & Time (s) & \begin{tabular}[c]{@{}c@{}}Postponement\\      Cost\end{tabular} & \begin{tabular}[c]{@{}c@{}}OR Opening\\      Cost\end{tabular} & \begin{tabular}[c]{@{}c@{}}Overtime\\      Cost\end{tabular} &  & LB & UB & GAP\% & Time (s) & \begin{tabular}[c]{@{}c@{}}Postponement\\      Cost\end{tabular} & \begin{tabular}[c]{@{}c@{}}OR Opening\\      Cost\end{tabular} & \begin{tabular}[c]{@{}c@{}}Overtime\\      Cost\end{tabular} & Delta \% \\ \hline
\multirow{5}{*}{40} & 1 &  & 7500 & 7500 & 0 & 366.68 & 500 & 7000 & 0 &  & 7500 & 7500 & 0 & 8.32 & 500 & 7000 & 0 & 0 \\
 & 2 &  & 9500 & 9500 & 0 & 27.82 & 500 & 9000 & 0 &  & 10000 & 10000 & 0 & 11.63 & 0 & 10000 & 0 & 5.26 \\
 & 3 &  & 8500 & 8500 & 0 & 2663.43 & 500 & 8000 & 0 &  & 9000 & 9000 & 0 & 1569.51 & 0 & 9000 & 0 & 5.88 \\
 & 4 &  & 7500 & 7500 & 0 & 1052.33 & 500 & 7000 & 0 &  & 8000 & 8000 & 0 & 24.18 & 0 & 8000 & 0 & 6.67 \\
 & 5 &  & 9000 & 9000 & 0 & 1610.53 & 1000 & 8000 & 0 &  & 9500 & 9500 & 0 & 19.15 & 500 & 9000 & 0 & 5.56 \\ \hline
\multicolumn{2}{c}{Average} &  & 8400 & 8400 & 0 & 1144.16 & 600 & 7800 & 0 &  & 8800 & 8800 & 0 & 326.558 & 200 & 8600 & 0 & 4.674 \\ \hline
\multirow{5}{*}{60} & 1 &  & 12138.3 & 12500 & 2.98 & 7200 & 500 & 12000 & 0 &  & 13500 & 13500 & 0 & 7.48 & 500 & 13000 & 0 & 8 \\
 & 2 &  & 16187.5 & 16500 & 1.93 & 7200 & 1500 & 15000 & 0 &  & 17232.8 & 17500 & 1.55 & 7200 & 1500 & 16000 & 0 & 6.06 \\
 & 3 &  & 14708.3 & 15000 & 1.98 & 7200 & 1000 & 14000 & 0 &  & 15500 & 15500 & 0 & 231.72 & 1500 & 14000 & 0 & 3.33 \\
 & 4 &  & 11875 & 12000 & 1.05 & 7200 & 0 & 12000 & 0 &  & 12875 & 12875 & 0 & 530.68 & 500 & 12000 & 375 & 7.29 \\
 & 5 &  & 12802.1 & 13000 & 1.55 & 7200 & 1000 & 12000 & 0 &  & 13125 & 13125 & 0 & 1728.73 & 0 & 13000 & 125 & 0.96 \\ \hline
\multicolumn{2}{c}{Average} &  & 13542.24 & 13800 & 1.898 & 7200 & 800 & 13000 & 0 &  & 14446.6 & 14500 & 0.31 & 1939.72 & 800 & 13600 & 100 & 5.128 \\ \hline
\multirow{5}{*}{80} & 1 &  & 18541.7 & 18541.7 & 0 & 345.18 & 2500 & 16000 & 41.6667 &  & 19041.7 & 19041.7 & 0 & 5743.49 & 3000 & 16000 & 41.7 & 2.7 \\
 & 2 &  & 18062.5 & 18250 & 1.04 & 7200 & 1000 & 17000 & 250 &  & 19333.3 & 19333.3 & 0 & 6838.67 & 2000 & 17000 & 333.3 & 5.94 \\
 & 3 &  & 17927.1 & 18000 & 0.41 & 7200 & 0 & 18000 & 0 &  & 18625 & 19000 & 2.01 & 7200 & 1000 & 18000 & 0 & 5.56 \\
 & 4 &  & 20812.5 & 21000 & 0.9 & 7200 & 2000 & 19000 & 0 &  & 21281.2 & 21500 & 1.03 & 7200 & 3500 & 18000 & 0 & 2.38 \\
 & 5 &  & 18666.7 & 19000 & 1.79 & 7200 & 3000 & 16000 & 0 &  & 19193.2 & 19500 & 1.6 & 7200 & 2500 & 17000 & 0 & 2.63 \\ \hline
\multicolumn{2}{c}{Average} &  & 18802.1 & 18958.3 & 0.828 & 5829.04 & 1700 & 17200 & 58.33334 &  & 19494.9 & 19675 & 0.928 & 6836.43 & 2400 & 17200 & 75 & 3.842 \\ \hline
\multirow{5}{*}{100} & 1 &  & 23500 & 23500 & 0 & 2219.44 & 3500 & 20000 & 0 &  & 24500 & 24500 & 0 & 4852.02 & 5500 & 19000 & 0 & 4.26 \\
 & 2 &  & 23312.5 & 23500 & 0.8 & 7200 & 2500 & 21000 & 0 &  & 24233.6 & 24500 & 1.1 & 7200 & 2500 & 22000 & 0 & 4.26 \\
 & 3 &  & 20937.5 & 21000 & 0.3 & 7200 & 2000 & 19000 & 0 &  & 22035.7 & 22125 & 0.41 & 7200 & 3000 & 19000 & 125 & 5.36 \\
 & 4 &  & 24083.3 & 24333.3 & 1.04 & 7200 & 2000 & 22000 & 333.333 &  & 24968.8 & 25208.3 & 0.96 & 7200 & 3000 & 22000 & 208.3 & 3.6 \\
 & 5 &  & 22000 & 22000 & 0 & 897.97 & 3000 & 19000 & 0 &  & 22500 & 22500 & 0 & 1108.05 & 3500 & 19000 & 0 & 2.27 \\ \hline
\multicolumn{2}{c}{Average} &  & 22766.66 & 22866.7 & 0.428 & 4943.48 & 2600 & 20200 & 66.6666 &  & 23647.6 & 23766.7 & 0.494 & 5512.01 & 3500 & 20200 & 66.66 & 3.95 \\ \hline
\multirow{5}{*}{120} & 1 &  & 29718.8 & 30541.7 & 2.77 & 7200 & 6500 & 24000 & 41.6667 &  & 30958.3 & 31041.7 & 0.27 & 7200 & 8000 & 23000 & 41.7 & 1.64 \\
 & 2 &  & 27989.6 & 28500 & 1.82 & 7200 & 5500 & 23000 & 0 &  & 29376.7 & 29500 & 0.42 & 7200 & 6500 & 23000 & 0 & 3.51 \\
 & 3 &  & 28760.4 & 29000 & 0.83 & 7200 & 8000 & 21000 & 0 &  & 29756.2 & 30000 & 0.82 & 7200 & 8000 & 22000 & 0 & 3.45 \\
 & 4 &  & 29770.8 & 30000 & 0.77 & 7200 & 6000 & 24000 & 0 &  & 32000 & 32000 & 0 & 2377.48 & 9000 & 23000 & 0 & 6.67 \\
 & 5 &  & 27583.3 & 28000 & 1.51 & 7200 & 5000 & 23000 & 0 &  & 28000 & 28000 & 0 & 1862.06 & 5000 & 23000 & 0 & 0 \\ \hline
\multicolumn{2}{c}{Average} &  & 28764.58 & 29208.3 & 1.54 & 7200 & 6200 & 23000 & 8.33334 &  & 30018.2 & 30108.3 & 0.302 & 5167.91 & 7300 & 22800 & 8.34 & 3.054 \\ \hline
\multicolumn{2}{c}{Total Average} &  & 18455.12 & 18646.7 & 0.94 & 5263.34 & 2380 & 16240 & 26.67 &  & 19281.5 & 19370 & 0.41 & 3956.53 & 2840 & 16480 & 50 & 4.13 \\ \hline
\end{tabular}%
}
\label{true_emergency_arrivals_table}
\end{table}

\paragraph{Priority Escalation via Due-Date Tightening}\label{Priority_Escalation_via_Due_Date_Tightening}
Beyond exogenous emergency arrivals, a common operational shock is endogenous reprioritization where the clinical status of an elective patients may deteriorate, effectively tightening their due dates and forcing resequencing. This mechanism stresses surgeon calendars and OR capacity, and is therefore a distinct robustness dimension relative to same-day insertions. Using the Naples (real‐world) instances, we first obtained a baseline schedule by solving the deterministic model with nominal parameters. We then constructed reprioritization scenarios by selecting a subset of $\{ 3, 5, 10\}$ surgeries and tightening their due dates to Day 4 (which means adding more mandatory surgeries to the problem). Same as the previous section, rescheduling followed a rolling‐horizon methodology where decisions for Days 1–2 were frozen at baseline values, and the remaining horizon was re‐optimized. Tables \ref{Priority_Escalation_via_Due_Date_Tightening_table1} and \ref{Priority_Escalation_via_Due_Date_Tightening_table2} show the detailed results for emergency surgery disruption via priority escalation. 

Excluding one instance (INST27) where the baseline stopped with a 2.51\% optimality gap while the reprioritized schedule solved to optimality, and yielded a negative $\Delta\%$ by construction, in all other scenarios, $\Delta\% < 0.5\%$, which is an indication for high robustness to moderate priority escalation. For small perturbations (3–5 tightened cases), the model predominantly exploits latent capacity within existing OR‐days, with minimal overtime and no additional ORs. At 10 tightened cases, these reserves are largely exhausted and the impact is absorbed through a combination of limited overtime and targeted postponements. As the number of reprioritizations increases, we expect the same hierarchy of levers observed for true emergencies which is opening additional ORs, expanding overtime, and selective deferrals. Feasibility can be lost if the cumulative mandatory load for a given surgeon exceeds their daily availability, since surgeries are preassigned to surgeons. Operational remedies include surgeon reassignment or temporary extension of surgeon availability. Incorporating these recourse options (e.g., surgeon‐switch variables with penalties) is a natural extension for future work.

\begin{table}[h]
\caption{Hospital KPIs for priority escalation via due-date tightening for 3 emergency arrivals}
\resizebox{\textwidth}{!}{%
\begin{tabular}{ccccccccccccccccccc}
\hline
\multirow{2}{*}{Instance Name} & \multirow{2}{*}{$\lvert P \lvert$} &  & \multicolumn{7}{c}{Initial Schedule} &  & \multicolumn{8}{c}{Rescheduling considering 3 emergency surgeries} \\ \cline{4-10} \cline{12-19} 
 &  &  & LB & UB & GAP\% & Time (s) & \begin{tabular}[c]{@{}c@{}}Postponement\\      Cost\end{tabular} & \begin{tabular}[c]{@{}c@{}}OR Opening\\      Cost\end{tabular} & \begin{tabular}[c]{@{}c@{}}Overtime\\      Cost\end{tabular} &  & LB & UB & GAP\% & Time (s) & \begin{tabular}[c]{@{}c@{}}Postponement\\      Cost\end{tabular} & \begin{tabular}[c]{@{}c@{}}OR Opening\\      Cost\end{tabular} & \begin{tabular}[c]{@{}c@{}}Overtime\\      Cost\end{tabular} & Delta   \% \\ \hline
INST1 & 143 &  & 17775 & 18000 & 1.27 & 7200 & 1000 & 17000 & 0 &  & 18000 & 18000 & 0 & 4036.54 & 1000 & 17000 & 0 & 0 \\
INST2 & 139 &  & 18000 & 18000 & 0 & 2391 & 1000 & 17000 & 0 &  & 18000 & 18000 & 0 & 167.65 & 1000 & 17000 & 0 & 0 \\
INST3 & 164 &  & 19708.3 & 20125 & 2.11 & 7200 & 2000 & 18000 & 125 &  & 20125 & 20125 & 0 & 5732.33 & 2000 & 18000 & 125 & 0 \\
INST4 & 125 &  & 16000 & 16000 & 0 & 1919.27 & 0 & 16000 & 0 &  & 16000 & 16000 & 0 & 34.67 & 0 & 16000 & 0 & 0 \\
INST5 & 160 &  & 18500 & 18625 & 0.68 & 7200 & 500 & 18000 & 125 &  & 18625 & 18625 & 0 & 32.5 & 500 & 18000 & 125 & 0 \\
INST6 & 146 &  & 18675 & 19000 & 1.74 & 7200 & 2000 & 17000 & 0 &  & 18766.7 & 19000 & 1.24 & 7200 & 2000 & 17000 & 0 & 0 \\
INST7 & 162 &  & 20941.7 & 21583.3 & 3.06 & 7200 & 3500 & 18000 & 83.3 &  & 21467.1 & 21583.3 & 0.54 & 7200 & 3500 & 18000 & 83.3 & 0 \\
INST8 & 136 &  & 17991.7 & 18500 & 2.83 & 7200 & 1500 & 17000 & 0 &  & 19000 & 19000 & 0 & 5786 & 1000 & 18000 & 0 & 2.7 \\
INST9 & 164 &  & 23033.3 & 23500 & 2.03 & 7200 & 4500 & 19000 & 0 &  & 23083.3 & 23500 & 1.81 & 7200 & 4500 & 19000 & 0 & 0 \\
INST10 & 155 &  & 19316.7 & 19500 & 0.95 & 7200 & 2500 & 17000 & 0 &  & 19425 & 19500 & 0.39 & 7200 & 2500 & 17000 & 0 & 0 \\
INST11 & 136 &  & 14758.3 & 15000 & 1.64 & 7200 & 1000 & 14000 & 0 &  & 15000 & 15000 & 0 & 1549.15 & 1000 & 14000 & 0 & 0 \\
INST12 & 149 &  & 18091.7 & 18500 & 2.26 & 7200 & 500 & 18000 & 0 &  & 18308.3 & 18500 & 1.05 & 7200 & 500 & 18000 & 0 & 0 \\
INST13 & 129 &  & 14408.3 & 14500 & 0.64 & 7200 & 500 & 14000 & 0 &  & 14500 & 14500 & 0 & 1420.27 & 500 & 14000 & 0 & 0 \\
INST14 & 133 &  & 17716.7 & 18000 & 1.6 & 7200 & 2000 & 16000 & 0 &  & 18000 & 18000 & 0 & 2893.92 & 2000 & 16000 & 0 & 0 \\
INST15 & 137 &  & 16308.3 & 16500 & 1.18 & 7200 & 500 & 16000 & 0 &  & 16500 & 16500 & 0 & 2249.2 & 500 & 16000 & 0 & 0 \\
INST16 & 141 &  & 18541.7 & 18541.7 & 0 & 6143.16 & 2500 & 16000 & 41.7 &  & 19000 & 19000 & 0 & 4370.23 & 3000 & 16000 & 0 & 2.47 \\
INST17 & 157 &  & 22116.7 & 22500 & 1.73 & 7200 & 3500 & 19000 & 0 &  & 22500 & 22500 & 0 & 2351.39 & 3500 & 19000 & 0 & 0 \\
INST18 & 162 &  & 21566.7 & 22000 & 2.01 & 7200 & 3000 & 19000 & 0 &  & 21769.7 & 22000 & 1.06 & 7200 & 3000 & 19000 & 0 & 0 \\
INST19 & 159 &  & 19333.3 & 19500 & 0.86 & 7200 & 1500 & 18000 & 0 &  & 19408.3 & 19500 & 0.47 & 7200 & 1500 & 18000 & 0 & 0 \\
INST20 & 137 &  & 16325 & 16500 & 1.07 & 7200 & 1500 & 15000 & 0 &  & 16500 & 16500 & 0 & 666.53 & 1500 & 15000 & 0 & 0 \\
INST21 & 110 &  & 14025 & 14125 & 0.71 & 7200 & 1000 & 13000 & 125 &  & 14125 & 14166.7 & 0.3 & 7200 & 1000 & 13000 & 166.7 & 0.3 \\
INST22 & 149 &  & 19675 & 20000 & 1.65 & 7200 & 3000 & 17000 & 0 &  & 19843.9 & 20000 & 0.79 & 7200 & 3000 & 17000 & 0 & 0 \\
INST23 & 128 &  & 14500 & 14500 & 0 & 605.17 & 500 & 14000 & 0 &  & 14500 & 14500 & 0 & 33.33 & 500 & 14000 & 0 & 0 \\
INST24 & 130 &  & 15600 & 16000 & 2.56 & 7200 & 1000 & 15000 & 0 &  & 16000 & 16000 & 0 & 1830.89 & 1000 & 15000 & 0 & 0 \\
INST25 & 122 &  & 14316.7 & 14500 & 1.28 & 7200 & 500 & 14000 & 0 &  & 14500 & 14500 & 0 & 1566.02 & 500 & 14000 & 0 & 0 \\
INST26 & 139 &  & 16500 & 16500 & 0 & 3024.9 & 500 & 16000 & 0 &  & 16500 & 16500 & 0 & 41.43 & 500 & 16000 & 0 & 0 \\
INST27 & 142 &  & 16258.3 & 16666.7 & 2.51 & 7200 & 500 & 16000 & 166.7 &  & 16666.7 & 16666.7 & 0 & 2804.48 & 500 & 16000 & 166.7 & 0 \\
INST28 & 132 &  & 14808.3 & 15000 & 1.29 & 7200 & 1000 & 14000 & 0 &  & 15000 & 15000 & 0 & 2722.84 & 1000 & 14000 & 0 & 0 \\
INST29 & 143 &  & 18366.7 & 19000 & 3.45 & 7200 & 2000 & 17000 & 0 &  & 19000 & 19000 & 0 & 2165.69 & 2000 & 17000 & 0 & 0 \\
INST30 & 137 &  & 16227.6 & 16500 & 1.68 & 7200 & 500 & 16000 & 0 &  & 16275 & 16500 & 1.38 & 7200 & 500 & 16000 & 0 & 0 \\ \hline
\multicolumn{2}{c}{Average} &  & 17646.2 & 17905.6 & 1.43 & 6469.45 & 1516.67 & 16366.67 & 22.22 &  & 17879.6 & 17938.9 & 0.3 & 3815.17 & 1516.67 & 16400 & 22.22 & 0.18 \\ \hline
\end{tabular}%
}
\label{Priority_Escalation_via_Due_Date_Tightening_table1}
\end{table}

\begin{table}[H]
\caption{Hospital KPIs for priority escalation via due-date tightening for 5 and 10 emergency arrivals}
\resizebox{\textwidth}{!}{%
\begin{tabular}{cccccccccccccccccccc}
\hline
\multirow{2}{*}{Instance Name} & \multirow{2}{*}{$\lvert P \lvert$} &  & \multicolumn{8}{c}{Rescheduling considering 5 emergency surgeries} &  & \multicolumn{8}{c}{Rescheduling considering 10 emergency surgeries} \\ \cline{4-11} \cline{13-20} 
 &  &  & LB & UB & GAP\% & Time (s) & \begin{tabular}[c]{@{}c@{}}Postponement\\      Cost\end{tabular} & \begin{tabular}[c]{@{}c@{}}OR Opening\\      Cost\end{tabular} & \begin{tabular}[c]{@{}c@{}}Overtime\\      Cost\end{tabular} & Delta   \% &  & LB & UB & GAP\% & Time (s) & \begin{tabular}[c]{@{}c@{}}Postponement\\      Cost\end{tabular} & \begin{tabular}[c]{@{}c@{}}OR Opening\\      Cost\end{tabular} & \begin{tabular}[c]{@{}c@{}}Overtime\\      Cost\end{tabular} & Delta   \% \\ \hline
INST1 & 143 &  & 18000 & 18000 & 0 & 5640.12 & 1000 & 17000 & 0 & 0 &  & 17883.3 & 18000 & 0.65 & 7200 & 1000 & 17000 & 0 & 0 \\
INST2 & 139 &  & 18000 & 18000 & 0 & 126.81 & 1000 & 17000 & 0 & 0 &  & 18000 & 18000 & 0 & 317.06 & 1000 & 17000 & 0 & 0 \\
INST3 & 164 &  & 20022.1 & 20125 & 0.51 & 7200 & 2000 & 18000 & 125 & 0 &  & 20222.1 & 20625 & 1.99 & 7200 & 2500 & 18000 & 125 & 2.48 \\
INST4 & 125 &  & 16000 & 16000 & 0 & 94.72 & 0 & 16000 & 0 & 0 &  & 16000 & 16000 & 0 & 33.2 & 0 & 16000 & 0 & 0 \\
INST5 & 160 &  & 18625 & 18625 & 0 & 2973.65 & 500 & 18000 & 125 & 0 &  & 18500 & 18625 & 0.68 & 7200 & 500 & 18000 & 125 & 0 \\
INST6 & 146 &  & 19000 & 19000 & 0 & 1495.4 & 2000 & 17000 & 0 & 0 &  & 18766.7 & 19000 & 1.24 & 7200 & 2000 & 17000 & 0 & 0 \\
INST7 & 162 &  & 21441.7 & 21583.3 & 0.66 & 7200 & 3500 & 18000 & 83.3 & 0 &  & 21479.7 & 21583.3 & 0.48 & 7200 & 3500 & 18000 & 83.3 & 0 \\
INST8 & 136 &  & 19000 & 19000 & 0 & 6832 & 1000 & 18000 & 0 & 2.7 &  & 19000 & 19208.3 & 1.1 & 7200 & 1000 & 18000 & 208.3 & 3.83 \\
INST9 & 164 &  & 23083.3 & 23500 & 1.81 & 7200 & 4500 & 19000 & 0 & 0 &  & 23083.3 & 23500 & 1.81 & 7200 & 4500 & 19000 & 0 & 0 \\
INST10 & 155 &  & 19416.7 & 19500 & 0.43 & 7200 & 2500 & 17000 & 0 & 0 &  & 19416.7 & 19500 & 0.43 & 7200 & 2500 & 17000 & 0 & 0 \\
INST11 & 136 &  & 15000 & 15000 & 0 & 2646.1 & 1000 & 14000 & 0 & 0 &  & 15000 & 15000 & 0 & 6272.6 & 1000 & 14000 & 0 & 0 \\
INST12 & 149 &  & 18315.5 & 18500 & 1.01 & 7200 & 500 & 18000 & 0 & 0 &  & 18500 & 18500 & 0 & 5047.48 & 500 & 18000 & 0 & 0 \\
INST13 & 129 &  & 14500 & 14500 & 0 & 3044.41 & 500 & 14000 & 0 & 0 &  & 14500 & 14500 & 0 & 6467.49 & 500 & 14000 & 0 & 0 \\
INST14 & 133 &  & 18000 & 18000 & 0 & 1815.02 & 2000 & 16000 & 0 & 0 &  & 18000 & 18000 & 0 & 623.91 & 2000 & 16000 & 0 & 0 \\
INST15 & 137 &  & 16358.3 & 16500 & 0.87 & 7200 & 500 & 16000 & 0 & 0 &  & 16500 & 16500 & 0 & 2181.77 & 500 & 16000 & 0 & 0 \\
INST16 & 141 &  & 19000 & 19000 & 0 & 3809.14 & 3000 & 16000 & 0 & 2.47 &  & 19500 & 19500 & 0 & 3809.14 & 3500 & 16000 & 0 & 5.17 \\
INST17 & 157 &  & 22391.7 & 22500 & 0.48 & 7200 & 3500 & 19000 & 0 & 0 &  & 22500 & 22500 & 0 & 1807.87 & 3500 & 19000 & 0 & 0 \\
INST18 & 162 &  & 21766.7 & 22000 & 1.07 & 7200 & 3000 & 19000 & 0 & 0 &  & 21793.4 & 22000 & 0.95 & 7200 & 3000 & 19000 & 0 & 0 \\
INST19 & 159 &  & 19500 & 19500 & 0 & 7200 & 1500 & 18000 & 0 & 0 &  & 19500 & 19500 & 0 & 2131.16 & 1500 & 18000 & 0 & 0 \\
INST20 & 137 &  & 16500 & 16500 & 0 & 5756.72 & 1500 & 15000 & 0 & 0 &  & 16500 & 16500 & 0 & 1293.52 & 1500 & 15000 & 0 & 0 \\
INST21 & 110 &  & 14125 & 14166.7 & 0.3 & 7200 & 1000 & 13000 & 166.7 & 0.3 &  & 14125 & 14166.7 & 0.3 & 7200 & 1000 & 13000 & 166.7 & 0.3 \\
INST22 & 149 &  & 19907.4 & 20000 & 0.47 & 7201 & 3000 & 17000 & 0 & 0 &  & 20000 & 20041.667 & 0.21 & 7200 & 3000 & 17000 & 41.7 & 0.21 \\
INST23 & 128 &  & 14500 & 14500 & 0 & 88.47 & 500 & 14000 & 0 & 0 &  & 14500 & 14500 & 0 & 34.84 & 500 & 14000 & 0 & 0 \\
INST24 & 130 &  & 16000 & 16000 & 0 & 778.56 & 1000 & 15000 & 0 & 0 &  & 16000 & 16000 & 0 & 3150.68 & 1000 & 15000 & 0 & 0 \\
INST25 & 122 &  & 14500 & 14500 & 0 & 1214.44 & 500 & 14000 & 0 & 0 &  & 14500 & 14500 & 0 & 3121.61 & 500 & 14000 & 0 & 0 \\
INST26 & 139 &  & 16500 & 16500 & 0 & 171.33 & 500 & 16000 & 0 & 0 &  & 16500 & 16500 & 0 & 49.6 & 500 & 16000 & 0 & 0 \\
INST27 & 142 &  & 16500 & 16500 & 0 & 2266.76 & 500 & 16000 & 0 & -1 &  & 16500 & 16500 & 0 & 3955.97 & 500 & 16000 & 0 & -1 \\
INST28 & 132 &  & 15000 & 15000 & 0 & 5699.82 & 1000 & 14000 & 0 & 0 &  & 15000 & 15000 & 0 & 1742.37 & 1000 & 14000 & 0 & 0 \\
INST29 & 143 &  & 19000 & 19000 & 0 & 4032.2 & 2000 & 17000 & 0 & 0 &  & 19000 & 19000 & 0 & 3311.05 & 2000 & 17000 & 0 & 0 \\
INST30 & 137 &  & 16500 & 16500 & 0 & 3715.55 & 500 & 16000 & 0 & 0 &  & 16500 & 16500 & 0 & 1135.91 & 500 & 16000 & 0 & 0 \\ \hline
\multicolumn{2}{c}{Average} &  & 17881.8 & 17933.3 & 0.25 & 4380.07 & 1516.67 & 16400 & 16.67 & 0.15 &  & 17909.01 & 17975 & 0.33 & 4189.57 & 1550 & 16400 & 25 & 0.37 \\ \hline
\end{tabular}%
}
\label{Priority_Escalation_via_Due_Date_Tightening_table2}
\end{table}

\subsubsection{Robustness to Surgery Duration Variability and Buffer-Time Analysis}\label{Robustness_Duration_Variability_Buffer_Analysis}
Deterministic schedules that ignore inherent variability in surgery durations systematically misestimate capacity, which in practice manifests as overtime, postponements, and reduced on-time performance. We therefore quantify the operational impact of duration uncertainty and evaluate time-buffer policies that mitigate these effects while aligning with hospital KPIs. Let $\tilde{d}_i$ denote the realized surgery durations for surgery $i$ which is modeled as:

\begin{flalign}\label{Eq_buffer}
 \tilde{d}_i = d_i(1+\epsilon_i),   \qquad \qquad \qquad \epsilon_i \sim U[-0.2, 0.2] &&
\end{flalign}

For fair comparisons across policies, the same realization of $\epsilon_i$ is used for all buffer times $\in B$ on a given instance. We study buffer policies $B \in \{0, 30, 60, 90, 120 \}$ minutes per OR on each day. These policies are implemented by reducing the regular scheduling capacity of each OR–day from $\lvert T \lvert$ to $\lvert T \lvert - B$. All other inputs (surgeons' availability, cost parameters and due dates, etc.) are held fixed. For each $B$, the deterministic MILP is solved once with original values for $d_i$ to obtain $x_{idt}$, $y_d$, and $z_i$. Holding these decisions fixed, we recompute realized objective components under $\tilde{d}_i$ while treating buffer time as consumed before any overtime. This design isolates the policy effect of buffering without re-optimization and reflects how deterministic plans perform when executed under stochastic durations. Finally, we report results on hospital-aligned KPIs including total cost, overtime cost, postponements cost, and OR-opening cost. We consider buffer policies $B \in \{0, 30, 60, 90, 120 \}$ minutes per OR–day to span a practically relevant range with interpretable increments. The baseline $B=0$ quantifies the cost of ignoring variability; 30–90 minutes reflect common scheduling granularities and turnover cushions in hospital practice; and $B=120$ represents a conservative upper bound. Beyond this point, additional buffering becomes operationally unrealistic and more likely reflects systematic misestimation of case durations rather than genuine stochastic variability. Values exceeding 120 minutes would excessively depress throughput and are better addressed by improving duration forecasts (e.g., via historical-case modeling, statistical estimation, or machine-learning predictors and estimators) rather than by enlarging buffers. This discrete grid thus balances external validity (aligning with realistic OR-day operations) and experimental sensitivity (revealing the trade-off between overtime reduction and postponements) without resorting to implausible buffer magnitudes.

\paragraph{Synthetic Instances}\label{Synthetic_Instances_Buffer_Analysis}
We assess buffer policies $B \in \{0, 30, 60, 90, 120 \}$ minutes per OR–day on synthetic instances. Tables \ref{Synthetic_Instances_Buffer_Analysis_table1} and \ref{Synthetic_Instances_Buffer_Analysis_table2} report per-instance outcomes and per-size means to facilitate interpretation. Across all sizes, the mean total cost decreases monotonically as $B$ increases, with the dominant driver being a marked reduction in overtime as buffers grow. Specifically, postponement cost rises initially and peaks at $B=90$ minutes, then declines slightly (from 2740 to 2680 on average), while OR-opening cost changes modestly (e.g., 16,560 to 16,640). The largest marginal gain up to this point, occurs between $B=60$ and $B=90$, where overtime drops by roughly 40\% (relative to $B=90$), which shows a favorable trade-off in which limited additional postponements and slightly fewer opened ORs are exchanged for substantial overtime reduction.

At $B=120$, overtime is reduced further, and the model achieves the lowest total mean cost with a small increase in OR openings and a reduction in postponements relative to previous buffer value. With no buffer the system absorbs duration variability via residual OR capacity and overtime, and adding 30-60 minutes yields only modest total-cost improvements (sub-1\%) as slight overtime relief is offset by new postponements. Beyond 60 minutes, buffers increasingly substitute for overtime, and by 120 minutes nearly all variability is accommodated within regular hours, with at most a few time slots (i.e., 25 minutes) of overtime in the worst cases. Under the tested $\pm20\%$ duration variability, $B=120$ minutes emerges as the preferred policy on synthetic instances, balancing overtime reduction against minimal increase in postponement costs and the same OR opening cost.

\begin{table}[H]
\caption{Hospital KPIs of buffer analysis for synthetic instances considering 30 and 60 minutes of buffers}
\resizebox{\textwidth}{!}{%
\begin{tabular}{ccccccccccccccccc}
\hline
\multirow{2}{*}{$\lvert P \lvert$} & \multirow{2}{*}{\# Instance} &  & \multicolumn{4}{c}{No Buffer - Baseline} &  & \multicolumn{4}{c}{Buffer = 30 minutes} &  & \multicolumn{4}{c}{Buffer = 60 minutes} \\ \cline{4-7} \cline{9-12} \cline{14-17} 
 &  &  & Total Cost & \begin{tabular}[c]{@{}c@{}}Postponement\\      Cost\end{tabular} & \begin{tabular}[c]{@{}c@{}}OR   Opening\\      Cost\end{tabular} & \begin{tabular}[c]{@{}c@{}}Overtime\\      Cost\end{tabular} &  & Total Cost & \begin{tabular}[c]{@{}c@{}}Postponement\\      Cost\end{tabular} & \begin{tabular}[c]{@{}c@{}}OR Opening\\      Cost\end{tabular} & \begin{tabular}[c]{@{}c@{}}Overtime\\      Cost\end{tabular} &  & Total Cost & \begin{tabular}[c]{@{}c@{}}Postponement\\      Cost\end{tabular} & \begin{tabular}[c]{@{}c@{}}OR Opening\\      Cost\end{tabular} & \begin{tabular}[c]{@{}c@{}}Overtime\\      Cost\end{tabular} \\ \hline
\multirow{5}{*}{40} & 1 &  & 8791.7 & 500 & 7000 & 1291.7 &  & 8791.7 & 500 & 7000 & 1291.7 &  & 8875 & 500 & 7000 & 1375 \\
 & 2 &  & 11583.3 & 0 & 10000 & 1583.3 &  & 11250 & 0 & 10000 & 1250 &  & 11125 & 0 & 10000 & 1125 \\
 & 3 &  & 9916.7 & 0 & 9000 & 916.7 &  & 9541.7 & 0 & 9000 & 541.7 &  & 10083.3 & 0 & 9000 & 1083.3 \\
 & 4 &  & 8625 & 0 & 8000 & 625 &  & 9041.7 & 0 & 8000 & 1041.7 &  & 9083.3 & 1000 & 7000 & 1083.3 \\
 & 5 &  & 10708.3 & 500 & 9000 & 1208.3 &  & 10333.3 & 500 & 9000 & 833.3 &  & 10291.7 & 500 & 9000 & 791.7 \\ \hline
\multicolumn{2}{c}{Average} &  & 9925 & 200 & 8600 & 1125 &  & 9791.7 & 200 & 8600 & 991.7 &  & 9891.7 & 400 & 8400 & 1091.7 \\ \hline
\multirow{5}{*}{60} & 1 &  & 14833.3 & 0 & 13000 & 1833.3 &  & 14875 & 0 & 13000 & 1875 &  & 14958.3 & 0 & 13000 & 1958.3 \\
 & 2 &  & 20166.7 & 1000 & 16000 & 3166.7 &  & 20291.7 & 1000 & 16000 & 3291.7 &  & 19041.7 & 2000 & 15000 & 2041.7 \\
 & 3 &  & 17666.7 & 1500 & 14000 & 2166.7 &  & 17333.3 & 1500 & 14000 & 1833.3 &  & 16916.7 & 1500 & 14000 & 1416.7 \\
 & 4 &  & 15500 & 500 & 12000 & 3000 &  & 14833.3 & 500 & 12000 & 2333.3 &  & 14791.7 & 500 & 12000 & 2291.7 \\
 & 5 &  & 14708.3 & 0 & 13000 & 1708.3 &  & 14333.3 & 0 & 13000 & 1333.3 &  & 14625 & 0 & 13000 & 1625 \\ \hline
\multicolumn{2}{c}{Average} &  & 16575 & 600 & 13600 & 2375 &  & 16333.3 & 600 & 13600 & 2133.3 &  & 16066.7 & 800 & 13400 & 1866.7 \\ \hline
\multirow{5}{*}{80} & 1 &  & 21000 & 3000 & 16000 & 2000 &  & 21333.3 & 3000 & 16000 & 2333.3 &  & 20916.7 & 3000 & 16000 & 1916.7 \\
 & 2 &  & 22458.3 & 2000 & 17000 & 3458.3 &  & 21875 & 2000 & 17000 & 2875 &  & 21791.7 & 2000 & 17000 & 2791.7 \\
 & 3 &  & 21041.7 & 1000 & 18000 & 2041.7 &  & 21875 & 1000 & 18000 & 2875 &  & 20666.7 & 1000 & 18000 & 1666.7 \\
 & 4 &  & 24125 & 2500 & 19000 & 2625 &  & 24166.7 & 2500 & 19000 & 2666.7 &  & 24333.3 & 2500 & 19000 & 2833.3 \\
 & 5 &  & 21375 & 2500 & 17000 & 1875 &  & 21541.7 & 2500 & 17000 & 2041.7 &  & 21291.7 & 2500 & 17000 & 1791.7 \\ \hline
\multicolumn{2}{c}{Average} &  & 22000 & 2200 & 17400 & 2400 &  & 22158.3 & 2200 & 17400 & 2558.3 &  & 21800 & 2200 & 17400 & 2200 \\ \hline
\multirow{5}{*}{100} & 1 &  & 27250 & 4500 & 20000 & 2750 &  & 27291.7 & 4500 & 20000 & 2791.7 &  & 27458.3 & 4500 & 20000 & 2958.3 \\
 & 2 &  & 27291.7 & 2500 & 22000 & 2791.7 &  & 27416.7 & 2500 & 22000 & 2916.7 &  & 27458.3 & 2500 & 22000 & 2958.3 \\
 & 3 &  & 25666.7 & 3000 & 19000 & 3666.7 &  & 25666.7 & 3000 & 19000 & 3666.7 &  & 25625 & 3000 & 19000 & 3625 \\
 & 4 &  & 27625 & 3000 & 22000 & 2625 &  & 26958.3 & 3000 & 22000 & 1958.3 &  & 27583.3 & 3000 & 22000 & 2583.3 \\
 & 5 &  & 25833.3 & 3500 & 19000 & 3333.3 &  & 25833.3 & 3500 & 19000 & 3333.3 &  & 24875 & 3500 & 19000 & 2375 \\ \hline
\multicolumn{2}{c}{Average} &  & 26733.3 & 3300 & 20400 & 3033.3 &  & 26633.3 & 3300 & 20400 & 2933.3 &  & 26600 & 3300 & 20400 & 2900 \\ \hline
\multirow{5}{*}{120} & 1 &  & 34125 & 7000 & 24000 & 3125 &  & 34333.3 & 7000 & 24000 & 3333.3 &  & 34166.7 & 7000 & 24000 & 3166.7 \\
 & 2 &  & 32541.7 & 6500 & 23000 & 3041.7 &  & 33000 & 6500 & 23000 & 3500 &  & 32666.7 & 5500 & 24000 & 3166.7 \\
 & 3 &  & 32208.3 & 8000 & 22000 & 2208.3 &  & 32208.3 & 8000 & 22000 & 2208.3 &  & 32833.3 & 8000 & 22000 & 2833.3 \\
 & 4 &  & 35291.7 & 7500 & 24000 & 3791.7 &  & 35250 & 8000 & 24000 & 3250 &  & 34541.7 & 8000 & 24000 & 2541.7 \\
 & 5 &  & 30708.3 & 5000 & 23000 & 2708.3 &  & 30333.3 & 5000 & 23000 & 2333.3 &  & 30750 & 5000 & 23000 & 2750 \\ \hline
\multicolumn{2}{c}{Average} &  & 32975 & 6800 & 23200 & 2975 &  & 33025 & 6900 & 23200 & 2925 &  & 32991.7 & 6700 & 23400 & 2891.7 \\ \hline
\multicolumn{2}{c}{Total Average} &  & 21641.7 & 2620 & 16640 & 2381.7 &  & 21588.3 & 2640 & 16640 & 2308.3 &  & 21470 & 2680 & 16600 & 2190 \\ \hline
\end{tabular}%
}
\label{Synthetic_Instances_Buffer_Analysis_table1}
\end{table}

\begin{table}[H]
\caption{Hospital KPIs of buffer analysis for synthetic instances considering 90 and 120 minutes of buffers}
\resizebox{\textwidth}{!}{%
\tiny
\begin{tabular}{cccccccccccc}
\hline
\multirow{2}{*}{$\lvert P \lvert$} & \multirow{2}{*}{\# Instance} &  & \multicolumn{4}{c}{Buffer = 90 minutes} &  & \multicolumn{4}{c}{Buffer = 120 minutes} \\ \cline{4-7} \cline{9-12} 
 &  &  & Total Cost & \begin{tabular}[c]{@{}c@{}}Postponement\\      Cost\end{tabular} & \begin{tabular}[c]{@{}c@{}}OR Opening\\      Cost\end{tabular} & \begin{tabular}[c]{@{}c@{}}Overtime\\      Cost\end{tabular} &  & Total Cost & \begin{tabular}[c]{@{}c@{}}Postponement\\      Cost\end{tabular} & \begin{tabular}[c]{@{}c@{}}OR Opening\\      Cost\end{tabular} & \begin{tabular}[c]{@{}c@{}}Overtime\\      Cost\end{tabular} \\ \hline
\multirow{5}{*}{40} & 1 &  & 7916.7 & 500 & 7000 & 416.7 &  & 7541.7 & 500 & 7000 & 41.7 \\
 & 2 &  & 10750 & 0 & 10000 & 750 &  & 10000 & 0 & 10000 & 0 \\
 & 3 &  & 9666.7 & 0 & 9000 & 666.7 &  & 9041.7 & 0 & 9000 & 41.7 \\
 & 4 &  & 8416.7 & 0 & 8000 & 416.7 &  & 8041.7 & 0 & 8000 & 41.7 \\
 & 5 &  & 10000 & 500 & 9000 & 500 &  & 9583.3 & 500 & 9000 & 83.3 \\ \hline
\multicolumn{2}{c}{Average} &  & 9350 & 200 & 8600 & 550 &  & 8841.7 & 200 & 8600 & 41.7 \\ \hline
\multirow{5}{*}{60} & 1 &  & 14541.7 & 0 & 13000 & 1541.7 &  & 13000 & 0 & 13000 & 0 \\
 & 2 &  & 18208.3 & 2000 & 15000 & 1208.3 &  & 17000 & 1000 & 16000 & 0 \\
 & 3 &  & 16875 & 1500 & 14000 & 1375 &  & 15541.7 & 1500 & 14000 & 41.7 \\
 & 4 &  & 14041.7 & 500 & 12000 & 1541.7 &  & 12541.7 & 500 & 12000 & 41.7 \\
 & 5 &  & 13916.7 & 1000 & 12000 & 916.7 &  & 13125 & 0 & 13000 & 125 \\ \hline
\multicolumn{2}{c}{Average} &  & 15516.7 & 1000 & 13200 & 1316.7 &  & 14241.7 & 600 & 13600 & 41.7 \\ \hline
\multirow{5}{*}{80} & 1 &  & 20541.7 & 3000 & 16000 & 1541.7 &  & 19041.7 & 3000 & 16000 & 41.7 \\
 & 2 &  & 20083.3 & 2000 & 17000 & 1083.3 &  & 19166.7 & 2000 & 17000 & 166.7 \\
 & 3 &  & 19916.7 & 1000 & 18000 & 916.7 &  & 19041.7 & 1000 & 18000 & 41.7 \\
 & 4 &  & 23250 & 2500 & 19000 & 1750 &  & 21500 & 2500 & 19000 & 0 \\
 & 5 &  & 20375 & 2500 & 17000 & 875 &  & 19708.3 & 2500 & 17000 & 208.3 \\ \hline
\multicolumn{2}{c}{Average} &  & 20833.3 & 2200 & 17400 & 1233.3 &  & 19691.7 & 2200 & 17400 & 91.7 \\ \hline
\multirow{5}{*}{100} & 1 &  & 26666.7 & 4500 & 20000 & 2166.7 &  & 24708.3 & 4500 & 20000 & 208.3 \\
 & 2 &  & 25416.7 & 2500 & 22000 & 916.7 &  & 24625 & 2500 & 22000 & 125 \\
 & 3 &  & 23750 & 3000 & 19000 & 1750 &  & 22041.7 & 3000 & 19000 & 41.7 \\
 & 4 &  & 26666.7 & 3000 & 22000 & 1666.7 &  & 25125 & 3000 & 22000 & 125 \\
 & 5 &  & 24083.3 & 3500 & 19000 & 1583.3 &  & 22625 & 3500 & 19000 & 125 \\ \hline
\multicolumn{2}{c}{Average} &  & 25316.7 & 3300 & 20400 & 1616.7 &  & 23825 & 3300 & 20400 & 125 \\ \hline
\multirow{5}{*}{120} & 1 &  & 33708.3 & 7000 & 24000 & 2708.3 &  & 31666.7 & 7500 & 24000 & 166.7 \\
 & 2 &  & 32333.3 & 7000 & 23000 & 2333.3 &  & 30166.7 & 6000 & 24000 & 166.7 \\
 & 3 &  & 31791.7 & 8000 & 22000 & 1791.7 &  & 30041.7 & 9000 & 21000 & 41.7 \\
 & 4 &  & 33166.7 & 8000 & 24000 & 1166.7 &  & 32041.7 & 8000 & 24000 & 41.7 \\
 & 5 &  & 29541.7 & 5000 & 23000 & 1541.7 &  & 28083.3 & 5000 & 23000 & 83.3 \\ \hline
\multicolumn{2}{c}{Average} &  & 32108.3 & 7000 & 23200 & 1908.3 &  & 30400 & 7100 & 23200 & 100 \\ \hline
\multicolumn{2}{c}{Total Average} &  & 20625 & 2740 & 16560 & 1325 &  & 19400 & 2680 & 16640 & 80 \\ \hline
\end{tabular}%
}
\label{Synthetic_Instances_Buffer_Analysis_table2}
\end{table}

\begin{figure}[H]
\centering
\resizebox{0.85\textwidth}{!}{
\begin{minipage}{\textwidth}
  \begin{subfigure}{0.48\textwidth}
    \includegraphics[width=\linewidth]{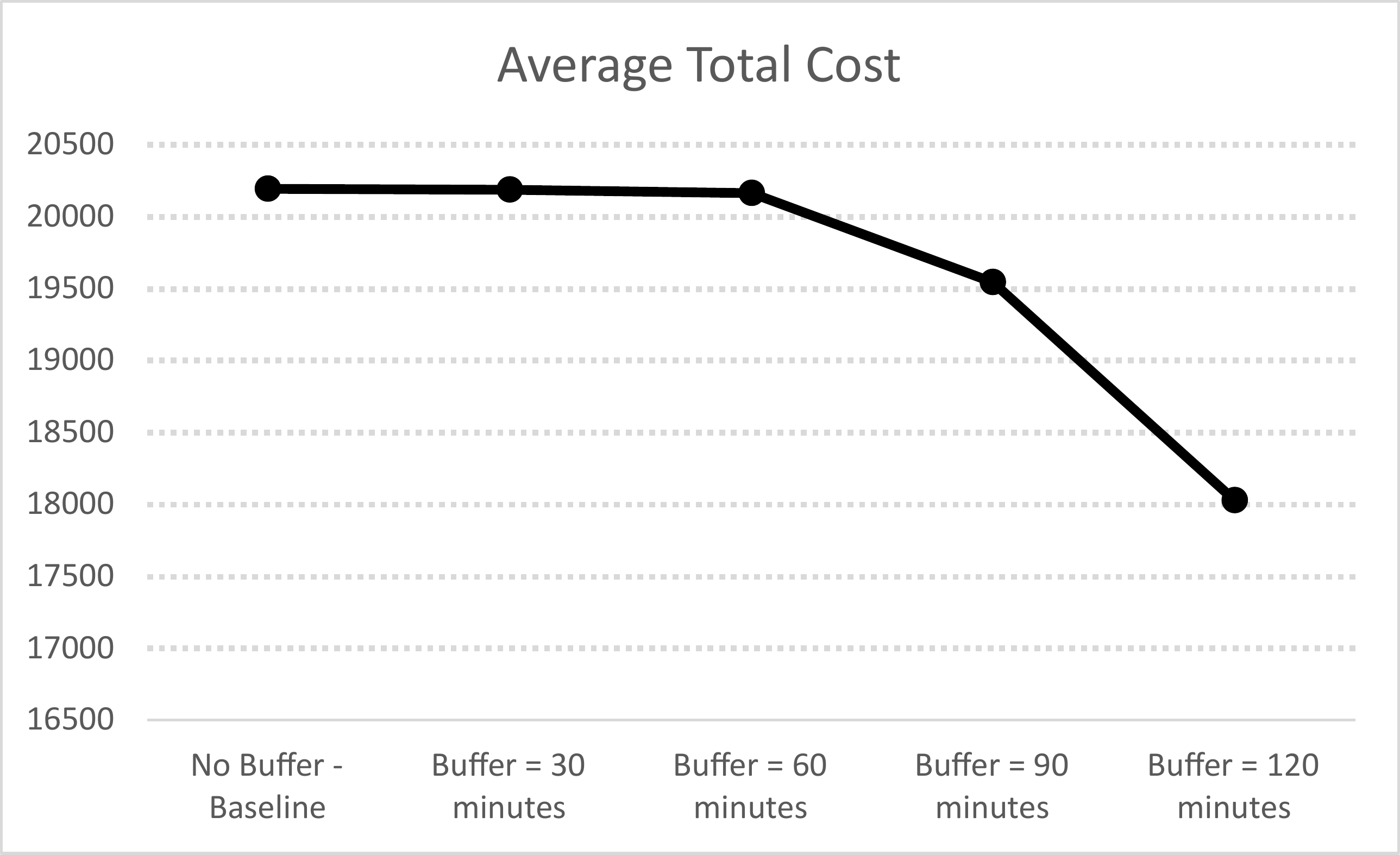}
    \caption{Total cost}
  \end{subfigure}\hfill
  \begin{subfigure}{0.48\textwidth}
    \includegraphics[width=\linewidth]{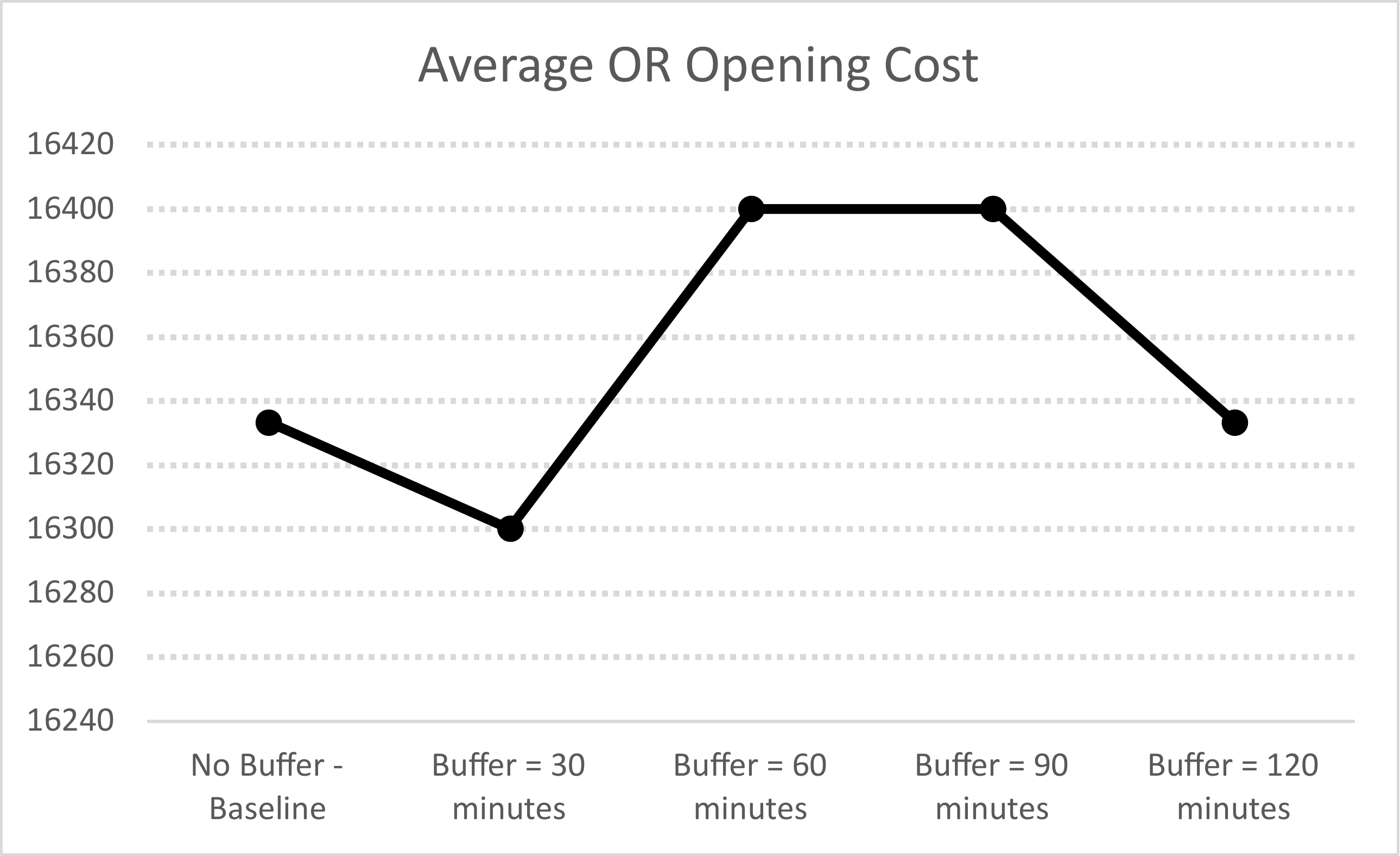}
    \caption{OR opening cost}
  \end{subfigure}

  \vspace{0.6em}

  \begin{subfigure}{0.48\textwidth}
    \includegraphics[width=\linewidth]{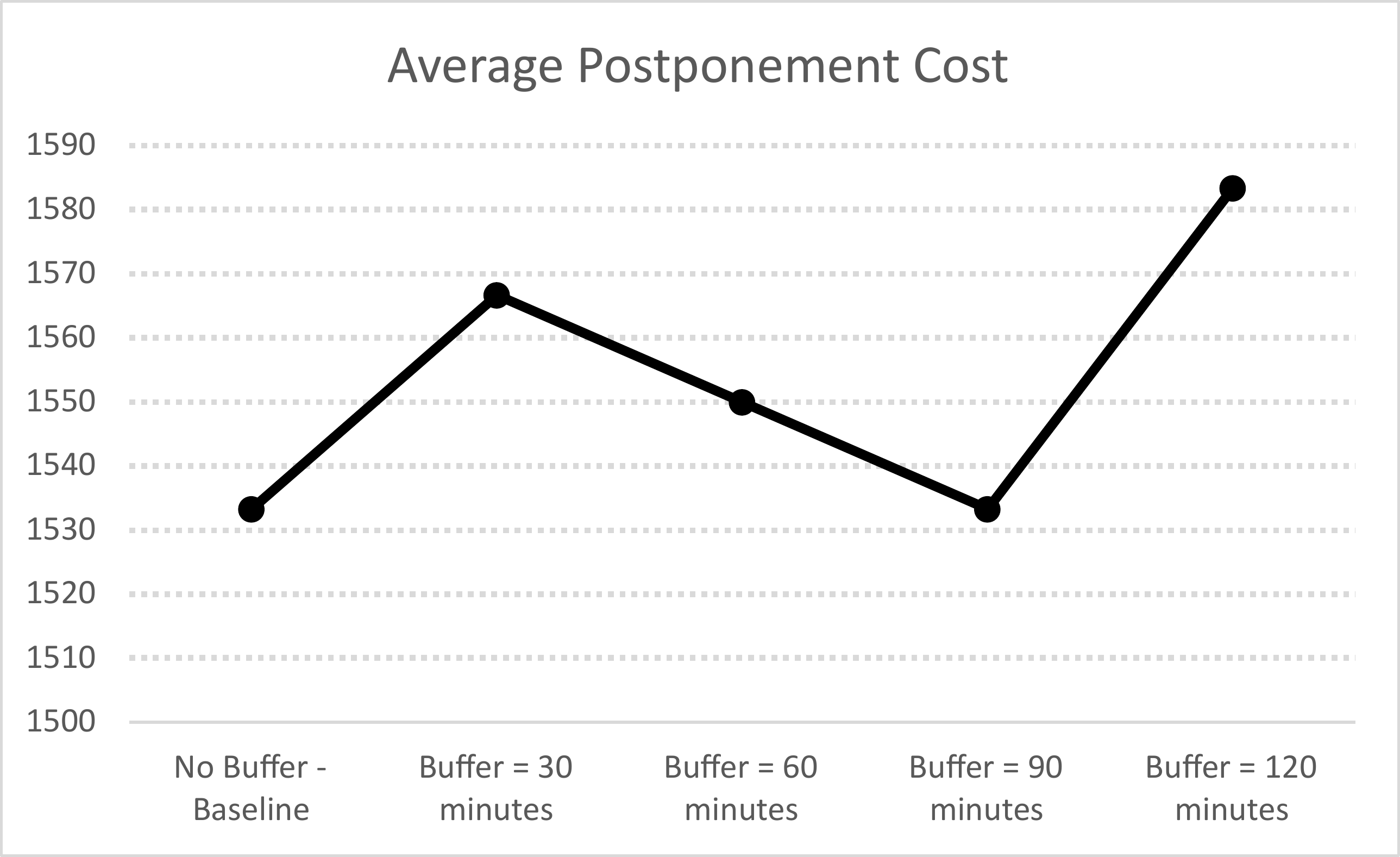}
    \caption{Postponement cost}
  \end{subfigure}\hfill
  \begin{subfigure}{0.48\textwidth}
    \includegraphics[width=\linewidth]{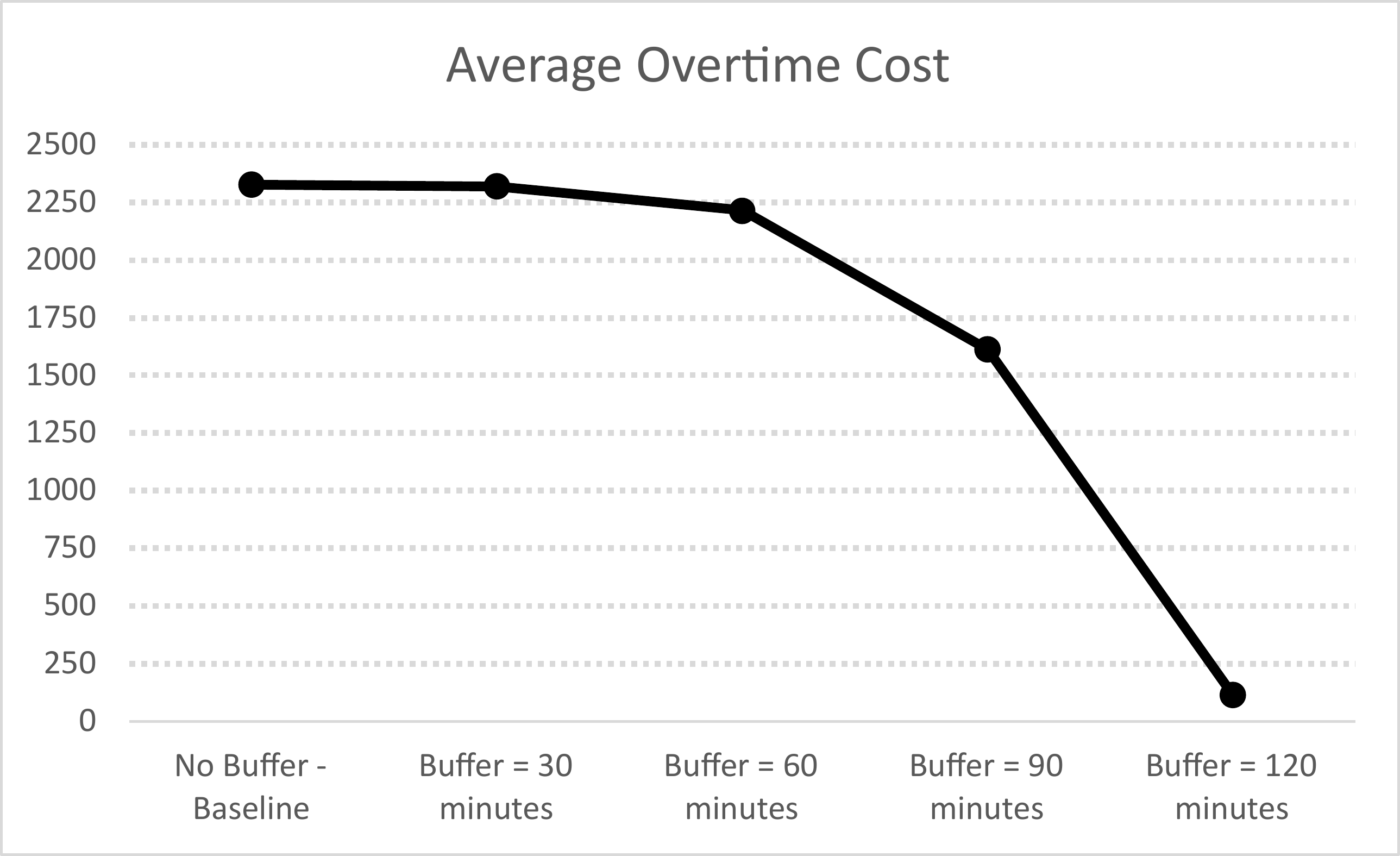}
    \caption{Overtime cost}
  \end{subfigure}
\end{minipage}
}
\caption{Impact of buffer-time policies for Naples instances}
\label{fig_naples_buffer}
\end{figure}

\paragraph{Real-world Instances}\label{Naples_Instances_Buffer_Analysis}
We replicated the buffer-time analysis on the real-world (Naples) instances. Figure \ref{fig_naples_buffer} shows the summary of the results, and the detailed results for each instance, is provided in Tables \ref{Naples_Instances_Buffer_Analysis_table1} and \ref{Naples_Instances_Buffer_Analysis_table2}. The mean total cost attains its minimum at $B=120$, overtime cost decreases monotonically as $B$ increases and overtime cost is nearly eliminated at $B=120$. OR-opening and postponement costs move in opposite directions as buffer grows. Notably, when increasing $B$ from 60 to 90 minutes, OR-opening cost remains essentially unchanged while postponement cost declines, indicating that the model exploits latent capacity within already opened OR-days to absorb variability. When moving from $B=90$ to $B=120$, postponement cost reaches its highest average level, OR-opening cost returns close to its baseline value, and overtime drops sharply toward zero which shows cost trade-off that favors buffering regular time over incurring overtime. To verify that these patterns are not driven by hard resource ceilings, we note (Section~\ref{Definition_Instances_Real}) that four ORs are available per day (24 OR-days over the horizon), and Table \ref{Naples_Instances_Buffer_Analysis_table2} shows a maximum OR-opening cost of 20{,}000. Therefore, the model could open all available OR-days when beneficial. Hence, the observed similar behavior of the model in both synthetic and Naples instances, reflects internal cost trade-offs rather than binding capacity limits. Considering the tested $\pm20\%$ duration variability, a 120-minute time buffer has been selected as the best policy in the Naples cases, delivering the lowest total cost by virtually eliminating overtime with only modest adjustments in OR-opening and postponements. \\

\begin{table}[h]
\caption{Hospital KPIs of buffer analysis for real‐world instances considering 30 and 60 minutes of buffers}
\resizebox{\textwidth}{!}{%
\begin{tabular}{ccccccccccccccccc}
\hline
\multirow{2}{*}{\begin{tabular}[c]{@{}c@{}}Instance\\ Name\end{tabular}} & \multirow{2}{*}{$\lvert P \lvert$} &  & \multicolumn{4}{c}{No Buffer - Baseline} &  & \multicolumn{4}{c}{Buffer = 30 minutes} &  & \multicolumn{4}{c}{Buffer = 60 minutes} \\ \cline{4-7} \cline{9-12} \cline{14-17} 
 &  &  & Total Cost & \begin{tabular}[c]{@{}c@{}}Postponement\\      Cost\end{tabular} & \begin{tabular}[c]{@{}c@{}}OR Opening\\      Cost\end{tabular} & \begin{tabular}[c]{@{}c@{}}Overtime\\      Cost\end{tabular} &  & Total Cost & \begin{tabular}[c]{@{}c@{}}Postponement\\      Cost\end{tabular} & \begin{tabular}[c]{@{}c@{}}OR Opening\\      Cost\end{tabular} & \begin{tabular}[c]{@{}c@{}}Overtime\\      Cost\end{tabular} &  & Total Cost & \begin{tabular}[c]{@{}c@{}}Postponement\\      Cost\end{tabular} & \begin{tabular}[c]{@{}c@{}}OR Opening\\      Cost\end{tabular} & \begin{tabular}[c]{@{}c@{}}Overtime\\      Cost\end{tabular} \\ \hline
INST1 & 143 &  & 20583.3 & 1000 & 17000 & 2583.3 &  & 20958.3 & 1000 & 17000 & 2958.3 &  & 20750 & 1000 & 17000 & 2750 \\
INST2 & 139 &  & 20333.3 & 1000 & 17000 & 2333.3 &  & 20541.7 & 1000 & 17000 & 2541.7 &  & 20458.3 & 1000 & 17000 & 2458.3 \\
INST3 & 164 &  & 22791.7 & 2000 & 18000 & 2791.7 &  & 22125 & 1500 & 19000 & 1625 &  & 23166.7 & 2000 & 19000 & 2166.7 \\
INST4 & 125 &  & 18833.3 & 0 & 16000 & 2833.3 &  & 18750 & 0 & 16000 & 2750 &  & 18916.7 & 0 & 16000 & 2916.7 \\
INST5 & 160 &  & 21791.7 & 500 & 18000 & 3291.7 &  & 21666.7 & 500 & 18000 & 3166.7 &  & 21208.3 & 1000 & 18000 & 2208.3 \\
INST6 & 146 &  & 21791.7 & 2000 & 17000 & 2791.7 &  & 21750 & 2000 & 17000 & 2750 &  & 22041.7 & 2000 & 17000 & 3041.7 \\
INST7 & 162 &  & 23833.3 & 3500 & 18000 & 2333.3 &  & 24833.3 & 3500 & 18000 & 3333.3 &  & 24541.7 & 3500 & 18000 & 3041.7 \\
INST8 & 136 &  & 20458.3 & 2500 & 16000 & 1958.3 &  & 20375 & 3000 & 15000 & 2375 &  & 19875 & 2500 & 16000 & 1375 \\
INST9 & 164 &  & 26125 & 4500 & 19000 & 2625 &  & 26083.3 & 4500 & 19000 & 2583.3 &  & 26333.3 & 4500 & 19000 & 2833.3 \\
INST10 & 155 &  & 21625 & 2500 & 17000 & 2125 &  & 22416.7 & 2500 & 17000 & 2916.7 &  & 22041.7 & 2500 & 17000 & 2541.7 \\
INST11 & 136 &  & 15875 & 1000 & 14000 & 875 &  & 16041.7 & 1000 & 14000 & 1041.7 &  & 16541.7 & 1000 & 14000 & 1541.7 \\
INST12 & 149 &  & 21541.7 & 500 & 18000 & 3041.7 &  & 20958.3 & 500 & 18000 & 2458.3 &  & 20583.3 & 500 & 18000 & 2083.3 \\
INST13 & 129 &  & 16458.3 & 500 & 14000 & 1958.3 &  & 15875 & 500 & 14000 & 1375 &  & 16250 & 500 & 14000 & 1750 \\
INST14 & 133 &  & 20125 & 2000 & 16000 & 2125 &  & 20291.7 & 3000 & 15000 & 2291.7 &  & 19750 & 3000 & 15000 & 1750 \\
INST15 & 137 &  & 18666.7 & 500 & 16000 & 2166.7 &  & 18833.3 & 500 & 16000 & 2333.3 &  & 18458.3 & 500 & 16000 & 1958.3 \\
INST16 & 141 &  & 20458.3 & 2500 & 16000 & 1958.3 &  & 20416.7 & 2500 & 16000 & 1916.7 &  & 20500 & 2500 & 16000 & 2000 \\
INST17 & 157 &  & 24583.3 & 3500 & 19000 & 2083.3 &  & 24708.3 & 3500 & 19000 & 2208.3 &  & 24666.7 & 3500 & 19000 & 2166.7 \\
INST18 & 162 &  & 25291.7 & 3000 & 19000 & 3291.7 &  & 25125 & 3000 & 19000 & 3125 &  & 25333.3 & 3000 & 20000 & 2333.3 \\
INST19 & 159 &  & 21125 & 1500 & 18000 & 1625 &  & 21041.7 & 1500 & 18000 & 1541.7 &  & 21666.7 & 1500 & 18000 & 2166.7 \\
INST20 & 137 &  & 19125 & 1500 & 15000 & 2625 &  & 18666.7 & 1500 & 15000 & 2166.7 &  & 19250 & 1500 & 15000 & 2750 \\
INST21 & 110 &  & 15750 & 1000 & 13000 & 1750 &  & 16208.3 & 1000 & 13000 & 2208.3 &  & 16375 & 1000 & 13000 & 2375 \\
INST22 & 149 &  & 23125 & 3000 & 17000 & 3125 &  & 23166.7 & 3000 & 17000 & 3166.7 &  & 22625 & 2000 & 18000 & 2625 \\
INST23 & 128 &  & 15916.7 & 500 & 14000 & 1416.7 &  & 16500 & 500 & 14000 & 2000 &  & 16708.3 & 500 & 14000 & 2208.3 \\
INST24 & 130 &  & 17750 & 1000 & 15000 & 1750 &  & 17583.3 & 1000 & 15000 & 1583.3 &  & 17291.7 & 1000 & 15000 & 1291.7 \\
INST25 & 122 &  & 16083.3 & 500 & 14000 & 1583.3 &  & 15750 & 500 & 14000 & 1250 &  & 15708.3 & 500 & 14000 & 1208.3 \\
INST26 & 139 &  & 19000 & 500 & 16000 & 2500 &  & 18875 & 500 & 16000 & 2375 &  & 18458.3 & 500 & 16000 & 1958.3 \\
INST27 & 142 &  & 19458.3 & 500 & 16000 & 2958.3 &  & 19041.7 & 500 & 16000 & 2541.7 &  & 18708.3 & 500 & 16000 & 2208.3 \\
INST28 & 132 &  & 16583.3 & 1000 & 14000 & 1583.3 &  & 16291.7 & 1000 & 14000 & 1291.7 &  & 16583.3 & 1000 & 14000 & 1583.3 \\
INST29 & 143 &  & 21291.7 & 1500 & 17000 & 2791.7 &  & 21500 & 1500 & 17000 & 3000 &  & 21625 & 1500 & 17000 & 3125 \\
INST30 & 137 &  & 19416.7 & 500 & 16000 & 2916.7 &  & 19250 & 500 & 16000 & 2750 &  & 18541.7 & 500 & 16000 & 2041.7 \\ \hline
\multicolumn{2}{c}{Average} &  & 20193.1 & 1533.3 & 16333.3 & 2326.4 &  & 20187.5 & 1566.7 & 16300 & 2320.8 &  & 20165.3 & 1550 & 16400 & 2215.3 \\ \hline
\end{tabular}%
}
\label{Naples_Instances_Buffer_Analysis_table1}
\end{table}

\begin{table}[H]
\caption{Hospital KPIs of buffer analysis for real‐world instances considering 90 and 120 minutes of buffers}
\resizebox{\textwidth}{!}{%
\tiny
\begin{tabular}{cccccccccccc}
\hline
\multirow{2}{*}{\begin{tabular}[c]{@{}c@{}}Instance\\ Name\end{tabular}} & \multirow{2}{*}{$\lvert P \lvert$} &  & \multicolumn{4}{c}{Buffer = 90 minutes} &  & \multicolumn{4}{c}{Buffer = 120 minutes} \\ \cline{4-7} \cline{9-12} 
 &  &  & Total Cost & \begin{tabular}[c]{@{}c@{}}Postponement\\      Cost\end{tabular} & \begin{tabular}[c]{@{}c@{}}OR Opening\\      Cost\end{tabular} & \begin{tabular}[c]{@{}c@{}}Overtime\\      Cost\end{tabular} &  & Total Cost & \begin{tabular}[c]{@{}c@{}}Postponement\\      Cost\end{tabular} & \begin{tabular}[c]{@{}c@{}}OR Opening\\      Cost\end{tabular} & \begin{tabular}[c]{@{}c@{}}Overtime\\      Cost\end{tabular} \\ \hline
INST1 & 143 &  & 20125 & 1000 & 17000 & 2125 &  & 18125 & 1000 & 17000 & 125 \\
INST2 & 139 &  & 19125 & 1000 & 17000 & 1125 &  & 18125 & 1000 & 17000 & 125 \\
INST3 & 164 &  & 22000 & 1500 & 19000 & 1500 &  & 20666.7 & 1500 & 19000 & 166.7 \\
INST4 & 125 &  & 17166.7 & 0 & 16000 & 1166.7 &  & 16166.7 & 0 & 16000 & 166.7 \\
INST5 & 160 &  & 20791.7 & 1000 & 18000 & 1791.7 &  & 19125 & 1000 & 18000 & 125 \\
INST6 & 146 &  & 20375 & 2000 & 17000 & 1375 &  & 19041.7 & 2000 & 17000 & 41.7 \\
INST7 & 162 &  & 24000 & 3500 & 18000 & 2500 &  & 21791.7 & 3500 & 18000 & 291.7 \\
INST8 & 136 &  & 19625 & 2500 & 16000 & 1125 &  & 18625 & 2500 & 16000 & 125 \\
INST9 & 164 &  & 25708.3 & 4500 & 19000 & 2208.3 &  & 24125 & 5000 & 19000 & 125 \\
INST10 & 155 &  & 21291.7 & 3000 & 17000 & 1291.7 &  & 19750 & 2500 & 17000 & 250 \\
INST11 & 136 &  & 16208.3 & 1000 & 14000 & 1208.3 &  & 15208.3 & 1000 & 14000 & 208.3 \\
INST12 & 149 &  & 20291.7 & 500 & 18000 & 1791.7 &  & 18541.7 & 500 & 18000 & 41.7 \\
INST13 & 129 &  & 16166.7 & 500 & 14000 & 1666.7 &  & 14666.7 & 500 & 14000 & 166.7 \\
INST14 & 133 &  & 19875 & 3000 & 15000 & 1875 &  & 18083.3 & 3000 & 15000 & 83.3 \\
INST15 & 137 &  & 18000 & 500 & 16000 & 1500 &  & 16666.7 & 500 & 16000 & 166.7 \\
INST16 & 141 &  & 19875 & 2500 & 16000 & 1375 &  & 18541.7 & 2500 & 16000 & 41.7 \\
INST17 & 157 &  & 25125 & 3000 & 20000 & 2125 &  & 22666.7 & 3500 & 19000 & 166.7 \\
INST18 & 162 &  & 23666.7 & 3000 & 19000 & 1666.7 &  & 22125 & 3000 & 19000 & 125 \\
INST19 & 159 &  & 21791.7 & 1500 & 18000 & 2291.7 &  & 19666.7 & 1500 & 18000 & 166.7 \\
INST20 & 137 &  & 18875 & 500 & 16000 & 2375 &  & 16541.7 & 1500 & 15000 & 41.7 \\
INST21 & 110 &  & 15250 & 1000 & 13000 & 1250 &  & 14083.3 & 1000 & 13000 & 83.3 \\
INST22 & 149 &  & 21416.7 & 3000 & 17000 & 1416.7 &  & 20041.7 & 3000 & 17000 & 41.7 \\
INST23 & 128 &  & 16166.7 & 500 & 14000 & 1666.7 &  & 14541.7 & 500 & 14000 & 41.7 \\
INST24 & 130 &  & 17166.7 & 1000 & 15000 & 1166.7 &  & 16000 & 1000 & 15000 & 0 \\
INST25 & 122 &  & 15583.3 & 500 & 14000 & 1083.3 &  & 14541.7 & 500 & 14000 & 41.7 \\
INST26 & 139 &  & 18000 & 500 & 16000 & 1500 &  & 16625 & 500 & 16000 & 125 \\
INST27 & 142 &  & 18041.7 & 500 & 16000 & 1541.7 &  & 16541.7 & 500 & 16000 & 41.7 \\
INST28 & 132 &  & 15583.3 & 1000 & 14000 & 583.3 &  & 15000 & 1000 & 14000 & 0 \\
INST29 & 143 &  & 21000 & 1500 & 17000 & 2500 &  & 18625 & 1500 & 17000 & 125 \\
INST30 & 137 &  & 18041.7 & 500 & 16000 & 1541.7 &  & 16583.3 & 500 & 16000 & 83.3 \\ \hline
\multicolumn{2}{c}{Average} &  & 19544.5 & 1533.3 & 16400 & 1611.1 &  & 18027.8 & 1583.3 & 16333.3 & 111.1 \\ \hline
\end{tabular}%
}
\label{Naples_Instances_Buffer_Analysis_table2}
\end{table}

\subsection{Implications for Healthcare Operations Management}\label{Implications_Operations_Management}
The numerical results presented in the previous section demonstrate the superiority of the PMIORPS model and the RGA-CG solution methodology. In this section, we discuss how the proposed model can be applied in real-world healthcare settings, particularly in hospitals and surgical centers.
The PMIORPS model is specifically designed for IORPS problems, which fall under the scope of operational healthcare management. Given that operational decision-making typically addresses short-term planning, PMIORPS is well-suited for such scenarios. As shown in Table \ref{tab_real_1}, the model is capable of efficiently handling planning problems involving up to 200 surgeries, four operating rooms (ORs), and a six-day planning horizon. For problem sizes within these bounds, we recommend solving the model directly using commercial solvers such as IBM CPLEX, as they provide high-quality solutions within a reasonable computation time. However, for larger instances that exceed these thresholds, the PMIORPS model should be combined with the RGA-CG approach. This hybrid solution methodology enhances scalability and provides near-optimal results for more complex planning scenarios.

Another important application of this work lies at the tactical decision-making level, where planning horizons typically span two to four weeks. In such cases, as demonstrated in Table \ref{tab_real_2}, the RGA-CG method is particularly effective for producing high-quality solutions. A recommended strategy for hospital managers is to adopt a two-stage planning approach: use the RGA-CG method to solve the tactical-level problem and assign surgeries to specific weeks, and then apply commercial solvers like CPLEX to solve each week’s operational subproblem independently. This integrated approach enables healthcare administrators to coordinate tactical and operational planning simultaneously, improving both resource utilization and responsiveness to scheduling demands.

Although the model is formulated for deterministic elective surgery planning and scheduling, the analyses in Sections \ref{Emergency_Surgery_Disruption_Analysis} and \ref{Robustness_Duration_Variability_Buffer_Analysis} demonstrate that it reliably accommodates disruption from emergency arrivals and duration variability with limited deviation from the initial plan. In practice, the effectiveness of these responses depends critically on input quality. Therefore, hospitals should invest in accurate, routinely recalibrated prediction and estimation models of surgery duration and clinically priority scores (with uncertainty estimates) to adjust scheduling and buffer sizing. With credible forecasts in place, managers can select buffer-time policies that mitigate variability (e.g., the 120-minute policy identified in our tests) while maintaining high utilization of key resources (ORs and surgeons) and controlling overtime, postponements, and OR-opening costs. Institutionalizing a short freeze window (e.g., Days~1–2) and daily re-optimization within a rolling horizon further enables timely schedule updates, making the proposed model and solution approach a robust, operationally viable tool for day-to-day theatre management aligned with hospital KPIs.

Even though the proposed model is designed for day-by-day advance planning and scheduling, its time-indexed decision variable $x_{idt}$ (with a 5-minute discretization) enables a natural adaptation to dynamic, within-day surgery scheduling. When the horizon is restricted to a single day and operational decisions fixed from upstream planning (e.g., the number of open ORs and preassigned resources), the model dimension shrinks substantially since many variables and constraints become fixed or redundant, so that reoptimization can be performed rapidly. Consequently, following the procedure outlined in Section~\ref{Emergency_Surgery_Disruption_Analysis}, the same formulation can serve as an online decision support tool. Managers can update estimated surgery durations, incorporate new emergency arrivals, and reoptimize on a rolling basis to generate feasible, and responsive schedules with minimal disruption to previously fixed assignments. Beyond this, the formulation and its associated time horizon are well suited to exploit predictive information in a lookahead dynamic policy. In practice, time-sensitive forecasts, such as expected surgery durations, or emergency arrivals, could be readily integrated to anticipate disruptions and proactively adjust schedules. This predictive integration further enhances the model’s applicability to operational control by supporting near–real-time adjustments while maintaining schedule stability and coordinated use of resources. In this way, the methodological developments presented here naturally lend themselves to decision-support tools for operational settings where flexibility and responsiveness are critical.

Our approach is designed for short operational horizons (less than 20 working days), since surgeon availability, patient status, and OR capacity change frequently. Also, longer horizons are typically handled via block scheduling. Our proposed model targets typical hospital scales (few ORs,  less than 200 surgeries/week), as substantially larger systems are uncommon. Finally, while instances up to 480 surgeries are solved efficiently for the extended time horizon using RGA-CG, 800-surgery, 40-day cases exceeded memory on our hardware. Therefore, handling such extremes would require parallel/distributed computation and memory-optimized data structures.

\section{Conclusion and Future Research}\label{Conclusion}
Operating rooms produce a significant portion of the hospital’s income, and OR utilization has a direct effect on the hospital’s cost management. In this paper, we focused on an integrated operating room management problem that addresses the planning and scheduling problems simultaneously. We developed a new mixed integer programming formulation for this problem that has fewer constraints and variables in comparison to the commonly used baseline mathematical model in the literature. We used two sets of synthetic and real-world instances for our analysis. Solving our developed MIP model yields an average optimality gap of 1.23\% for synthetic scenarios, outperforming the current state-of-the-art model that resulting an average optimality gap of 17.69\%. 
An analysis of the real-world instances sourced from a local hospital in Naples demonstrates the strong performance of our proposed model, which achieved an average optimality gap of just 1.49\%. In contrast, the baseline MIORPS model was unable to find any feasible solution with a GAP\% below 100 in 28 out of the 30 instances, indicating a significant shortcoming in both solution quality and robustness. These results clearly highlight the superior efficiency and reliability of our approach when applied to complex, real-world healthcare scheduling scenarios.
For solving large-scale instances, we went one step further and developed a new reinforcement-learning-based column algorithm that relies on the proposed mathematical model. In this algorithm, named RGA-CG algorithm, we integrated reinforcement learning and genetic algorithms for generating the initial columns of the CG algorithms. The results show the RGA-CG is a powerful method to find a feasible solution in a real-world setting significantly outperforms an existing branch-and-cut algorithm from the literature.
In this study, we proposed a model that, while not explicitly encoding emergency surgeries, integrates naturally with a rolling-horizon scheme to accommodate on-the-fly rescheduling triggered by emergency arrivals and priority escalations. Computational experiments indicate that the incremental cost of such rescheduling can be kept as low as 4.13\%, and that priority-escalation adjustments are even less costly, at roughly 0.5\%. We further examined robustness to surgery duration uncertainty and found that when surgery durations vary by up to $\pm20\%$ from their estimates, a 120-minute buffer suffices to absorb the imposed variability.

For future research, it would be interesting to consider that emergency arrivals surgeries could happen anytime during the planning horizon. Therefore, studying an OR rescheduling or replanning-rescheduling problem integrated with PMIORPS would be a practical extension of this study. Another practical extension could be to study the stochastic version of PMIORPS where the durations of surgeries are not deterministic. In fact, based on the type of chosen uncertainty, stochastic optimization methods such as sample average approximation or robust optimization methods could be investigated.

\vspace{20 pt}

\noindent
{\large \bf{Data Availability Statement}}\\[-10pt]

\noindent
The data that support the findings of this study are available upon request.

\noindent
{\large \bf{Conflict of Interests Statement}}\\[-10pt]

\noindent
The authors declare that they have no conflict of interest.

\printbibliography

\newpage
\appendix
\setcounter{page}{1}
\section*{Appendix}\label{Appendix}
\section{Proof of Theorem \ref{Tr_1}}\label{AppendixA}
In explaining the proof of Theorem \ref{Tr_1}, we introduce a new set, a supplementary variable, and an alternative mathematical model named APMIORPS. Denoted as $T'_i$, the set encompasses time slots wherein surgery $i$ is potentially underway. The binary variable $x'_{idkt}$ takes a value of one if surgery $i$ is being performed during time slot $t$ within operating room $k$ on day $d$.

\begin{flalign}\label{Eq_AA0}
 \textbf{APMIORPS: }\;\; \min \sum_{i \in I_2} c^{Pos} z_i + \sum_{d \in D} c^{OR} y_d + \sum_{i \in I}\sum_{d \in D:d \leq d_i}\sum_{t \in T_i} c_{it} x_{idt} &&
\end{flalign}

\noindent \mbox{Subject to:}
\begin{alignat}{2}
& Constraints \: (\ref{Eq_28}) \; to \;(\ref{Eq_33})   && \nonumber \\
& \sum_{t' \in T_i:\{t' \leq t \;and\; t'+t_i>t\}} x_{idt'} = \sum_{k \in K} x'_{idkt}   \;\;                                   && \forall\, i \in I, \forall\, d \in D, \forall\, t \in T'_i\label{Eq_AA1}\\
& x'_{idk(t-1)}-x_{id(t-t_{i})} \leq x'_{idkt}                                   && \forall\, i \in I, \forall\, d \in D, \forall\, k \in K, \forall\, t \in T'_i:\{t \geq 2 \;,\; t-t_{i}>0\}\label{Eq_AA2}\\
& x'_{idk(t-1)} \leq x'_{idkt}                                   && \forall\, i \in I, \forall\, d \in D, \forall\, k \in K, \forall\, t \in T'_i:\{t \geq 2 \;,\; t-t_{i}\leq0\}\label{Eq_AA3}\\
& x'_{idkt} \in \{0,1\}                                                                && \forall\, i \in I, \forall\, d \in D, \forall\, k \in K, \forall\, t \in T'_i \nonumber
\end{alignat}

our proof relies on the following two steps:
\begin{enumerate}
    \item Step 1: In this initial stage, we prove the APMIORPS model is a valid representation of the MIORPS model.
    \item Step 2: Subsequently, we demonstrate in the APMIORPS model that for any feasible solution of $x_{idt}$ obtained from constraints (\ref{Eq_28}) to (\ref{Eq_33}), there always exists a feasible solution for $x'_{idkt}$ variables such that Constraints (\ref{Eq_AA1}) to (\ref{Eq_AA3}) are always satisfied.
\end{enumerate}

In the following, we provide the details of the above-mentioned two steps.\\

\textbf{Step1}:\\
As evidenced in the APMIORPS model, the variable $x'_{idkt}$ solely features within the newly introduced constraints and not within the objective function. Constraint (\ref{Eq_AA1}) ensures that if a surgery commences before time slot $t$ and extends beyond $t$ on day $d$, it must be fulfilled within an operating room during time slot $t$. For instance, consider surgery 1 scheduled to start at the beginning of time slot 3 on day 2, with a duration spanning 5 time slots. Accordingly, $x_{123}=1$. At time slot 3, the left-hand side of constraint (\ref{Eq_AA1}) equals 1, necessitating the right-hand side to also equal one for one of the operating rooms $k \in K$. This condition holds true for all time slots during which surgery 1 is ongoing, from time slots 3 to 7. However, the APMIORPS model remains invalid as the $x'_{idkt}$ variables could be assigned to different operating rooms for the same surgery. Constraints (\ref{Eq_AA2}) and (\ref{Eq_AA3}) rectify this by stipulating that all surgeries must be carried out within the same operating room for the entire duration of the surgery. Thus, the APMIORPS model now stands as a valid reformulation of the MIORPS model.

\textbf{Step2}:\\
If we could show that the constraints (\ref{Eq_AA1}) to (\ref{Eq_AA3}) are always satisfied regardless of the values for $x_{idt}$, we can remove the constraints from the APMIORPS model since $x'_{idkt}$ not appeared in the objective function, and we will have PMIORPS model.
Random binary variables can be allocated to $x'_{idkt}$ such that constraint (\ref{Eq_AA1}) aligns with the constraints (\ref{Eq_28}) to (\ref{Eq_33}) outlined in this paper. Subsequently, regardless of the initial values of $x'_{idkt}$, they can be adjusted using Algorithm \ref{Alg_1} to ensure compliance with both constraints (\ref{Eq_AA2}) and (\ref{Eq_AA3}).

In Algorithm \ref{Alg_1}, we utilize $\gamma_i$ as a flag variable, which is set to 1 if the necessary adjustments for surgery $i$ are made. Initially, in the first line, all $\gamma_i$ variables are initialized to 0. Subsequently, the modification process commences from the last available time slot to the first one, essentially operating in reverse (line 2). Lines 3 and 4 encompass loops over all operating rooms and surgeries. In line 5, we ascertain whether modifications for surgery $i$ have already been conducted, proceeding only if it marks the first instance of modification. Line 6 involves verifying whether $x'_{idkt}$ equals one or not. It's worth noting that since the algorithm operates in reverse, the initial occurrence of $x'_{idkt}$ being one corresponds to the final time slot during which surgery $i$ performs in operating room $k$ on day $d$. 
In line 7, we loop over the auxiliary parameter $t'$ which is used to make $(t_i-1)$ modifications for surgery $i$.

In line 8, $i^*$ denotes the index of the surgery allocated to the preceding time slot within operating room $k$ on day $d$. Similarly, $k^*$ signifies the index of the operating room where surgery $i$ was performed during the previous time slot. Lines 10 to 13 encompass a process wherein we exchange values for various $z$ variables in the time slot $(t-t')$ to uphold the integrity of the solution to satisfy the constraints (\ref{Eq_AA2}) and (\ref{Eq_AA3}). Line 10 indicates that surgery $i$ should not be performed in operating room $k^*$ during time slot $(t-t')$, ensuring it remains in the same operating room as in time slot $t$, as enforced by line 11. In the 12th line, we specify that the surgery $i^*$ must not be performed in operating room $k$ during time slot $(t-t')$ since it is already occupied by surgery $i$. Line 13 assigns surgery $i^*$ to operating room $k^*$ during time slot $(t-t')$ on day $d$. Upon completing the necessary changes for surgery $i$, we set $\gamma_i$ to 1 in line 15. 

By implementing the aforementioned adjustments across all time slots, operating rooms, and surgeries, constraints (\ref{Eq_AA2}) and (\ref{Eq_AA3}) are effectively satisfied. Thus, as previously noted, given that $x'_{idkt}$ solely appears in the additional constraints of the PMIORPS model and not within the objective functions, and considering that all three constraints (\ref{Eq_AA1}) to (\ref{Eq_AA3}) are invariably fulfilled for any feasible solution of $x_{idt}$, constraints (\ref{Eq_AA1}) to (\ref{Eq_AA3}) can be omitted from the APMIORPS model. Consequently, since the APMIORPS model serves as a valid formulation of the MIORPS model, the PMIORPS model likewise will be a valid formulation of the MIORPS model.

\begin{algorithm}
\caption{$x'_{idkt}$ modification algorithm}\label{Alg_1}
\begin{algorithmic}[1]
\State $Set \; all \; \gamma_i = 0$
\For{($t=|T| \;to\; 1$)} \Comment{moving backward from $t=|T|$ to $t=1$}
  \For{($k=1 \;to\; |K|$)}
    \For{($i=1 \;to\; |I|$)}
      \If{$\gamma_i=0$}
        \If{$x'_{idkt}=1$}
          \For{($t'=1 \;to\; t_i-1$)}
             \State $i^*:=arg\{ a \in I : x'_{adk(t-t')}=1\}$ 
             \State $k^*:=arg\{ a \in K : x'_{ida(t-t')}=1\}$ 
             \State $x'_{idk^*(t-t')}:=0$                     
             \State $x'_{idk(t-t')}:=1$                       
             \State $x'_{i^*dk(t-t')}:=0$                     
             \State $x'_{i^*dk^*(t-t')}:=1$                  
          \EndFor 
          \State $\gamma_i=1$
        \EndIf
      \EndIf
    \EndFor
  \EndFor
\EndFor
\end{algorithmic}
\end{algorithm}

Table \ref{tab_APMIORPS} presents the computational comparison between the APMIORPS and PMIORPS models, based on the experimental setup described in Section \ref{Computational_Experiments}. The numerical results clearly indicate the superior performance of the PMIORPS model.
Notably, the average number of variables in APMIORPS, as reported in the “\# Vars” column, is more than eight times that of PMIORPS. Similarly, the average number of constraints, presented in the “\# Cons” column, is over 52 times higher for APMIORPS compared to PMIORPS. This substantial increase in model size for APMIORPS significantly limits its computational performance.
As a result of this complexity, none of the test instances could be solved to optimality using APMIORPS. On average, the final optimality gap for APMIORPS exceeds 53\%. In contrast, the PMIORPS model demonstrates significantly enhanced solvability, and achieved an average gap of only 1.23\%.
This significant difference in solution quality and tractability underscores the effectiveness of the PMIORPS formulation. By reducing model size and improving numerical stability, PMIORPS facilitates efficient solution procedures and enables practical applicability to real-world problem instances.

\begin{sidewaystable}
\doublespacing
\caption{Evaluation of the APMIORPS model, and PMIORPS model}
\resizebox{\textwidth}{!}{%
\tiny
\begin{tabular}{ccccccccccccccccccc}
\hline
\multirow{2}{*}{\begin{tabular}[c]{@{}c@{}}$\lvert D \lvert$\end{tabular}} &  & \multicolumn{8}{c}{APMIORPS} &  & \multicolumn{8}{c}{PMIORPS} \\ \cline{3-10} \cline{12-19} 
 &  & IS\% & OSI\% & LB & UB & Gap\% & \# Vars. & \# Cons. & Time &  & IS\% & OSI\% & LB & UB & Gap\% & \# Vars. & \# Cons. & Time \\ \hline
40 &  & 100 & 0 & 8645.832 & 8800 & 1.79 & 139097.2 & 146986.2 & 3600 &  & 100 & 40 & 8738.91 & 8800 & \textbf{0.67} & 16984 & 5530 & 2369.42 \\
60 &  & 100 & 0 & 14183.34 & 14700 & 3.76 & 214984.4 & 224546.4 & 3600 &  & 100 & 20 & 14206.52 & 14400 & \textbf{1.37} & 25934 & 5550 & 2900 \\
80 &  & 100 & 0 & 19329.16 & 20900 & 8.13 & 280399.2 & 291341.8 & 3600 &  & 100 & 0 & 19341.66 & 19666.66 & \textbf{1.69} & 33704.4 & 5570 & 3600 \\
100 &  & 100 & 0 & 23513.8 & 46016.68 & 96.43 & 359505 & 371761.4 & 3600 &  & 100 & 0 & 23533.34 & 23700 & \textbf{0.7} & 43404.6 & 5590 & 3600 \\
120 &  & 100 & 0 & 29510.44 & 75200 & 154.93 & 425461.6 & 439229.6 & 3600 &  & 100 & 0 & 29737.5 & 30250 & \textbf{1.71} & 51136 & 5610 & 3600 \\ \hline
Ave. &  & 100 & 0 &  &  & 53 & 283889.5 & 294773.1 & 3600 &  & 100 & \textbf{12} &  &  & \textbf{1.23} & \textbf{34232.6} & \textbf{5570} & \textbf{3213.89} \\ \hline
\end{tabular}%
}
\label{tab_APMIORPS}
\end{sidewaystable}

\newpage
\section{Proposed two-phase RGA algorithm}\label{AppendixB}
\begin{enumerate}
   \item Read instance and calculate $c_{it}$
   \begin{enumerate}
     \item Read instance data from the file
     \item Calculate $c_{it}$ based on the formulation given in Table \ref{tab_1} 
   \end{enumerate}
   
   \item RGA Phase1
   \begin{enumerate}
     \item Initialize parameters
     \begin{enumerate}
       \item Initialize GA parameters (nPop, nGen, pMu, Selection Method, Termination conditions)
       \item Initialize Reinforcement Learning parameters (RLA\_algorithm, L\_Reward, F\_Reward)
     \end{enumerate}
     
    \item Initialize RGA Phase1 population
    \begin{enumerate}
       \item Create random permutation of mandatory surgeries
       \item Create random permutation of non-mandatory surgeries
       \item Evaluate the uniqueness of the created permutations (lists) based on evaluated solution list
       \item Calculate $x_{idt}$, $y_d$, $z_i$ matrices according to the surgery lists for unique combinations
       \item Calculate fitness value
       \item Add the individual to the population and to the evaluated solution list
       \item Update the Best solution $x_{idt}^B$, $y_d^B$, $z_i^B$ if a better fitness value found
       \item Add a new column to the columns' pool if a better fitness value found
     \end{enumerate}

    \item For Gen=1 to Gen= nGen:
    \begin{enumerate}
       \item Terminate the for loop if the termination condition met 
       \item For nc=1 to nCr:
       \begin{enumerate}
         \item Select parent indexes
         \item Perform crossover separately on mandatory or non-mandatory lists and the updated probability rates for each operator based on RLA
         \item Perform mutation separately on mandatory or non-mandatory lists of each child considering pMu and the updated probability rates for each operator based on RLA
         \item Evaluate the uniqueness of both children based on the evaluated solution list
         \item Calculate $x_{idt}$, $y_d$, $z_i$ matrices according to the surgery lists for unique children
         \item Calculate fitness values for unique children
         \item Add the child to the population and to the evaluated solution list
         \item Update the Best solution $x_{idt}^B$, $y_d^B$, $z_i^B$ if a better fitness values found
         \item Add a new column to the columns' pool if better fitness values found
         \item Update the accumulated reward values for each crossover and mutation method
       \end{enumerate}
       
       \item For nm=1 to nMu:
       \begin{enumerate}
         \item Select candidate for mutation
         \item Perform mutation separately on mandatory or non-mandatory lists considering the updated probability rates for each operator based on RLA
         \item Evaluate the uniqueness of the mutated candidate
         \item Calculate$ x_{idt}$, $y_d$, $z_i$ matrices according to the surgery lists for unique candidate
         \item Calculate fitness values for a unique candidate
         \item Add the candidate to the population and to the evaluated solution list
         \item Update the Best solution $x_{idt}^B$, $y_d^B$, $z_i^B$ if a better fitness values found
         \item Add a new column to the columns' pool if a better fitness value found
         \item Update the accumulated reward values for mutation methods
       \end{enumerate}

       \item Sort the population based on the descending order of fitness values
       \item Trim the population (keep the first nPop individuals)
     \end{enumerate}
\end{enumerate}

   \item RGA Phase1 enhancement
   \begin{enumerate}
    \item For np=1 to nPop:
    \begin{enumerate}
       \item Select the np-th individual from the population of the RGA Phase1
       \item Calculate $x_{idt}$, $y_d$, $z_i$ matrices according to the surgery lists of each individual
       
       \item For d=1 to D:
       \begin{enumerate}
         \item Create $y^\prime$ matrix where$ y^\prime=y$
         \item Set $y_d^\prime=y_d-1$
         \item Calculate new $x_{idt}^{\prime\prime}$, $y_d^{\prime\prime}$, $z_i^{\prime\prime}$ matrices considering  $y^\prime$ as an upper bound for the total available ORs on day $d$
         \item Calculate fitness value
         \item Update the Best solution $x_{idt}^B$, $y_d^B$, $z_i^B$ if a better fitness values found
         \item Update the stored fitness value for the current individual if a better fitness value found
         \item Add a new column to the columns' pool if a better fitness value found
       \end{enumerate}
     \end{enumerate}
    \end{enumerate}

   \item RGA Phase2
   \begin{enumerate}
    \item Initialize RGA Phase2 population
    \begin{enumerate}
       \item For np=1 to nPop:
       \begin{enumerate}
         \item Select the np-th individual from the population of the RGA Phase1
         \item Create an integrated surgery list by joining mandatory and non-mandatory surgery lists together
         \item Add the individual to the RGA Phase2 population with the same fitness value
       \end{enumerate}
     \end{enumerate}

    \item For Gen=1 to Gen= nGen:
    \begin{enumerate}
       \item Terminate the for loop if the termination condition met 
       \item For nc=1 to nCr:
       \begin{enumerate}
         \item Select parent indexes
         \item Perform crossover separately on mandatory or non-mandatory lists and the updated probability rates for each operator based on RLA
         \item Perform mutation separately on mandatory or non-mandatory lists of each child considering pMu and the updated probability rates for each operator based on RLA
         \item Evaluate the uniqueness of both children based on the evaluated solution list
         \item Calculate $x_{idt}$, $y_d$, $z_i$ matrices according to the surgery lists for unique children
         \item Calculate fitness values for unique children
         \item Add the child to the population and to the evaluated solution list
         \item Update the Best solution $x_{idt}^B$, $y_d^B$, $z_i^B$ if a better fitness values found
         \item Add a new column to the columns' pool if better fitness values found
         \item Update the accumulated reward values for each crossover and mutation method
       \end{enumerate}
       
       \item For nm=1 to nMu:
       \begin{enumerate}
         \item Select candidate for mutation
         \item Perform mutation considering the updated probability rates for each operator based on RLA
         \item Evaluate the uniqueness of the mutated candidate
         \item Calculate $x_{idt}$, $y_d$, $z_i$ matrices according to the surgery lists for unique candidate
         \item Calculate fitness values for a unique candidate
         \item Add the candidate to the population and to the evaluated solution list
         \item Update the Best solution $x_{idt}^B$, $y_d^B$, $z_i^B$ if a better fitness values found
         \item Add a new column to the columns' pool if a better fitness value found
         \item Update the accumulated reward values for mutation methods
       \end{enumerate}

       \item Sort the population based on the descending order of fitness values
       \item Trim the population (keep the first nPop individuals)
     \end{enumerate}
\end{enumerate}

   \item RGA Phase2 enhancement
   \begin{enumerate}
    \item For np=1 to nPop:
    \begin{enumerate}
       \item Select the np-th individual from the population of the RGA Phase2
       \item Calculate $x_{idt}$, $y_d$, $z_i$ matrices according to the surgery lists of each individual
       
       \item For d=1 to D:
       \begin{enumerate}
         \item Create $y^\prime$ matrix where$ y^\prime=y$
         \item Set $y_d^\prime=y_d-1$
         \item Calculate new $x_{idt}^{\prime\prime}$, $y_d^{\prime\prime}$, $z_i^{\prime\prime}$ matrices considering  $y^\prime$ as an upper bound for the total available ORs on day $d$
         \item Calculate fitness value
         \item Update the Best solution $x_{idt}^B$, $y_d^B$, $z_i^B$ if a better fitness values found
         \item Update the stored fitness value for the current individual if a better fitness value found
         \item Add a new column to columns pool if a better fitness values found
       \end{enumerate}
     \end{enumerate}
    \end{enumerate}

   \item Report the $x_{idt}^B$, $y_d^B$, $z_i^B$ and the corresponding fitness value as the best solution
   \item Construct the relaxed M-CG where $x_j \in (0,1)$ and calculate the reduced cost
   \item Iteratively solve the S-CG and M-CG until termination condition(s) met
   \item Solve the M-CG with binary $x_j \in \{0,1\}$ variables and report the solution as the best found upper bound of the problem

\end{enumerate}

\section{BCMIORPS mathematical model}\label{AppendixC}
Like the IORPS, we modified the proposed Branch-and-Cut algorithm which is referred to as BCMIORPS to be compatible with the new cost minimization objective function. Table \ref{tab_9} shows the new notation of the BCMIORPS. Roshanaei et al. \cite{roshanaei2021solving} proposed a Branch-and-Cut algorithm to solve a large-scale IORPS. Like the IORPS, we modified the proposed Branch-and-Cut algorithm which is referred to as BCMIORPS to be compatible with the new cost minimization objective function.

In the BCMIORPS, the MIORPS is decomposed into the master problem and subproblem. The planning phase will be tackled in the master problem, While sequencing and scheduling will be handled and validated in the subproblem. For the suboptimal or infeasible subproblems, the proper cut will be added to the master problem which will be discussed. The Master problem of BCMIORPS is as follows:

\begin{table}[H]
\centering
\small
\caption{New notations for the BCMIORPS} 
\begin{tabular}{R{1cm} L{13cm}}
\hline
                     &                \\
\multicolumn{2}{l}{New sets:}       \\
  $\hat{I}_d:$       &  Set of patients that assigned to day $d$              \\
                     &                \\

\multicolumn{2}{l}{New variables:} \\
  $x_{idk}:$         &  1 if surgery $i$ assigned to operating room $k$ on day $d$, 0 otherwise              \\
  $o_{dk}:$          &  1 if operating room $k$ opened on day $d$, 0 otherwise              \\
  $\alpha_d:$        &  Auxiliary variable to capture the overtime cost of the day $d$               \\
  $\beta_d:$         &  Auxiliary variable to capture the subproblem cost of the day $d$              \\
  $x_i:$             &  0 if we could successfully schedule surgery $i$, 1 otherwise               \\
  $c_i:$             &  Completion time of surgery $i$                \\
  $y_{ii^\prime}:$   &  1 if surgery $i$ scheduled after surgery $i^\prime$, 0 otherwise               \\
                     &                \\
\hline
\end{tabular}
\label{tab_9}
\end{table}

\begin{flalign}\label{Eq_8}
 \textbf{MIORPS: }\;\; \min \sum_{i \in I_2} c^{Pos} z_i + \sum_{d \in D}\sum_{k \in K_d} c^{OR} y_{dk} + \sum_{d \in D} \alpha_d &&
\end{flalign}

\noindent \mbox{Subject to:}
\begin{alignat}{2}
& \sum_{d \in D:d \leq d_i}\sum_{k \in K_d} x_{idk} = 1                                             && \forall\, i \in I_1 \label{Eq_9}\\
& \sum_{d \in D:d \leq d_i}\sum_{k \in K_d} x_{idk} + z_i = 1                                       && \forall\, i \in I_2 \label{Eq_10}\\
& \sum_{i \in I:d \leq d_i} t_{i} x_{idk} \leq |T_1| y_{dk} + o_{dk} \qquad \qquad \qquad \qquad    && \forall\, d \in D, \forall\, k \in K_d \label{Eq_11}\\
& \sum_{i \in I^{t^\prime}_{l}:d \leq d_i} t_{i} x_{idk} \leq a_{ld}                                && \forall\, l \in L, \forall\, d \in D \label{Eq_12}\\
& \alpha_d \geq \sum_{k \in K_d} c^{Ovt} o_{dk}                                                     && \forall\, d \in D \label{Eq_13}\\
& o_{dk} \leq  |T_2| y_{dk}                                                                         && \forall\, d \in D, \forall\, k \in K_d \label{Eq_14}\\
& x_{idk}, y_{dk}, x_i \in \{0,1\}                                                                  && \forall\, i \in I, \forall\, d \in D:d \leq d_i, \forall\, k \in K_d \nonumber \\
& o_{dk}, \alpha_d, \beta_d \geq 0                                                                  && \forall\, i \in I, \forall\, d \in D:d \leq d_i, \forall\, k \in K_d \nonumber \\
& \mbox{Feasibility and optimality cuts}                                                            &&  \nonumber
\end{alignat}

The BCMIORPS master problem is aimed at minimizing the postponement cost, ORs’ opening cost, and total daily overtime costs. The concept of constraints (\ref{Eq_9}) and (\ref{Eq_10}) are similar to that of constraints (\ref{Eq_2}) and (\ref{Eq_3}) respectively. Constraint (\ref{Eq_11}) ensures that total planned surgery times are less equal to the available normal time plus overtime time. Constraint (\ref{Eq_12}) and (\ref{Eq_6}) are followed the same idea. Constraint (\ref{Eq_13}) captures the overtime cost of the planning phase, and constraint (\ref{Eq_14}) ensures the respect of the maximum amount of overtime availability. The feasibility and optimality of each solution found by the master problem should be validated by solving the subproblems for each day $d$ in the planning horizon. The subproblem of BCMIORPS is as follows:

\begin{flalign}\label{Eq_15}
 \textbf{MIORPS: }\;\; \min \beta_d &&
\end{flalign}

\noindent \mbox{Subject to:}
\begin{alignat}{2}
& \beta_d \geq \sum_{k \in K_d} c^{Ovt} o_{dk} + \sum_{i \in I_2 \cap \hat{I}_d} c^{Pos} x_i        && \label{Eq_16} \\
& c_i \geq t_{i}(1-x_i)                                                                             && \forall\, i \in \hat{I}_d \label{Eq_17}\\
& c_i \geq c_{i'} + t_{i} - M(1 - y_{ii^\prime} + x_i + x_{i'}) \qquad \qquad \qquad \qquad         && \forall\, k \in K_d, \forall\, i,i^\prime \in \hat{I}_d:i > i' \label{Eq_18}\\
& c_i' \geq c_i + t_{i'} - M(y_{ii^\prime} + x_i + x_{i'})                                          && \forall\, k \in K_d, \forall\, i,i^\prime \in \hat{I}_d:i > i' \label{Eq_19}\\
& c_i \geq c_{i'} + t_{i} - M(1 - y_{ii^\prime} + x_i + x_{i'})                                     && \forall\, i,i^\prime \in \hat{I}_d:i > i',and\, l_i = l_{i'} \label{Eq_20}\\
& c_i' \geq c_i + t_{i'} - M(y_{ii^\prime} + x_i + x_{i'})                                          && \forall\, i,i^\prime \in \hat{I}_d:i > i',and\, l_i = l_{i'} \label{Eq_21}\\
& x_i = 0                                                                                           && \forall\, i \in I_1 \cap \hat{I}_d \label{Eq_22}
\end{alignat}
\begin{alignat}{2}
& c_i \leq |T|                                                                                      && \forall\, i \in \hat{I}_d  \label{Eq_23}\\
& o_{dk} \leq  |T_2|                                                                                && \forall\, k \in K_d \label{Eq_24}\\
& x_i, y_{ii^\prime} \in \{0,1\} \qquad \qquad \qquad \qquad \qquad \qquad \qquad \qquad \qquad     && \forall\, i,i^\prime \in \hat{I}_d:i > i' \nonumber \\
& o_{dk}, c_i \geq 0                                                                                && \forall\, k \in K_d, \forall\, i,i^\prime \in \hat{I}_d:i > i' \nonumber 
\end{alignat}

The goal of the BCMIORPS subproblem is to schedule the assigned surgeries to the day $d$ with the minimum cost. Constraint (\ref{Eq_16}) captures the total cost of overtime and postponing the surgery if the subproblem can not be scheduled in the planning horizon. Constraint (\ref{Eq_17}) specifies the minimum completion time of each planned surgery on day $d$. Constraints (\ref{Eq_18}) and (\ref{Eq_19}) ensure the feasible sequence of surgeries in each OR. With the same goal, Constraints (\ref{Eq_20}) and (\ref{Eq_21}) ensure the feasible sequence of surgeries for each surgeon. Constraint (\ref{Eq_22}) implies that mandatory surgeries assigned to day $d$ are not allowed to be postponed to the next planning horizon. Constraints (\ref{Eq_23}) and (\ref{Eq_24}) enforce the upper bound for completion time and amount of overtime respectively.  

A feasible solution of the master problem will result in one of the following three different scenarios for the subproblem. We denote ${\bar{\alpha}}_d$ and ${\bar{\beta}}_d$ the calculated overtime cost by the master problem and the calculated objective function of the subproblem for day $d$ respectively.
\begin{enumerate}
  \item An optimal subproblem when ${\bar{\beta}}_d={\bar{\alpha}}_d$
  \item A suboptimal subproblem when ${\bar{\beta}}_d>{\bar{\alpha}}_d$
  \item An infeasible subproblem
\end{enumerate}
	
In the first scenario, which is the desirable scenario we expect to have for all subproblems, we do not need to add any cuts to the master problem. However, in the second and third scenarios, we must add optimality and feasibility cuts respectively to the master problem to transfer the achieved information on the proposed solution by the master problem. Equation (\ref{Eq_25}) defines the feasibility cut which is in the form of a no-good cut. When scheduling the assigned surgeries to the day $d$ is not possible, the subproblem will become infeasible. In this case, we should transfer this information to the master problem to avoid having the same combination of surgeries in the future. Equation (\ref{Eq_25}) ensures that at least one surgery should be removed from the list of assigned surgeries to day $d$, and the master problem will not produce the same infeasible combination.

\begin{flalign}\label{Eq_25}
 \sum_{k \in K_d} \sum_{i \in \hat{I}_d} (1 - x_{idk}) \geq 1 &&
\end{flalign}

When ${\bar{\beta}}_d>{\bar{\alpha}}_d$, the subproblem is feasible, but it is not optimal since the overtime cost for the day $d$ is not the same as the calculated overtime cost for the same day by the master problem. Therefore, we should add the optimality cut proposed by equation (\ref{Eq_26}) to the master problem to perform necessary modifications in the master problem.  

\begin{flalign}\label{Eq_26}
 \alpha_d \geq {\bar{\beta}}_d - {\bar{\beta}}_d \left( \sum_{k \in K_d} \sum_{i \notin \hat{I}_d} x_{idk} + \sum_{k \in K_d} \sum_{i \in \hat{I}_d} (1 - x_{idk}) \right) &&
\end{flalign}

Equation (\ref{Eq_26}) could transfer an if-then form of cut to the master problem without using an actual if-then constraint. The master problem could benefit from the optimality cut in two different ways which are as follows:  

\begin{enumerate}
  \item Producing the same results for the surgery allocation but with an updated value of $\alpha_d$
  \item Propose a new solution for allocating surgeries to days in the planning horizon
\end{enumerate}

By producing the same results for day $d$ in the master problem, the summation of $(1 - x_{idk})$ for surgeries assigned to day $d$ will be equal to zero (part B). Also, for those surgeries that are not assigned to day $d$, the first summation will be equal to zero(part A). Therefore, the optimality cut will be reduced to $\alpha_d \geq {\bar{\beta}}_d$ for day $d$. It means that, if the master problem insists on producing the same results for the day $d$ in the future, the related overtime cost of the solution should be at least equal to ${\bar{\beta}}_d$. 

\begin{flalign}
 \alpha_d \geq {\bar{\beta}}_d - {\bar{\beta}}_d \left( \overbrace{ \sum_{k \in K_d} \sum_{i \notin \hat{I}_d} x_{idk}}^\text{ Part A} + \overbrace{ \sum_{k \in K_d} \sum_{i \in \hat{I}_d} (1 - x_{idk})}^\text{ Part B} \right) && \nonumber
\end{flalign}

On the other hand, producing a new surgery allocation for day $d$ in the master problem will result in at least one difference in the values of $x_{idk}$. In case a surgery is deleted from the previous solution and postponed to the next planning horizon, one of the $x_{idk}$ will be changed from 1 to 0 and consequently, the value $(1 - x_{idk})$ will be one. Thus, part B will be equal to at least one (since we assumed one of the surgeries was removed and postponed) while part A remains unchanged. Therefore the optimality cut will be reduced to $\alpha_d \geq 0 $ which makes the cut redundant. Also, if one surgery is removed and a new one is added to the surgery plan for day $d$, both parts A and B will be equal to one, and optimality cut will be reduced to  $\alpha_d \geq -{\bar{\beta}}_d$. Therefore, in case of new surgery allocation for the day $d$, the optimality cut will be $\alpha_d \geq -n{\bar{\beta}}_d$ where $n$ could be any integer value.

\section{Hybrid RGA-CG algorithm detailed structure and descriptions}\label{AppendixD}
In this Appendix, we present a detailed description of the Hybrid RGA-CG components. 

\subsection{RGA chromosome structure}\label{chromosome_structure}
\begin{figure}[h]
    \centering
    \includegraphics[width=1\textwidth]{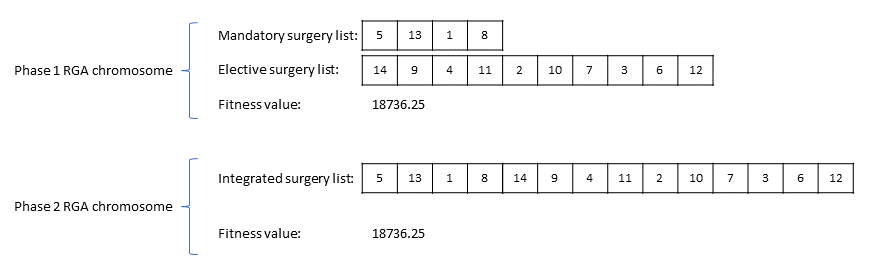}
    \caption{The structure of the RGA chromosomes for Phase1 and Phase2 respectively}
    \label{fig_4}
\end{figure}

To facilitate the coding procedure and also reduce the amount of memory needed to run the Hybrid RGA-CG, we designed different chromosome structures for each phase. In phase 1, the structure of the chromosome consists of the list of mandatory surgeries, the list of elective surgeries, and the fitness value. For phase 2, the structure has only one integrated list and the corresponding fitness value. Figure \ref{fig_4} shows the structure of chromosomes for different phases of the solution algorithm for a hypothetical instance. The first-fit strategy algorithm will calculate $x_{idt}$, $y_d$, $z_i$ matrices which are unique based on each chromosome. Therefore, it is not necessary to store the mentioned matrices alongside the mandatory and elective surgery list inside the chromosomes' structure.

\subsection{RGA crossover operators}\label{Crossover}
Crossover is a GA operator used to combine the information or chromosomes of two or more parents together and generate new children or offspring. Many crossover or permutation operators have been developed by researchers since the introduction of the GA by John Holland in the book Adaptation in Natural and Artificial Systems (1975). Each operator is designed to perform mutation or crossover procedures based on the structure of the solution. In this study, according to the definition provided in section \ref{chromosome_structure}, the chromosome's structure falls into the permutation category since we are dealing with the list of surgeries specified by unique integer numbers. Therefore, Partially Mapped Crossover (PMX) and Order crossover (OX1) are chosen based on the published book by Eiben and Smith \cite{eiben2015introduction}. In addition, we considered order-based crossover (OX2) since the OX2 and OX1 follow the same cut points, but they are different in repairing procedure \cite{back2018evolutionary}. For more information on the algorithms used to implement PMX, OX1, and OX2, please refer to \cite{eiben2015introduction} and \cite{back2018evolutionary}.

\subsection{RGA mutation operators}\label{Mutation}
The mutation is applied on a single chromosome and it has been added to GA to explore the solution area to avoid getting into the local optimum trap. Swap, Insertion, Scramble, and Inversion are four mutation operators for permutation structures introduced by Eiben and Smith \cite{eiben2015introduction}. In this study, we implemented all the mentioned four operators.

\subsection{RGA selection method}\label{selection}
Eiben and Smith \cite{eiben2015introduction} mentioned four selection methods to select the parent(s) for the crossover or mutation where each method has its own pros and cons. We chose the Fitness Proportional Selection (FPS) which is also referred to as Roulette Wheel Selection in literature. This method is a stochastic method that allows individuals with better fitness values to have a higher chance of selection. Therefore, no improvement opportunity is missed and all individuals in the population will have a chance to participate in the evolution process as a parent for crossover or individual in mutation.

\end{sloppypar}
\end{document}